\documentclass[%
amsmath,amssymb,
preprint,%
]{revtex4-1}

\usepackage{graphicx}
\usepackage{dcolumn}
\usepackage{bm}

\usepackage[utf8]{inputenc}
\usepackage[T1]{fontenc}
\usepackage{mathptmx}
\usepackage{etoolbox}
\usepackage{amssymb,amsmath,amsthm}
\usepackage{mathrsfs}
\usepackage{hyperref}
\usepackage{bm}

\newtheorem{theorem}{Theorem}[section]
\newtheorem{corollary}[theorem]{Corollary}
\newtheorem{lemma}[theorem]{Lemma}

\newtheorem{remark}[theorem]{Remark}

\usepackage{subfigure}

\makeatletter
\def\@email#1#2{%
 \endgroup
 \patchcmd{\titleblock@produce}
  {\frontmatter@RRAPformat}
  {\frontmatter@RRAPformat{\produce@RRAP{*#1\href{mailto:#2}{#2}}}\frontmatter@RRAPformat}
  {}{}
}%
\makeatother
\begin{document}

\preprint{ }

\title[Equivariant Hopf Bifurcation on a Circular Domain]{Equivariant Hopf Bifurcation in a Class of Partial Functional Differential Equations on a Circular Domain}

\author{Yaqi Chen}
\affiliation{Department of Mathematics, Harbin Institute of Technology, Weihai, Shandong 264209, P.R.China.}

\author{Xianyi Zeng}
\affiliation{Department of Mathematics, Lehigh University, Bethlehem, PA 18015, United States.}

\author{Ben Niu~*}%
 \email{niu@hit.edu.cn.}
\affiliation{Department of Mathematics, Harbin Institute of Technology, Weihai, Shandong 264209, P.R.China. 
}%

\date{\today}

\begin{abstract}
Circular domains frequently appear in the fields of ecology, biology and chemistry.
In this paper, we investigate the equivariant Hopf bifurcation of partial functional differential equations with Neumann boundary condition on a two-dimensional disk.
The properties of these bifurcations around equilibriums are analyzed rigorously by studying the equivariant normal forms.
Two reaction-diffusion systems with discrete time delays are selected as numerical examples to verify the theoretical results, in which spatially inhomogeneous periodic solutions including standing waves and rotating waves, and spatially homogeneous periodic solutions are found near the bifurcation points.

\end{abstract}

\maketitle

\section{\label{sec1}Introduction}
The research on the reaction-diffusion equation plays an important role in physics, chemistry, medicine, biology, and ecology.
Many mathematical problems, such as the existence, boundedness, regularity and stability of solutions and traveling waves have been raised in~\citep{Wang1996J,Wang2011J,Nefedov2016J,Hsu2018J,Jin2020J}.
Recently, the effect of time delay has drawn a lot of attention and Hopf bifurcation analysis becomes an effective tool to explain complex phenomena in reaction-diffusion systems, and even in more general partial functional differential equations (PFDEs).
Wu~\citep{Wu1996M}  gave a general Hopf bifurcation theorem for PFDEs by restricting the system to an eigenspace of the Laplacian.
Faria \citep{Faria2000J} gave a framework for directly calculating the normal form of PFDEs with parameters.
Based on these theories, many achievements have been made in the study of local Hopf bifurcation~\citep{Yi2009J,Hu2011J,Chen2013J,Guo2015J,Yang2022J,Song2021J,Song2021J2} or other codimension-two bifurcations~\citep{Cao2018J,Du2020J,Jiang2020J,Geng2022J}.

The phenomena of symmetry appears a lot in real-word models, which usually leads to multiple eigenvalues,
and the standard Hopf bifurcation theory of functional differential equations cannot be applied to solve such problems.
Golubitsky et al.~\citep{Golubitsky1989M} used group theory to characterize transitions in symmetric systems and worked out the bifurcation theory for a number of symmetry groups.
Based on these theories, there have been many subsequent studies on symmetry. Firstly, some researchers were concerned about nonlinear optical systems, which can effectively characterize optical problems such as circular diffraction~\citep{Razgulin2013J,Romanenko2014J,Budzinskiy2017J}. Besides, a Hopfield-Cohen-Grossberg network consisting of $n$ identical elements also has a certain symmetry, which has been studied in~\citep{Wu1998J,Guo2013M,Campbell2005J}. Furthermore, there have been many studies on the regions with $O(2)$ symmetry where the models are established. For example, Gils and Mallet-Paret\citep{Gils1986J} considered Hopf bifurcation in the presence of $O(2)$ symmetry and distinguished the phase portraits of the normal form into six cases.
In \citep{Schley2003J}, Schley studied a delay parabolic equation in a disk with the Neumann boundary conditions and proved the existence of rotating waves with methods of eigenfunction.

In recent years, a sequence of results about equivariant Hopf bifurcation in neutral functional differential equations~\citep{Guo2008J,Guo2010J} and functional differential equations of mixed type \citep{Guo2011J} have been established.
In particular, Guo~\citep{Guo2022J} applied the equivariant Hopf bifurcation theorem to study the Hopf bifurcation of a delayed Ginzburg-Landau equation on a two-dimensional disk with the homogeneous Dirichlet boundary condition.
More recently, Qu and Guo applied Lyapunov-Schmidt reduction to study the existence of inhomogeneous steady-state solutions on a unit disk~\citep{Qu2023J}, whereas different kinds of spatial-temporal solutions with symmetry have been detected by investigating isotropy subgroups of these equations~\citep{Golubitsky1989M,Wu1998J,Wu1999J,Guo2003J}.

In fact, the disk, a typical region with $O(2)$ symmetry, is usually used to describe many real-world problems.
 For example in physics, rotating waves were often observed in the case of a circular aperture~\citep{Akhmanov1992J,Ramazza1996J,Residori2007J}.
In chemical experiments, one usually studies chemical reactions in circular petri dishes, whose size may affect the existence and pattern of spiral waves~\citep{Dai2021J, Joseph1994J}.
And in the field of ecology, some lakes could be abstracted as circular domains to study the interaction between predator and prey,
and the mathematical modeling of predator-prey systems on the circular domain has been summarized in~\citep{Abid2015J,Yafia2016C}.
We found that a complete derivation of normal forms and bifurcation analysis in general partial functional differential equations on two-dimensional circular domains remains lacking.
Therefore, in this paper, we aim to consider general partial functional differential equations with homogeneous Neumann boundary condition defined on a disk and to improve the center manifold reduction technique established in \citep{Faria2000J,Faria1995J,Faria1995J2} to the normal form derivation for PFDEs on circular domains to fill the gap.

Compared to the results in \citep{Wu1996M}, due to the $O(2)$ symmetry leading to multiple pure imaginary eigenvalues, the eigenspace of the Laplacian  is sometimes two-dimensional, which gives rise to higher dimensional center subspace of the equilibrium at the bifurcation point.
By introducing similar operators as in~\citep{Faria2000J}, we derive the normal form of the equivariant Hopf bifurcation of general partial functional differential equations on a disk in explicit formulas, which can be directly applied to some models with practical significance or the model on other kind of circular domains, for example an annulus or a circular sector.
With the aid of normal forms, we find standing wave solutions and rotating wave solutions in a delayed predator-prey model. This is done after all the coefficients in the normal forms are explicitly computed.

The structure of the article is as follows. In Section \ref{sec2}, the eigenvalue problem of the Laplace operator on a circular domain is reviewed and the existence of Hopf bifurcation is explored. In Section \ref{sec3}, we study the properties of equivariant Hopf bifurcation on the center manifold. The normal forms are also rigorously derived in this section. In Section \ref{sec4}, two types of reaction-diffusion equations with discrete time delay are selected and numerically solved to verify the theoretical results.
\section{Preliminaries}\label{sec2}
\subsection{The eigenvalue problem of the Laplace operator on a circular domain}\label{sec2.1}

The eigenvalue problem associated with Laplace operators on a circular domain could be given in a standard way, see \citep{Murray2001M,Pinchover2005M}. For the convenience of our research, we use a similar method to treat the eigenvalue problem and state the main results here and consider a disk as follows
$$
\mathbb{D}=\{(r, \theta): 0 \leq r \leq R, 0 \leq \theta \leq 2 \pi\}.
$$
The Laplace operator defined in the cartesian coordinates is $\Delta \varphi=\frac{\partial^{2}}{\partial x^{2}} \varphi+\frac{\partial^{2}}{\partial y^{2}} \varphi$. Letting $x=r\cos(\theta),y=r\sin(\theta)$, it can be converted into the polar coordinates as $\Delta_{r \theta} \varphi=\frac{\partial^{2}}{\partial r^{2}} \varphi+\frac{1}{r} \cdot \frac{\partial}{\partial r} \varphi+\frac{1}{r^{2}} \cdot \frac{\partial^{2}}{\partial \theta^{2}} \varphi$.

One needs to consider the following eigenvalue problem and calculate the eigenvectors on the disk.
\begin{equation}\label{eigenvalue problem on the circular domain}
\left\{\begin{array}{l}
\Delta_{r \theta}\phi=-\lambda \phi, \\
\phi_r^{\prime}(R,\theta)=0, \theta \in[0,2 \pi].
\end{array}\right.
\end{equation}
Using the method of separation of variables and letting $\phi(r, \theta)=P(r) \Phi(\theta)$, we get that the eigenfunction corresponding to $\lambda_{mn}$ is
\begin{equation}\label{eigenfunction of the Laplace operator}
\phi_{n m}(r, \theta)=J_{n}\left(\sqrt{\lambda_{n m}} r\right)\Phi_{n}(\theta),
\end{equation}
with
\begin{equation}\label{Phi theta}
\Phi_{n}(\theta)=\begin{array}{c}
a_{n} \cos n \theta+ b_{n} \sin n \theta,
\end{array}
\end{equation}
and
\begin{equation}\label{Jn}
J_{n}(\rho)=\sum_{m=0}^{+\infty} \frac{(-1)^{m}}{m ! \Gamma(n+m+1)}\left(\frac{\rho}{2}\right)^{n+2 m}.
\end{equation}
$\lambda_{nm}$ is chosen such that the boundary condition $P'(R)=0$ is satisfied.

\begin{remark}\label{lambdanm}
Considering the Neumann boundary conditions, we have $J_n^{\prime}\left(\sqrt{\lambda_{n m}} R\right)=0$, which indicates that $\sqrt{\lambda_{n m}} R$ are roots of $J_n^{\prime}(\sqrt{\lambda}r)$. We use $\alpha_{n m}$ to represent these non-zero roots and assume that they are indexed in increasing order, i.e. $J_{n}^{\prime}\left(\alpha_{n m}\right)=0, \alpha_{n 1}<\alpha_{n 2}<\alpha_{n 3}<\cdots$, where $n \ge 0$ is the indices of the Bessel function and $m \ge 1$ are the indices for these roots. So $\lambda_{n m}=\left(\alpha_{n m} / R\right)^{2}$. For convience, we use $\alpha_{0 0}=0, \lambda_{0 0}=0$.
\end{remark}

\begin{remark}\label{basis}
From the standard Sturm-Liouville theorem, we know that for any given nonnegative integer $n$, $J_{n}\left(\alpha_{n m} r\right)$ are the orthogonal sets with weight $r$ on the interval $[0,R]$, where $J_{n}^{\prime}\left(\alpha_{n m}\right)=0$. That is, for any given $m, k$, we have
$$
\int_{0}^{R} r J_{n}\left(\frac{\alpha_{n m}}{R} r\right) J_{n}\left(\frac{\alpha_{n k}}{R} r\right) \rm{d} \it{r}=\left\{\begin{array}{cc}
0, & m \neq k,\\
\frac{R^2}{2}\left[1-\left(\frac{n}{\alpha_{nm}}\right)^2\right]J_n^2\left(\alpha_{n m}\right), & m=k.
\end{array} \right.
$$
Furthermore, for any given nonnegative integer $n$,  the $\mathbb{L}^2$ norm with weight $r$ of
function systems $\left\{J_{n}\left(\frac{\alpha_{n m}}{R} r\right)\right\}$ that include $\alpha_{0 0}$ are complete in the space $\mathbb{L}^{2}[0,R]$. Besides, the trigonometric function systems are orthogonal in the interval $[0,2 \pi]$ and complete in the space  $\mathbb{L}^{2}[0,2 \pi]$.
Therefore, for $n=0, 1, 2, \cdots$, $m=1, 2, \cdots$, we use a complexification of the space  and the system of functions
$$
\phi_0,~\phi_{nm}^{c},~\phi_{nm}^{s},
$$
constitutes an orthogonal basis with weight $r$ in the space $\mathbb{L}^{2}\{0 \leq \theta \leq 2 \pi, 0 \leq r \leq R\}$ with
$$
\phi_0=J_{0}\left(\frac{\alpha_{0 0}}{R} r\right)=1,~\phi_{nm}^{c}=J_{n}\left(\frac{\alpha_{n m}}{R} r\right) \mathrm{e}^{\mathrm{i} n \theta},~\phi_{nm}^{s}=\overline{\phi_{nm}^{c}}=J_{n}\left(\frac{\alpha_{n m}}{R} r\right)\mathrm{e}^{-\mathrm{i} n \theta}.
$$
\end{remark}

From the above analysis, we can draw the following conclusions.

\begin{theorem}\label{F-B}
 The solution of the Laplace equation on $\mathbb{D}$ with homogeneous Neumann condition at $r=R$ can be written as
$$
\varphi(r, \theta)=\phi_{0}(r, \theta)+A_{nm}\sum_{n=0}^{+\infty} \sum_{m=1}^{+\infty}\phi_{nm}^c(r, \theta)+B_{nm}\sum_{n=1}^{+\infty} \sum_{m=1}^{+\infty}\phi_{nm}^s(r, \theta),
$$
where\\
$$
A_{nm}=\frac{\delta_{n}}{ R^{2}\pi \left[ 1-\left( \frac{n}{\alpha_{nm}} \right)^2 \right] J_{n}^{2}\left(\alpha_{nm}\right)} \int_{0}^{R} \int_{0}^{2 \pi} r \varphi(r, \theta) J_{n}\left(\frac{\alpha_{nm}}{R} r\right) \mathrm{e}^{-\mathrm{i} n \theta} \mathrm{d} r \mathrm{d} \theta,
$$
$$
B_{nm}=\frac{2}{R^{2} \pi \left[ 1-\left( \frac{n}{\alpha_{nm}} \right)^2 \right] J_{n}^{2}\left(\alpha_{nm}\right)} \int_{0}^{R} \int_{0}^{2 \pi} r \varphi(r, \theta) J_{n}\left(\frac{\alpha_{nm}}{R} r\right) \mathrm{e}^{\mathrm{i} n \theta} \mathrm{d} r \mathrm{d} \theta,
$$
$$
\delta_{n}=\left\{\begin{array}{l}1, n=0, \\ 2, n \neq 0.\end{array}\right.
$$
This means, for $n=0$, the eigenspace corresponding to the eigenvalue $\lambda_{0m}$ is spanned by $\phi_{0m}^{c},~m=0,1,\cdots$. For $n> 0$, the eigenspace corresponding to the eigenvalue $\lambda_{nm}$ is spanned by $\phi_{nm}^{c}$ and $\phi_{nm}^{s},~m=1,2,\cdots$.
\end{theorem}

\begin{remark}
The above eigenvalue problems can be directly applied to an annulars domain
$$
\tilde{\mathbb{D}}=\{(r, \theta):~R_1 \leq r \leq R_2,~0 \leq \theta \leq 2 \pi\}.
$$
The difference is that the homogeneous Neumann condition is given at $r = R_1$ and $r = R_2$. Besides, $P(r)$ becomes
\begin{equation*}
  P_n(\sigma_{nm},r)=A_nJ_n(\frac{\sigma_{nm}}{R_2} r)+B_nN_n(\frac{\sigma_{nm}}{R_2} r),
\end{equation*}
with $P_n^{\prime}(\sigma_{nm},R_1)=P_n^{\prime}(\sigma_{nm},R_2)=0$,
where $N_{n}(\rho)=\frac{J_{n}(\rho) \cos n \pi-J_{-n}(\rho)}{\sin n \pi}$.
Then $\frac{P_n(\sigma_{nm},r)}{\|P_n(\sigma_{nm},r)\|_{2,2}}$ forms a standard orthogonal basis of the space $\mathbb{L}^{2}\{ R_1 \leq r \leq R_2\}$. In what follows, there will be no significant difference in the subsequent process of bifurcation analysis.

In fact, the above eigenvalue analysis in a circular sector domain
$$
\hat{\mathbb{D}}=\{(r, \theta):~0 \leq r \leq R,~0 \leq \theta \leq \Theta,~\Theta<\pi\},
$$
is also similar. The difference is that the homogeneous Neumann conditions make    $\Phi(\theta)$ become $\Phi_n(\theta)=\cos \frac{n\pi}{\Theta}\theta$.
The subsequent calculation process is even simpler, which is a direct extension of the case in one-dimensional intervals.
\end{remark}

\subsection{The existence of the Hopf bifurcation}\label{sec2.2}
We consider a general partial functional differential equations with homogeneous Neumann boundary conditions defined on a disk as follows:
\begin{equation}\label{PFDEs}
\frac{\partial U(t, x, y)}{\partial t}=D(\nu) \Delta U(t, x, y)+L(\nu)U_t(x, y)+F\left(U_t(x, y),\nu\right), \\
\end{equation}
where $t \in[0,+\infty),~\Omega=\left\{(x, y) \in \mathbb{R}^{2} \mid x^{2}+y^{2}<R^{2}\right\}$,
$$
U(t, x, y)=\left(\begin{array}{c}
 u_1(t, x, y)\\
 u_2(t, x, y)\\
 \vdots\\
 u_n(t, x, y)
\end{array}\right),
~ U_t(\vartheta)(x, y)=\left(\begin{array}{c}
 u_t^1(\vartheta)(x, y)\\
 u_t^2(\vartheta)(x, y)\\
 \vdots\\
 u_t^n(\vartheta)(x, y)
\end{array}\right),
~ D(\nu)=\left(\begin{array}{cccc}
 d_1(\nu) & 0 & \cdots & 0\\
 0 & d_2(\nu) & \cdots & 0\\
 \vdots & \vdots & \ddots & \vdots\\
 0 & 0 & \cdots & d_n(\nu)\\
\end{array}\right),
$$
$d_i(\nu)>0,i=1,2,\cdots,n,~\nu \in \mathbb{R},~U_t(\vartheta)(x, y)=U(t+\vartheta, x, y),~\vartheta \in \left[-\tau,0\right],~U_t(x, y)\in \tilde{\mathscr{C}}:=C([-\tau,0], \tilde{\mathscr{X}_{\mathbb{C}}}),~L:\mathbb{R}\times \tilde{\mathscr{C}}\rightarrow\mathscr{X}_{\mathbb{C}}$ is a bounded linear operator, and
${F}: \tilde{\mathscr{C}} \times \mathbb{R} \rightarrow \mathscr{X}_{\tilde{\mathscr{C}}}$ is a $C^k~(k\ge 3)$ function such that ${F}\left(0,  \nu\right)=0,~D_{\varphi}{F}\left(0, \nu\right)=0$ that stands for the Fr\'{e}chet derivative of ${F}\left(\varphi, \nu\right)$ with  respect to $\varphi$ at $\varphi=0$.

\begin{displaymath}
  \tilde{\mathscr{X}_{\mathbb{C}}}=\left\{\tilde{U}(x, y)\in {W}^{2,2}({\Omega}): \nabla\tilde{U}(x, y)\cdot \eta=0,(x, y)\in \partial \Omega \right\},
\end{displaymath}
where $\eta$ is the out-of-unit normal vector. Here we use the complexification space $\tilde{\mathscr{X}_{\mathbb{C}}}$, because  a complex form of eigenvector is more suitable to shorten the expressions in the normal form derivation.

Let us now explore the existing conditions of Hopf bifurcation based on system (\ref{PFDEs}). Again, we use $x=r \cos \theta$, $y=r \sin \theta$ and the domain $\Omega$ is transformed into $\mathbb{D}=\{(r, \theta): 0 \le r<R, 0 \le \theta<2 \pi\}$. For simplicity, we still use symbols in (\ref{PFDEs}).
System (\ref{PFDEs}) can be converted in polar coordinates to
\begin{equation}\label{PFDEs r theta}
\frac{\partial U(t, r, \theta)}{\partial t}=D(\nu) \Delta_{r \theta} U(t, r, \theta)+L(\nu)U_t(r, \theta)+F\left(U_t(r, \theta),\nu\right), \\
\end{equation}
and analogously, define the phase space
$$
{\mathscr{C}}:=C([-\tau,0], {\mathscr{X}_{\mathbb{C}}}),
$$
where
$$
{\mathscr{X}_{\mathbb{C}}}=\left\{\tilde{U}(r, \theta)\in {W}^{2,2}(\mathbb{D}): \partial_r \tilde{U}(R, \theta)=0,~\theta \in [0,2\pi) \right\},
$$
with inner product $\langle u(r,\theta),v(r,\theta)\rangle=\iint_{\mathbb{D}}r u(r,\theta) \bar{v}(r,\theta) \mathrm{d} r\mathrm{d}\theta$ weighted $r$ for $u(r,\theta),~v(r,\theta) \in \mathscr{X}_{\mathbb{C}}$.
Then, $U_t(r, \theta)\in {\mathscr{C}}$ .

Linearizing system (\ref{PFDEs r theta}) at the origin, we have
\begin{equation}\label{linearized system}
\frac{\partial U(t,r,\theta)}{\partial t}={D}(\nu)\Delta_{r \theta} U(t,r,\theta)+L(\nu) U_{t}(r,\theta).
\end{equation}
The characteristic equations of (\ref{linearized system}) are
\begin{equation}\label{characteristic equation}
\gamma \varphi -D(\nu) \Delta_{r \theta} \varphi-L(\nu)(\mathrm{e}^{\gamma\cdot}\varphi)=0,
\end{equation}
where $\mathrm{e}^{\gamma\cdot}(\vartheta)\varphi=\mathrm{e}^{\gamma\vartheta}\varphi$, for $\vartheta \in [-\tau,0]$, and $\gamma$ is an eigenvalue of equation (\ref{linearized system}).
By Theorem \ref{F-B}, we find that solving (\ref{characteristic equation}) is equivalent to solving the following two groups of characteristic equations. The first group is given as \begin{equation}\label{characteristic equation1}
 \mathrm{det}\left[\gamma I +\lambda_{0m}D(\nu)-L(\nu)(\mathrm{e}^{\gamma\cdot}I)\right]=0,~m=0,1,2,\cdots.
\end{equation}
The second groups of equations have multiple roots as the eigenspace of $\lambda_{nm},~n=1,2,\cdots,~m=1,2,\cdots$ is of two-dimensional. They are
\begin{equation}\label{characteristic equation2}
 \mathrm{det}\left[\gamma I +\lambda_{nm}D(\nu)-L(\nu)(\mathrm{e}^{\gamma\cdot}I)\right]^2=0,~n=1,2,\cdots,~m=1,2,\cdots.
\end{equation}

In order to consider the Hopf bifurcation, we assume that the following conditions hold for some $\nu_{\hat{\lambda}},~\hat{\lambda}=\lambda_{0m}$ or $\lambda_{nm}$.
\begin{itemize}
\item [$\bm{\mathrm{(H_1)}}$]
There exists a neighborhood $\mathscr{U}_1$ of $\nu_{\hat{\lambda}},~\hat{\lambda}=\lambda_{0m}$ such that for $\nu \in \mathscr{U}_1$, system (\ref{linearized system}) has a pair of complex simple conjugate eigenvalues $\alpha_{\hat{\lambda}}(\nu)\pm \mathrm{i}\omega_{\hat{\lambda}}(\nu)$ and the remaining eigenvalues of (\ref{linearized system}) have non-zero real part for $\nu \in \mathscr{U}_1$.
\item [$\bm{\mathrm{(H_2)}}$]
There exists a neighborhood $\mathscr{U}_2$ of $\nu_{\hat{\lambda}},~\hat{\lambda}=\lambda_{nm}$ such that for $\nu \in \mathscr{U}_2$, system (\ref{linearized system}) has a pairs of complex repeated conjugate eigenvalues $\alpha_{\hat{\lambda}}(\nu)\pm \mathrm{i}\omega_{\hat{\lambda}}(\nu)$ and the remaining eigenvalues of (\ref{linearized system}) have non-zero real part for $\nu \in \mathscr{U}_2$.
\item [$\bm{\mathrm{(H_3)}}$]
$\alpha_{\hat{\lambda}}(\nu)\pm \mathrm{i}\omega_{\hat{\lambda}}(\nu)$ are continuously differential in $\nu$ with $\alpha_{\hat{\lambda}}(\nu_{\hat{\lambda}})=0,~\omega_{\hat{\lambda}}(\nu_{\hat{\lambda}})=\omega_{\hat{\lambda}}>0$.
\end{itemize}

\begin{remark}\label{equicariant}
According to \citep{Golubitsky1989M}, problem (\ref{PFDEs r theta}) is $\Gamma$ equivariant in spatial dimension, with $\Gamma = O(2)$. For instance, write the right hand of system (\ref{PFDEs r theta}) as $\mathscr{F}(U(t,r,\theta))$, we have
$$
\mathscr{F}(\kappa U(t,r,\theta))=\kappa \mathscr{F}(U(t,r,\theta)), \forall \kappa \in \Gamma.
$$
Thus, in what follows, any solution after the action of this group is still a solution of the equation.
\end{remark}

By \citep{Golubitsky1989M,Guo2013M,Faria1995J,Ruan2003J} and Remark \ref{equicariant}, if $\bm{\mathrm{(H_1)}}$ and $\bm{\mathrm{(H_3)}}$ or $\bm{\mathrm{(H_2)}}$ and $\bm{\mathrm{(H_3)}}$ hold, noting $\hat{\nu}=\min\left\{\nu_{\hat{\lambda}}\right\}$, we know Hopf bifurcations occur at the critical values $\nu=\hat{\nu}$. When $\hat{\lambda}=\lambda_{0m},~m=0,1,2,\cdots$, the center subspace of the equilibrium is of two dimensional, so we call this a standard Hopf bifurcation. When $\hat{\lambda}=\lambda_{nm},~n=1,2,\cdots,~m=1,2,\cdots$, the center subspace of the equilibrium is of four dimensional, we say this is a (real) equivariant Hopf bifurcation. In the coming section, we will calculate the equivariant Hopf bifurcation around $E^{*}$.
\section{Hopf bifurcation analysis}\label{sec3}
\subsection{Normal form for PFDEs}
In this section, we will investigate the properties of the equivariant Hopf bifurcation around $E^{*}$, using the theory in  \citep{Wu1996M,Faria2000J,Gils1986J,Wu1999J,Faria1995J}.

Letting $\nu=\hat{\nu}+\mu$, where $\hat{\nu}$ is given in subsection \ref{sec2.2} and $\mu \in \mathbb{R}$, following the method proposed in \citep{Faria2000J}, and using $\mu$ as a new variable, the Taylor expansions of ${L}(\hat{\nu}+\mu)$ and ${D}(\hat{\nu}+\mu)$ are as follows
$$
{L}(\hat{\nu}+\mu)=\tilde{L}_0+\mu \tilde{L}_1+\frac{1}{2}\mu^2 \tilde{L}_2+\cdots,
$$
$$
{D}(\hat{\nu}+\mu)=\tilde{D}_0+\mu \tilde{D}_1+\frac{1}{2}\mu^2 \tilde{D}_2+\cdots,
$$
where $\tilde{D}_0=D(\hat{\nu}),~\tilde{L}_0(\cdot)={L}(\hat{\nu})(\cdot)$ is a linear operator from $\mathscr{C}$ to $\mathscr{X}_{\mathbb{C}}$.
Now, in the space $\mathscr{C}$, system (\ref{PFDEs r theta}) is equivalent to
\begin{equation}\label{abstract functional differential equation}
\frac{\mathrm{d} U(t)}{\mathrm{d} t}=\tilde{D}_0\Delta_{r \theta} U(t)+\tilde{L}_0 U_{t}+\tilde{F}\left(U_{t}, \mu\right),
\end{equation}
where
$\tilde{F}(\varphi, \mu)=[{D}(\hat{\nu}+\mu)-\tilde{D}_0]\Delta_{r\theta}\varphi(0)+[L(\hat{\nu}+\mu)-\tilde{L}_0](\varphi)+F(\varphi,\hat{\nu}+\mu)$,
and the linearization system of (\ref{abstract functional differential equation}) is obtained as follows
\begin{equation}\label{Linearizing abstract functional differential equation}
\frac{\mathrm{d} U(t)}{\mathrm{d} t}=\tilde{D}_0 \Delta_{r \theta} U(t)+\tilde{L}_0 U_{t}.
\end{equation}

\subsubsection{Decomposition of $\mathscr{C}$}
Let $A:{\mathscr{C}} \rightarrow \mathscr{X}_{\mathbb{C}}$ represent the infinitesimal generators of the semigroup induced by the solutions of (\ref{Linearizing abstract functional differential equation}), and $A^{*}$ is the adjoint operator of $A$,
which satisfy
\begin{equation}\label{A}
 A \varphi(\vartheta)=\left\{\begin{array}{cc}\varphi^{\prime}(\vartheta), & \vartheta \in[-\tau,0), \\
 \tilde{D}_0 \Delta_{r \theta} \varphi(0)+\tilde{L}_0 \varphi, & \vartheta=0,\end{array}\right.
\end{equation}
\begin{equation}\label{Astar}
A^{*} \psi(\varrho)=\left\{\begin{array}{cc}-\psi^{\prime}(\varrho), & \varrho \in(0,\tau], \\
-\tilde{D}_0 \Delta_{r \theta} \psi(0)-\tilde{L}_0 \psi, & \varrho=0.\end{array}\right.
\end{equation}
In addition, define a bilinear pairing
\begin{equation}\label{bilinear product}
\begin{aligned}
(\psi, \varphi) &=\langle \varphi(0),\psi(0)\rangle-\int_{-\tau}^{0} \langle  \varphi(\xi),\tilde{L}_0{\psi}(\xi+\tau)\rangle  \mathrm{d} \xi. \\
\end{aligned}
\end{equation}
From the discussion in subsection \ref{sec2.2}, we know that $A$ has a pair of repeated purely imaginary eigenvalues $\pm \mathrm{i} \omega_{\hat{\lambda}} $ which are also eigenvalues of $A^{*}$.
Let the central subspace $P$ and $P^{*}$ be the generalized eigenspace of $A$ and $A^{*}$ about $\Lambda_{0}=\{\pm \mathrm{i} \omega_{\hat{\lambda}} ,\pm \mathrm{i} \omega_{\hat{\lambda}} \}$, respectively. $P^{*}$ is the adjoint space of $P$.
An important task is to decompose the space $\mathscr{C}$ through the relationship of bases in $P$ and $P^{*}$, and we write $\mathscr{C}=P_{C N} \oplus Q_{S}$, where $P_{C N}$ is the central subspace and $Q_{S}$ is its complementary space.

Define
$$
\hat{\phi}_{nm}^c=\frac{\phi_{nm}^c}{\|\phi_{nm}^c\|_{2,2}},~
\hat{\phi}_{nm}^s=\frac{\phi_{nm}^s}{\|\phi_{nm}^s\|_{2,2}}.
$$

\begin{lemma}
Let the basis of $P$ is
\begin{equation}\label{basis of P}
\Phi_{r \theta}(\vartheta)=\left(\Phi_{r \theta}^1(\vartheta),\Phi_{r \theta}^2(\vartheta)\right)=\left(\Phi^1(\vartheta)\cdot \hat{\phi}_{nm}^c,\Phi^2(\vartheta)\cdot\hat{\phi}_{nm}^s\right),~\vartheta \in [-\tau,0],
\end{equation}
with
$$
\begin{aligned}
\Phi_{r \theta}^1(\vartheta)=\left(\Phi_1(\vartheta)\cdot \hat{\phi}_{nm}^c,\Phi_2(\vartheta)\cdot \hat{\phi}_{nm}^c\right)=\left(\mathrm{e}^{\mathrm{i}\omega_{\hat{\lambda}} \vartheta} \xi \hat{\phi}_{nm}^c, \mathrm{e}^{-\mathrm{i}\omega_{\hat{\lambda}} \vartheta} \bar{\xi}\hat{\phi}_{nm}^c\right),\\
\Phi_{r \theta}^2(\vartheta)=\left(\Phi_3(\vartheta)\cdot \hat{\phi}_{nm}^s,\Phi_4(\vartheta)\cdot \hat{\phi}_{nm}^s\right)=\left(\mathrm{e}^{\mathrm{i}\omega_{\hat{\lambda}} \vartheta} {\xi} \hat{\phi}_{nm}^s, \mathrm{e}^{-\mathrm{i}\omega_{\hat{\lambda}} \vartheta} \bar{\xi} \hat{\phi}_{nm}^s\right),
\end{aligned}
$$
where $\xi$ can be noted as $\xi=(p_{11},p_{12},\cdots,p_{1n})^{\mathrm{T}}$. A basis for the adjoint space $P^*$ is
\begin{equation}\label{basis of P*}
\Psi_{r \theta}(\varrho)=\left(\Psi_{r \theta}^1(\varrho),\Psi_{r \theta}^2(\varrho)\right)^{\mathrm{T}}=\left(\Psi^1(\varrho)\cdot\hat{\phi}_{nm}^c,\Psi^2(\varrho)\cdot\hat{\phi}_{nm}^s\right)^{\mathrm{T}},~\varrho \in [0,\tau],
\end{equation}
with
$$
\begin{aligned}
\Psi_{r \theta}^1(\varrho)=\left(\Psi_1(\varrho)\cdot \hat{\phi}_{nm}^c,\Psi_2(\varrho)\cdot \hat{\phi}_{nm}^c\right)^{\mathrm{T}}=\left({q}^{-1} \mathrm{e}^{\mathrm{i}\omega_{\hat{\lambda}} \varrho} {{\xi}}^{\mathrm{T}}\hat{\phi}_{nm}^c, \bar{q}^{-1} \mathrm{e}^{-\mathrm{i}\omega_{\hat{\lambda}} \varrho} \bar{\xi}^{\rm{T}} \hat{\phi}_{nm}^c\right)^{\rm{T}},\\
\Psi_{r \theta}^2(\varrho)=\left(\Psi_3(\varrho)\cdot \hat{\phi}_{nm}^s,\Psi_4(\varrho)\cdot \hat{\phi}_{nm}^s\right)^{\mathrm{T}}=\left({{q}}^{-1} \mathrm{e}^{\mathrm{i}\omega_{\hat{\lambda}} \varrho} {\xi}^{\mathrm{T}} \hat{\phi}_{nm}^s, \bar{q}^{-1} \mathrm{e}^{-\mathrm{i}\omega_{\hat{\lambda}} \varrho} \bar{{\xi}}^{\rm{T}} \hat{\phi}_{nm}^s\right)^{\rm{T}},
\end{aligned}
$$
where $q$ can be obtained by $(\Psi_{r\theta},\Phi_{r\theta})=I$, according to the adjoint bilinear form defined in (\ref{bilinear product}).
\end{lemma}

One can decompose $U_t$ into two parts:
\begin{equation}\label{fenjie}
U_t=U_t^P+U_t^Q=\sum_{k=1}^2 \Phi_{r \theta}^k\left(\Psi_{r \theta}^k,  U_t\right)+U_t^Q=\sum_{k=1}^2 \Phi_{r \theta}^k z_{r \theta}^k+y_t,
\end{equation}
where $z_{r \theta}^k=\left(\Psi_{r \theta}^k,  U_t\right),~y_t \in Q_s$. Define $z =(z_{r\theta}^1,z_{r\theta}^2)\equiv (z_1,z_2,z_3,z_4)^{\mathrm{T}}$ as the local coordinate system on the four-dimensional center manifold, which is induced by the basis $\Phi_{r \theta}$.

Then, we get that
\begin{equation}\label{zt2}
\begin{aligned}
\dot{z}(t) &= \tilde{B} z(t)+\left(\begin{array}{cc}
\left\langle \tilde{F}(\sum_{k=1}^2\Phi_{r \theta}^k z_{r\theta}^k+y,\mu), \Psi_{r \theta}^1(0)\right\rangle\\
\left\langle \tilde{F}(\sum_{k=1}^2\Phi_{r \theta}^k z_{r\theta}^k+y,\mu), \Psi_{r \theta}^2(0)\right\rangle
\end{array}\right),\\
\frac{\mathrm{d}y}{\mathrm{d}t}&=A_Qy+(I-\pi)X_0\tilde{F}\left(\sum_{k=1}^2\Phi_{r \theta}^k z_{r\theta}^k +y,\mu\right),
\end{aligned}
\end{equation}
where $A_Q$ is the restriction of $A$ on $Q_s$, $A_Q \varphi=A \varphi$ for $\varphi \in Q_s$, $\pi: \mathscr{C} \rightarrow P_{CN}$ is the projection, and
$$
\tilde{B}=\left(\begin{array}{cccc}
\mathrm{i}\omega_{\hat{\lambda}} & 0 & 0 & 0\\
0 & -\mathrm{i}\omega_{\hat{\lambda}} & 0 & 0\\
0 & 0 & \mathrm{i}\omega_{\hat{\lambda}} & 0\\
0 & 0 & 0 & -\mathrm{i}\omega_{\hat{\lambda}}\\
\end{array}
\right).
$$

According to the formal Taylor expansion $\tilde{F}(\varphi,\mu)=\sum_{j\ge 2}\frac{1}{j!}\tilde{F}_j(\varphi,\mu)$, (\ref{zt2}) can be written as
\begin{equation}\label{zt22}
\begin{aligned}
\dot{z}(t) &= \tilde{B} z(t)+\sum_{j\ge 2}\frac{1}{j!}f_j^1(z,y,\mu),\\
\frac{\mathrm{d}y}{\mathrm{d}t}&=A_Qy+\sum_{j\ge 2}\frac{1}{j!}f_j^2(z,y,\mu),
\end{aligned}
\end{equation}
where $f_j=(f_j^1,f_j^2),~j\ge 2$ are defined by
\begin{equation}\label{fj}
\begin{aligned}
f_j^1(z,y,\mu)&=\left(\begin{array}{cc}
\left\langle \tilde{F}(\sum_{k=1}^2\Phi_{r \theta}^k z_{r\theta}^k+y,\mu), \Psi_{r \theta}^1(0)\right\rangle\\
\left\langle \tilde{F}(\sum_{k=1}^2\Phi_{r \theta}^k z_{r\theta}^k+y,\mu), \Psi_{r \theta}^2(0)\right\rangle
\end{array}\right),\\
f_j^2(z,y,\mu)&=(I-\pi)X_0 \tilde{F}_j\left(\sum_{k=1}^2\Phi_{r \theta}^k z_{r\theta}^k +y,\mu\right).
\end{aligned}
\end{equation}

Referring to \citep{Faria2000J}, we get the normal form on the center manifold of the origin is
\begin{equation}\label{normal form}
\begin{aligned}
 \dot{z}(t)&=\tilde{B}z(t)+\frac{1}{2}g_2^1(z,y,\mu)+\frac{1}{6}g_3^1(z,y,\mu)+h.o.t.,\\
 \frac{\mathrm{d}y}{\mathrm{d}t}&=A_Qy+\frac{1}{2}g_2^2(z,y,\mu)+\frac{1}{6}g_3^2(z,y,\mu)+h.o.t.,
 \end{aligned}
\end{equation}
where $g=(g_j^1,g_j^2),~j\ge 2$ is given by
$$
g_j(z,y,\mu)=\bar{f}_j(z,y,\mu)-M_jU_j(z,\mu),
$$
where $\bar{f}_j^1$ is the terms of order $j$ in $(z, y)$
obtained after the computation of normal forms up to order $j-1$, $U_j=(U_j^1,U_j^2)$ denotes the change of variables about the transformation from $f_j$ to $g_j$, and the operator $M_j=(M_j^1,M_j^2)$ is defined by
\begin{equation}\label{Mj}
\begin{aligned}
M_j^1&: \mathbb{V}_j^5(\mathscr{C}^4) \rightarrow \mathbb{V}_j^5(\mathscr{C}^4),\\
M_j^1U_j^1&=D_zU_j^1(z,\mu)\tilde{B}z-\tilde{B}U_j^1(z,\mu),\\
M_j^2&: \mathbb{V}_j^5(Q_s) \rightarrow \mathbb{V}_j^5(\mathrm{Ker}\pi),\\
M_j^2U_j^2&=D_zU_j^2(z,\mu)\tilde{B}z-A_QU_j^2(z,\mu),
\end{aligned}
\end{equation}
where $ \mathbb{V}_j^5(Y)$ denotes the space homogeneous polynomials of $z=(z_1,z_2,z_3,z_4)^{\mathrm{T}}$ and $\mu$ with coefficients in $\mathscr{C}^4$.

It is easy to verify that
\begin{equation}\label{Mj1z}
M_j^1(\mu z^p e_k)=\mathrm{i} \omega_{\hat{\lambda}} \mu\left(p_1-p_2+p_3-p_4+(-1)^k\right)z^p e_k,~ |p|=j-1,
\end{equation}
where $j \ge 2,~ k=1,2,3,4$, and $\{e_1,e_2,e_3,e_4\}$ is the canonical basis for $\mathscr{C}^4$.

\subsubsection{Calculation of $g_2^1(z,0,\mu)$}

For $j=2$, similar to the results in \citep{Budzinskiy2017J,Wu1999J}, we have
$$
\mathrm{Ker}(M_2^1)=\mathrm{span}\left\{\left(\begin{array}{cccc}
\mu z_1\\
0\\
0\\
0
\end{array}
\right),
\left(\begin{array}{cccc}
0\\
\mu z_2\\
0\\
0
\end{array}
\right),
\left(\begin{array}{cccc}
0\\
0\\
\mu z_3\\
0
\end{array}
\right),
\left(\begin{array}{cccc}
0\\
0\\
0\\
\mu z_4
\end{array}
\right),
\left(\begin{array}{cccc}
\mu z_3\\
0\\
0\\
0
\end{array}
\right),
\left(\begin{array}{cccc}
0\\
\mu z_4\\
0\\
0
\end{array}
\right),
\left(\begin{array}{cccc}
0\\
0\\
\mu z_1\\
0
\end{array}
\right),
\left(\begin{array}{cccc}
0\\
0\\
0\\
\mu z_2
\end{array}
\right)
\right\},
$$
then
$$
\begin{aligned}
   &\mathrm{Ker}(M_2^1) \cap \mathrm{span}\left\{\mu z^p e_k; |p|=1,k=1,2,3,4 \right\}\\
=&\mathrm{span}\left\{\mu x_1 e_1,\mu x_3 e_1,\mu x_2 e_2,\mu x_4 e_2,\mu x_1 e_3,\mu x_3 e_3,\mu x_2 e_4,\mu x_4 e_4\right\}.
\end{aligned}
$$
Therefore, the second order term of $\tilde{F}(U_t,\mu)$ is
\begin{equation}\label{F2}
\tilde{F}_2(U_t,\mu)= \mu \tilde{D}_1 \Delta U_t(0)+\mu \tilde{L}_1 U_t+F_2(U_t,\mu),
\end{equation}
and
\begin{equation}\label{F2z}
\begin{aligned}
\tilde{F}_2(z,y,\mu)&=\tilde{F}(\Phi_{r\theta} z +y,\mu)\\
&=\mu \tilde{D}_1 \Delta \left(\Phi_{r\theta}(0) z +y(0)\right)+\mu \tilde{L}_1 (\Phi_{r\theta} z +y)+F_2(\Phi_{r\theta} z +y,\mu).
\end{aligned}
\end{equation}
Since $F(0,\mu)=0,~DF(0,\mu)=0$, $F_2(\Phi_{r\theta} z+y,\mu)$ can be written as follows
\begin{equation}\label{F2phi}
\begin{aligned}
F_2(\Phi_{r\theta} z+y,\mu)&=F_2(\Phi_{r\theta} z+y,0)\\
&=\sum_{p_1+p_2+p_3+p_4=2} A_{p_1p_2p_3p_4} \left(\hat{\phi}_{nm}^c\right)^{p_1+p_2}  \left({\hat{\phi}_{nm}^s}\right)^{p_3+p_4} z_1^{p_1}z_2^{p_2}z_3^{p_3}z_4^{p_4}+S_2(\Phi_{r\theta} z,y)+o(|y|^2),
\end{aligned}
\end{equation}
where $S_2$ represents the linear terms of $y$, which can be calculated by $DF_2(\Phi_{r\theta} z+y,0)|_{y=0}(y)$.

By (\ref{zt22})-(\ref{F2phi}), noticing the fact
$$
\int_0^R\int_0^{2\pi} r \hat{\phi}_{nm}^c \hat{\phi}_{nm}^s \mathrm{d} \theta \mathrm{d} r=1,
$$
and the relationship of $\Phi_{r\theta}$ and $\Psi_{r\theta}$,
we obtain

\begin{equation}\label{g21}
\frac{1}{2} g_2^1(z,0,\mu)=\frac{1}{2} \mathrm{Proj}_{\mathrm{Ker}(M_2^1)}f_2^1(z,0,\mu)=\left(
\begin{array}{cccc}
B_{11}\mu z_1\\
\overline{B_{11}} \mu z_2\\
B_{11}\mu z_3\\
\overline{B_{11}} \mu z_4
\end{array}
\right),
\end{equation}
with
\begin{equation}\label{B11B13}
\begin{aligned}
&B_{11}=\frac{1}{2}\overline{\Psi_1(0)}(-\lambda_{nm} \tilde{D}_1\Phi_1(0)+\tilde{L}_1\Phi_1).\\
\end{aligned}
\end{equation}

\subsubsection{Calculation of $g_3^1(z,0,\mu)$}

For $j=3$, we have
$$
\begin{aligned}
\mathrm{Ker}(M_3^1)=\mathrm{span}\left\{\left(\begin{array}{cccc}
z_1^2z_2\\
0\\
0\\
0
\end{array}
\right),
\left(\begin{array}{cccc}
0\\
0\\
z_1^2z_2\\
0
\end{array}
\right),
\left(\begin{array}{cccc}
z_1^2z_4\\
0\\
0\\
0
\end{array}
\right),
\left(\begin{array}{cccc}
0\\
0\\
z_1^2z_4\\
0
\end{array}
\right),
\left(\begin{array}{cccc}
z_3^2z_2\\
0\\
0\\
0
\end{array}
\right),
\left(\begin{array}{cccc}
0\\
0\\
z_3^2z_2\\
0
\end{array}
\right),
\left(\begin{array}{cccc}
z_3^2z_4\\
0\\
0\\
0
\end{array}
\right),
\left(\begin{array}{cccc}
0\\
0\\
z_3^2z_4\\
0
\end{array}
\right),\right.\\
\left(\begin{array}{cccc}
0\\
z_2^2z_1\\
0\\
0
\end{array}
\right),
\left(\begin{array}{cccc}
0\\
0\\
0\\
z_2^2z_1
\end{array}
\right),
\left(\begin{array}{cccc}
0\\
z_2^2z_3\\
0\\
0
\end{array}
\right),
\left(\begin{array}{cccc}
0\\
0\\
0\\
z_2^2z_4
\end{array}
\right),
\left(\begin{array}{cccc}
0\\
z_4^2z_1\\
0\\
0
\end{array}
\right),
\left(\begin{array}{cccc}
0\\
0\\
0\\
z_4^2z_1
\end{array}
\right),
\left(\begin{array}{cccc}
0\\
z_4^2z_3\\
0\\
0
\end{array}
\right),
\left(\begin{array}{cccc}
0\\
0\\
0\\
z_4^2z_3
\end{array}
\right),\\
\left.\left(\begin{array}{cccc}
z_1z_2z_3\\
0\\
0\\
0
\end{array}
\right),
\left(\begin{array}{cccc}
0\\
0\\
z_1z_2z_3\\
0
\end{array}
\right),
\left(\begin{array}{cccc}
z_1z_3z_4\\
0\\
0\\
0
\end{array}
\right),
\left(\begin{array}{cccc}
0\\
0\\
z_1z_3z_4\\
0
\end{array}
\right),
\left(\begin{array}{cccc}
0\\
z_1z_2z_4\\
0\\
0
\end{array}
\right),
\left(\begin{array}{cccc}
0\\
0\\
0\\
z_1z_2z_4
\end{array}
\right),\left(\begin{array}{cccc}
0\\
z_2z_3z_4\\
0\\
0
\end{array}
\right),
\left(\begin{array}{cccc}
0\\
0\\
0\\
z_2z_3z_4
\end{array}
\right)
\right\},
\end{aligned}
$$
see \citep{Budzinskiy2017J,Wu1999J} again. Then
$$
\mathrm{Ker}(M_3^1) \cap \mathrm{span}\left\{\mu z^p e_k; |p|=2,k=1,2,3,4 \right\}
= \emptyset .
$$
We define
\begin{equation}\label{f31}
\bar{f}_3^1(z,0,\mu)=f_3^1(z,0,\mu)+\frac{3}{2}\left[D_z f_2^1(z,0,\mu)U_2^1(z,\mu)+D_y f_2^1(z,0,\mu)U_2^2(z,\mu)-D_zU_2^1(z,\mu)g_2^1(z,0,\mu)\right].
\end{equation}
According to \citep{Faria2000J}, the normal form up to the third order is
$$
\begin{aligned}
g_3^1(z,0,\mu)&=\mathrm{Proj}_{\mathrm{Ker}(M_3^1)}\bar{f}_3^1(z,0,\mu)\\
&=\mathrm{Proj}_{\mathrm{Ker}(M_3^1)}\bar{f}_3^1(z,0,0)+o(\mu^2|x|).
\end{aligned}
$$
Since $g_2^1(z,0,0)=0$, we only need to calcalate three parts
$$
\mathrm{Proj}_{\mathrm{Ker}(M_3^1)} f_3^1(z,0,0),
$$
$$
\mathrm{Proj}_{\mathrm{Ker}(M_3^1)} \left(D_z f_2^1(z,0,0)U_2^1(z,0)\right),
$$
and
$$
\mathrm{Proj}_{\mathrm{Ker}(M_3^1)}\left(D_y f_2^1(z,0,0)U_2^2(z,0)\right).
$$

Through calculation, we obtain the main results as follows. Please refer to  Appendix \ref{ Appendix A} for the specific calculation process.

\begin{equation}\label{Part1}
\frac{1}{3!} \mathrm{Proj}_{\mathrm{Ker}(M_3^1)}f_3^1(z,0,0)=\left(
\begin{array}{cccc}
C_{2001}z_1^2z_4+C_{1110}z_1z_2z_3\\
\overline{C_{2001}}z_2^2z_3+\overline{C_{1110}}z_1z_2z_4\\
C_{2001}z_3^2z_2+C_{1110}z_1z_3z_4\\
\overline{C_{2001}}z_4^2z_1+\overline{C_{1110}}z_2z_3z_4
\end{array}
\right),
\end{equation}
where
\begin{equation}\label{C}
\begin{aligned}
&C_{2001}=\frac{1}{6}\overline{\Psi_1(0)}A_{2001}\mathrm{M}_{22},~C_{1110}=\frac{1}{6}\overline{\Psi_1(0)}A_{1110}\mathrm{M}_{22}.
\end{aligned}
\end{equation}

\begin{equation}\label{Part2}
\begin{aligned}
\frac{1}{3!}\mathrm{Proj}_{\mathrm{Ker(M_3^1)}}\left(D_zf_2^1(z,0,0)U_2^1(z,0)\right)=\bf{0}.
\end{aligned}
\end{equation}

\begin{equation}\label{Part3}
\begin{aligned}
&\frac{1}{3!}\mathrm{Proj}_{\mathrm{Ker(M_3^1)}}\left(D_yf_2^1(z,0,0)U_2^2(z,0)\right)\\
=&\left(
\begin{array}{cccc}
E_{2100}z_1^2z_2+E_{2001}z_1^2z_4+E_{0120}z_3^2z_2+E_{0021}z_3^2z_4+E_{1110}z_1z_2z_3+E_{1011}z_1z_3z_4\\
\overline{E_{2100}}z_1z_2^2+\overline{E_{2001}}z_2^2z_3+\overline{E_{0120}}z_4^2z_1+\overline{E_{0021}}z_4^2z_3+\overline{E_{1110}}z_1z_2z_4+\overline{E_{1011}}z_2z_3z_4\\
E_{2100}z_3^2z_4+E_{2001}z_3^2z_2+E_{0120}z_1^2z_4+E_{0021}z_1^2z_2+E_{1110}z_1z_3z_4+E_{1011}z_1z_2z_3\\
\overline{E_{2100}}z_3z_4^2+\overline{E_{2001}}z_4^2z_1+\overline{E_{0120}}z_2^2z_3+\overline{E_{0021}}z_1^2z_2+\overline{E_{1110}}z_2z_3z_4+\overline{E_{1011}}z_1z_2z_4
\end{array}
\right),
\end{aligned}
\end{equation}
where
$$
\begin{aligned}
E_{2100}=\frac{1}{6}\overline{\Psi_1(0)}&\left[\mathrm{M}_{0kcs}^c\left(S_{yz_1}(h_{0k1100}^{ccs})+S_{yz_2}(h_{0k2000}^{ccs})\right)\right],\\
E_{2001}=\frac{1}{6}\overline{\Psi_1(0)}&\left[\mathrm{M}_{0kcs}^cS_{yz_1}(h_{0k1001}^{ccs})+\mathrm{M}_{2nkss}^cS_{yz_4}(h_{2nk2000}^{css})\right],\\
E_{0120}=\frac{1}{6}\overline{\Psi_1(0)}&\left[\mathrm{M}_{0kcs}^cS_{yz_2}(h_{0k0020})^{ccs})+\mathrm{M}_{2nkss}^cS_{yz_3}(h_{2nk0110}^{css})\right],\\
E_{0021}=\frac{1}{6}\overline{\Psi_1(0)}&\left[\mathrm{M}_{2nkss}^c\left(S_{yz_3}(h_{2nk0011}^{css})+S_{yz_4}(h_{2nk0020}^{css})\right)\right],\\
\end{aligned}
$$
$$
\begin{aligned}
E_{1110}=\frac{1}{6}\overline{\Psi_1(0)}&\left[\mathrm{M}_{0kcs}^c\left(S_{yz_1}(h_{0k0110}^{ccs})+S_{yz_2}(h_{0k1010}^{ccs})\right)+\mathrm{M}_{2nkss}^cS_{yz_3}(h_{2nk1100}^{css})\right],\\
E_{1011}=\frac{1}{6}\overline{\Psi_1(0)}&\left[\mathrm{M}_{0kcs}^cS_{yz_1}(h_{0k0011}^{ccs})
+\mathrm{M}_{2nkss}^c\left(S_{yz_3}(h_{2nk1001}^{css})+S_{yz_4}(h_{2nk1010}^{css})\right)\right].\\
\end{aligned}
$$

Hence, by (\ref{f31}), (\ref{Part1}), (\ref{Part2}),and (\ref{Part3}), we have

\begin{equation}\label{g31}
\begin{aligned}
\frac{1}{3!}g_3^1(z,0,0)=&\frac{1}{3!}\mathrm{Proj}_{\mathrm{Ker(M_3^1)}}\bar{f}_3^1(z,0,0)\\
=&\left(
\begin{array}{cccc}
B_{2100}z_1^2z_2+B_{2001}z_1^2z_4+B_{0120}z_3^2z_2+B_{0021}z_3^2z_4+B_{1110}z_1z_2z_3+B_{1011}z_1z_3z_4\\
\overline{B_{2100}}z_1z_2^2+\overline{B_{2001}}z_2^2z_3+\overline{B_{0120}}z_4^2z_1+\overline{B_{0021}}z_4^2z_3+\overline{B_{1110}}z_1z_2z_4+\overline{B_{1011}}z_2z_3z_4\\
B_{2100}z_3^2z_4+B_{2001}z_3^2z_2+B_{0120}z_1^2z_4+B_{0021}z_1^2z_2+B_{1110}z_1z_3z_4+B_{1011}z_1z_2z_3\\
\overline{B_{2100}}z_3z_4^2+\overline{B_{2001}}z_4^2z_1+\overline{B_{0120}}z_2^2z_3+\overline{B_{0021}}z_1^2z_2+\overline{B_{1110}}z_2z_3z_4+\overline{B_{1011}}z_1z_2z_4
\end{array}
\right),
\end{aligned}
\end{equation}
with
$$
B_{p_1p_2p_3p_4}=C_{p_1p_2p_3p_4}+\frac{3}{2}(D_{p_1p_2p_3p_4}+E_{p_1p_2p_3p_4}).
$$

\subsubsection{The normal form}
Based on the above analysis, the normal form truncated to the third order on the center manifold can be summarized as follows
\begin{equation}\label{Normal form}
\begin{aligned}
&\dot{z}_1=\mathrm{i}\omega_{\hat{\lambda}}z_1+B_{11}z_1\mu+B_{2100}z_1^2z_2+B_{2001}z_1^2z_4+B_{0120}z_3^2z_2+B_{0021}z_3^2z_4+B_{1110}z_1z_2z_3+B_{1011}z_1z_3z_4,\\
&\dot{z}_2=-\mathrm{i}\omega_{\hat{\lambda}}z_2+\overline{B_{11}}z_2\mu+\overline{B_{2100}}z_1z_2^2+\overline{B_{2001}}z_2^2z_3+\overline{B_{0120}}z_4^2z_1+\overline{B_{0021}}z_4^2z_3+\overline{B_{1110}}z_1z_2z_4+\overline{B_{1011}}z_2z_3z_4,\\
&\dot{z}_3=\mathrm{i}\omega_{\hat{\lambda}}z_3+B_{11}z_3\mu+B_{2100}z_3^2z_4+B_{2001}z_3^2z_2+B_{0120}z_1^2z_4+B_{0021}z_1^2z_2+B_{1110}z_1z_3z_4+B_{1011}z_1z_2z_3,\\
&\dot{z}_4=-\mathrm{i}\omega_{\hat{\lambda}}z_4+\overline{B_{11}}z_4\mu+\overline{B_{2100}}z_3z_4^2+\overline{B_{2001}}z_4^2z_1+\overline{B_{0120}}z_2^2z_3+\overline{B_{0021}}z_1^2z_2+\overline{B_{1110}}z_2z_3z_4+\overline{B_{1011}}z_1z_2z_4.
\end{aligned}
\end{equation}

\begin{lemma}\label{reduce}
By \citep{Gils1986J}, the normal form truncated to the third order can be reduced to
\begin{equation}\label{z1z2z3z4}
\begin{aligned}
&\dot{z}_1=\mathrm{i}\omega_{\hat{\lambda}}z_1+B_{11}z_1\mu+B_{2001}z_1^2z_4+B_{1110}z_1z_2z_3,\\
&\dot{z}_2=-\mathrm{i}\omega_{\hat{\lambda}}z_2+\overline{B_{11}}z_2\mu+\overline{B_{2001}}z_3z_2^2+\overline{B_{1110}}z_1z_2z_4,\\
&\dot{z}_3=\mathrm{i}\omega_{\hat{\lambda}}z_3+B_{11}z_3\mu+B_{2001}z_3^2z_2+B_{1110}z_1z_3z_4,\\
&\dot{z}_4=-\mathrm{i}\omega_{\hat{\lambda}}z_4+\overline{B_{11}}z_4\mu+\overline{B_{2001}}z_1z_4^2+\overline{B_{1110}}z_2z_3z_4.
\end{aligned}
\end{equation}
\end{lemma}

The proof is given in Appendix \ref{Proof of Lemma}.

Introducing double sets of polar coordinates
\begin{equation}\label{chi}
\begin{aligned}
z_1=\rho_1\mathrm{e}^{\mathrm{i}\chi_1},~z_4=\rho_1\mathrm{e}^{-\mathrm{i}\chi_1},\\
z_3=\rho_2\mathrm{e}^{\mathrm{i}\chi_2},~z_2=\rho_2\mathrm{e}^{-\mathrm{i}\chi_2},
\end{aligned}
\end{equation}
we can obtain that
\begin{equation}\label{rho}
\begin{aligned}
&\dot{\rho}_1=(a_1\mu+a_2\rho_1^2+a_3\rho_2^2)\rho_1,\\
&\dot{\chi}_1=\omega_{\hat{\lambda}},\\
&\dot{\rho}_2=(a_1\mu+a_2\rho_2^2+a_3\rho_1^2)\rho_2,\\
&\dot{\chi}_2=\omega_{\hat{\lambda}},
\end{aligned}
\end{equation}
with
$$
a_1=\mathrm{Re}\{B_{11}\},~a_2=\mathrm{Re}\{B_{2001}\},~a_3=\mathrm{Re}\{B_{1110}\}.\\
$$

Based on the above analysis, by \citep{Gils1986J,Guckenheimer1983M}, we get, when $a_1\mu<0(>0)$, system (\ref{rho}) has six unfoldings(see Table \ref{tab1}) and their dynamical classifications are shown in Table \ref{tab2}.

\begin{table}
\caption{The six unfoldings of system (\ref{rho}). }
\label{tab1}
\begin{ruledtabular}
\begin{tabular}{ccccccc}
 {Case} &  1& 2 & 3 & 4 & 5 & 6 \\
   \hline
    $a_2$ & -- & -- & -- & + & + & +  \\
    $a_2+a_3$ & -- & -- & + & -- & + & +\\
   $a_2-a_3$ & -- & + & -- & + & + & -- \\
\end{tabular}
\end{ruledtabular}
\end{table}

\begin{table}
\caption{The dynamical classifications of system (\ref{rho}) in each case.}
\begin{ruledtabular}
\label{tab2}
\begin{tabular}{ccccccc}

{}       &                                                              {Case 1} &  {Case 2}& {Case 3} & {Case 4} & {Case 5} & {Case 6}  \\ \hline
{$a_1\mu<0$}
                                                      & \begin{minipage}[b]{0.13\columnwidth}\raisebox{-.5\height}{\centerline{\includegraphics[width=0.9\linewidth]{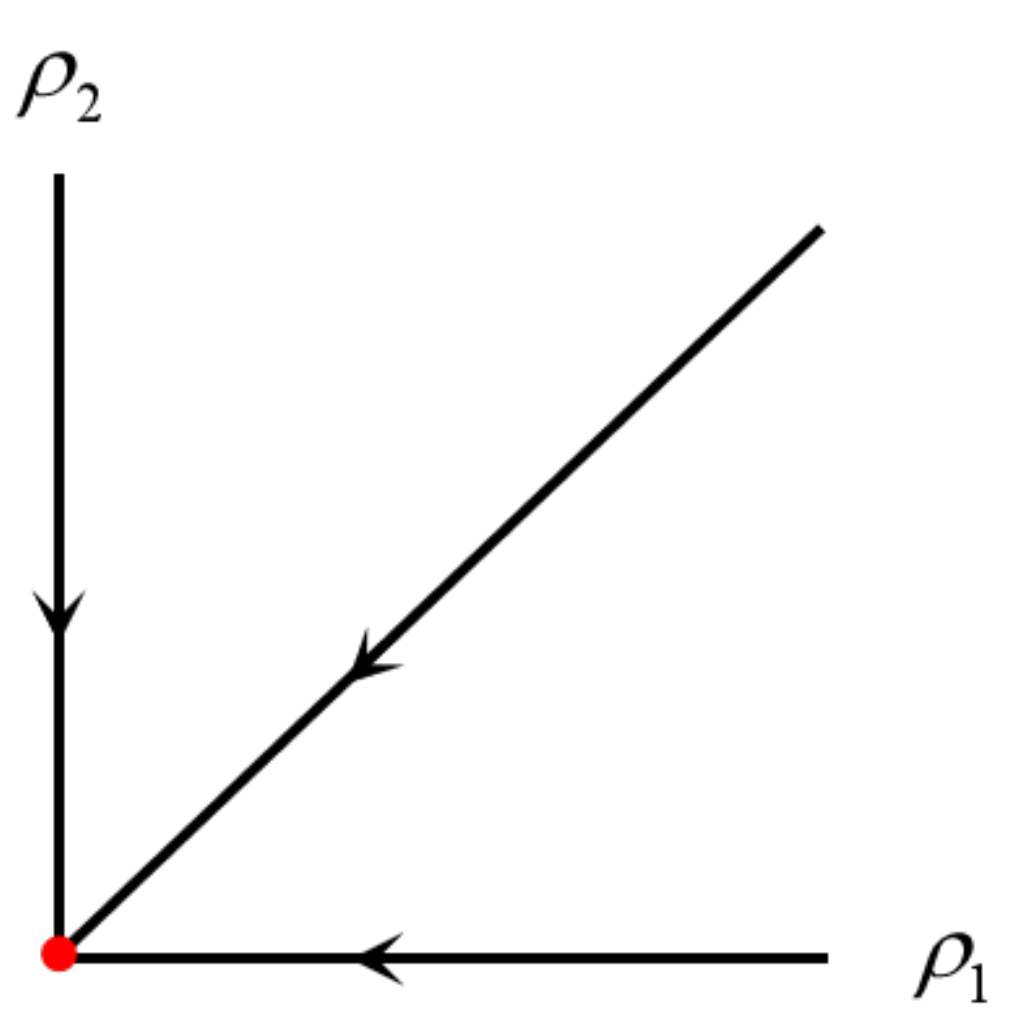}}}\end{minipage}                    &  \begin{minipage}[b]{0.13\columnwidth}\raisebox{-.5\height}{\centerline{\includegraphics[width=0.9\linewidth]{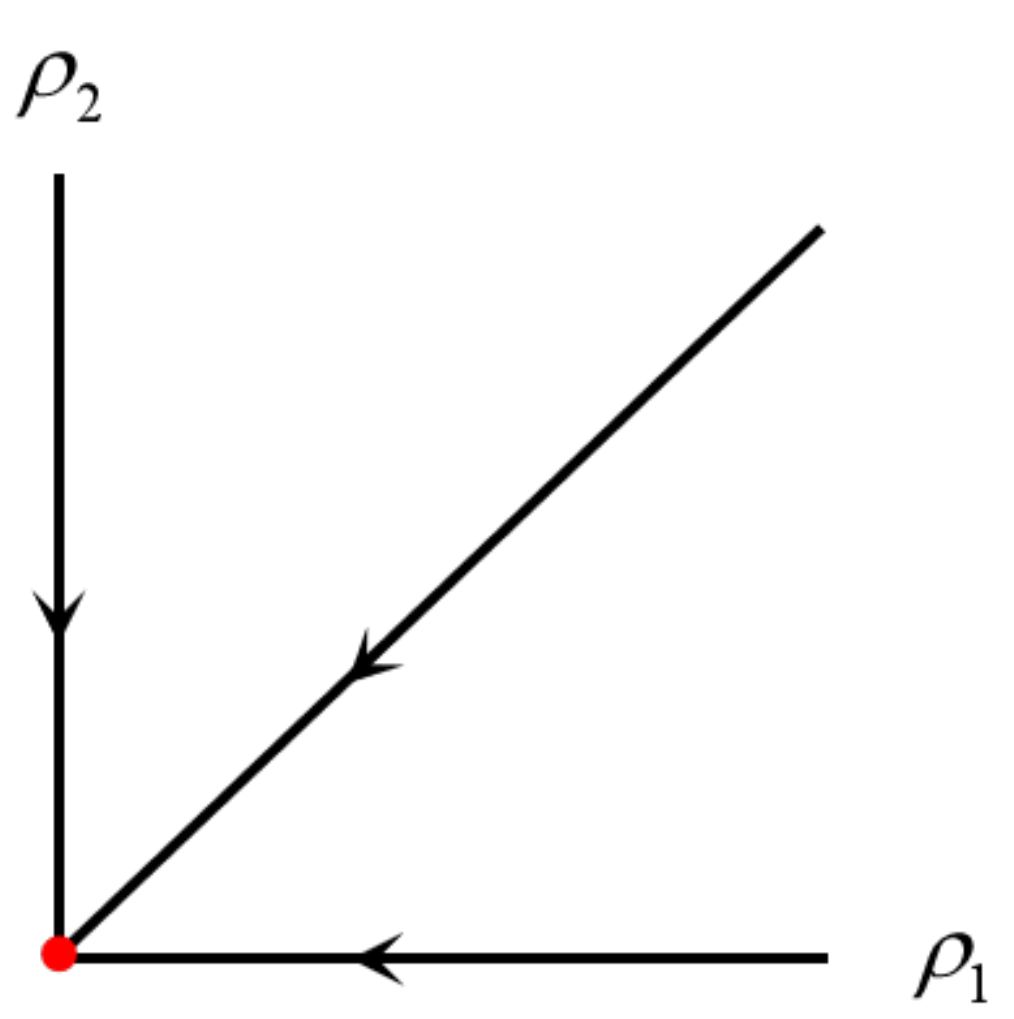}}}\end{minipage}
                                                      & \begin{minipage}[b]{0.13\columnwidth}\raisebox{-.5\height}{\centerline{\includegraphics[width=0.9\linewidth]{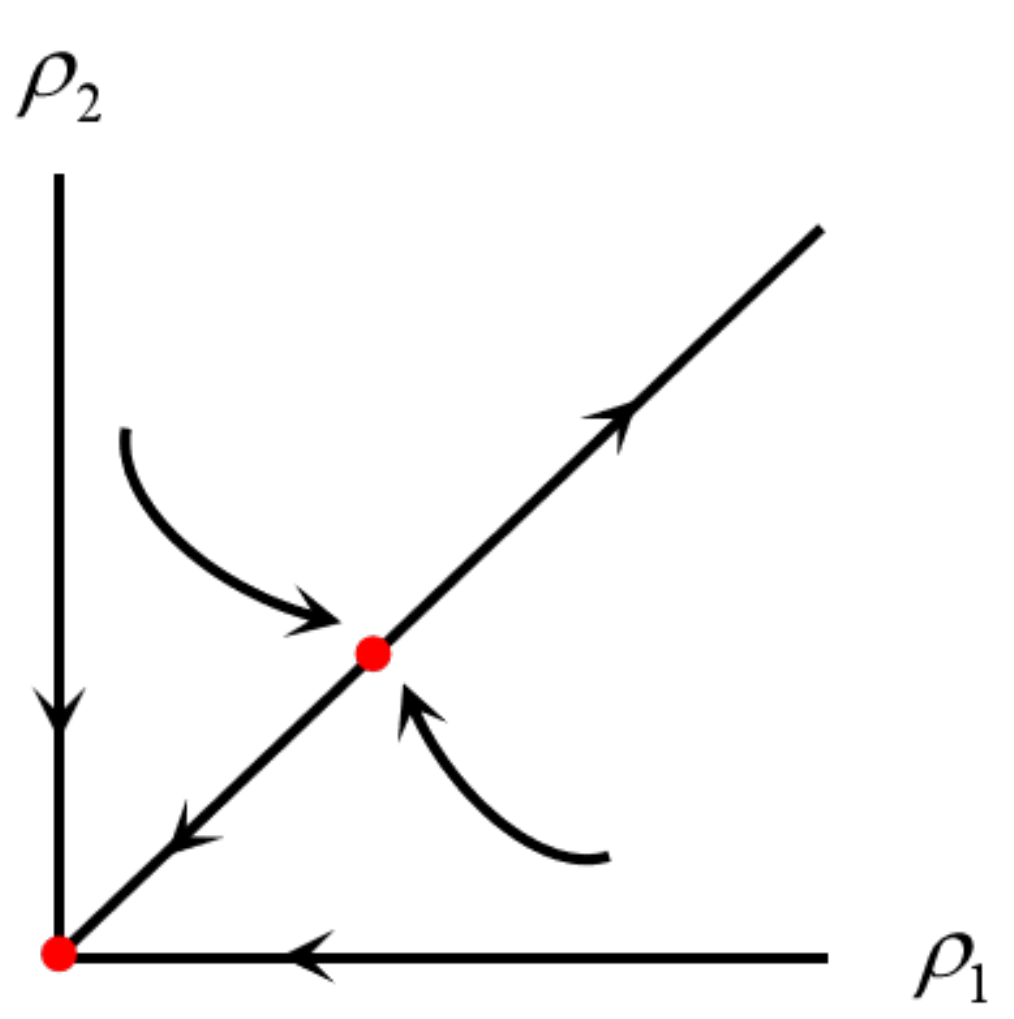}}}\end{minipage}                      & \begin{minipage}[b]{0.13\columnwidth}\raisebox{-.5\height}{\centerline{\includegraphics[width=0.9\linewidth]{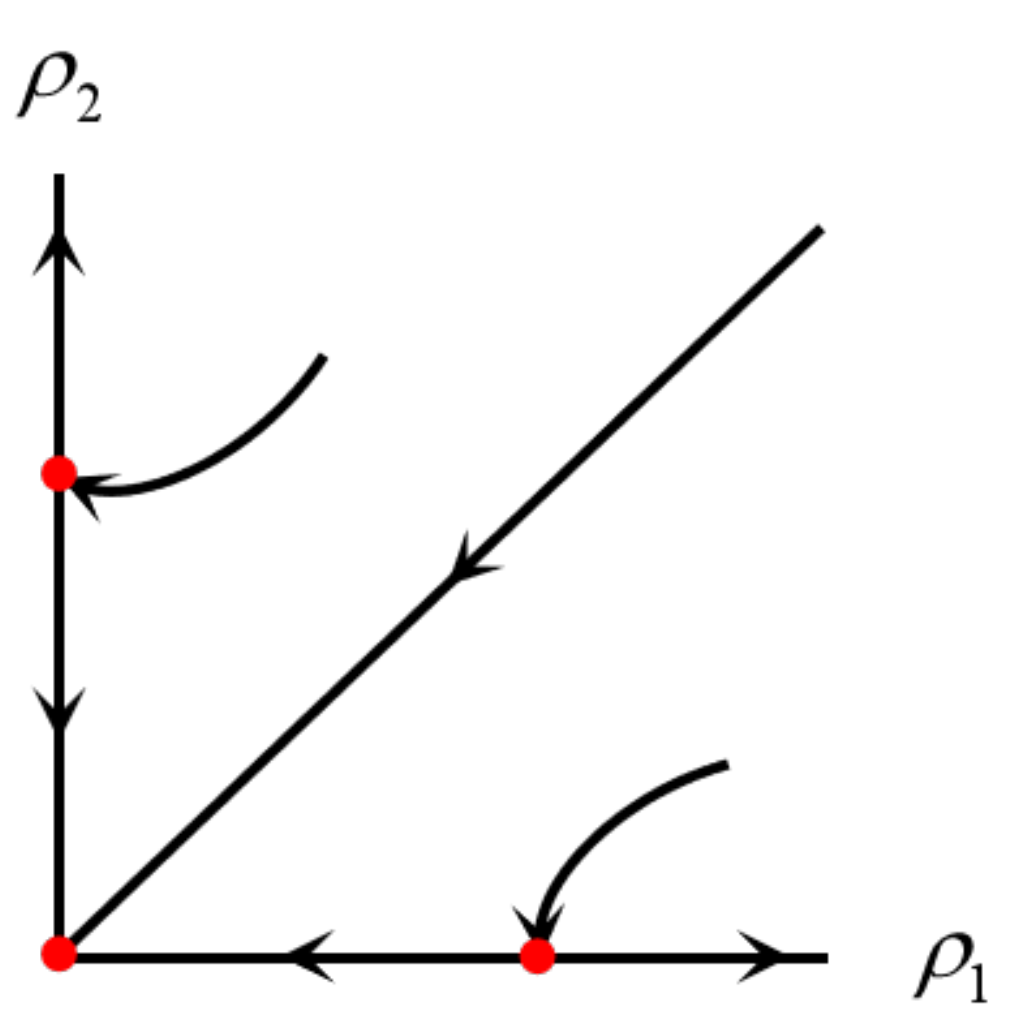}}}\end{minipage}
                                                      & \begin{minipage}[b]{0.13\columnwidth}\raisebox{-.5\height}{\centerline{\includegraphics[width=0.9\linewidth]{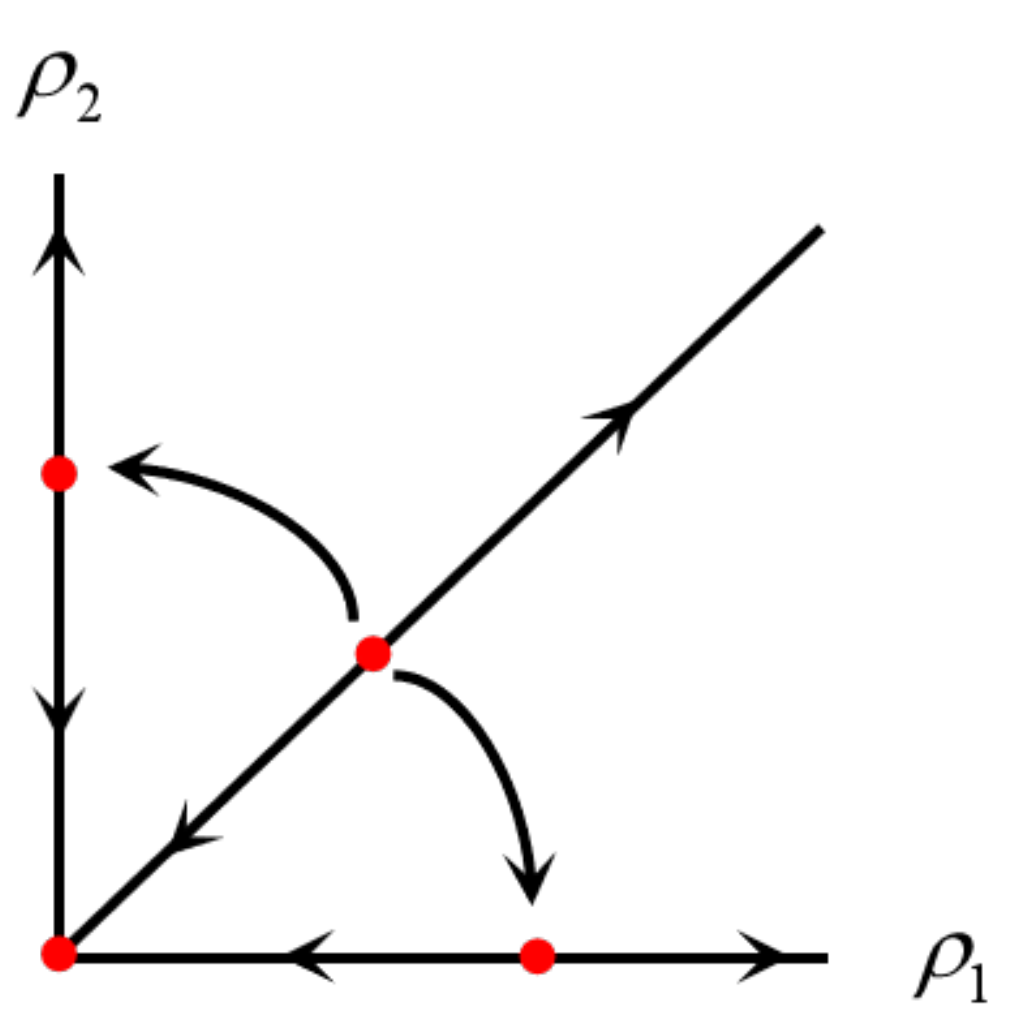}}}\end{minipage}                      & \begin{minipage}[b]{0.13\columnwidth}\raisebox{-.5\height}{\centerline{\includegraphics[width=0.9\linewidth]{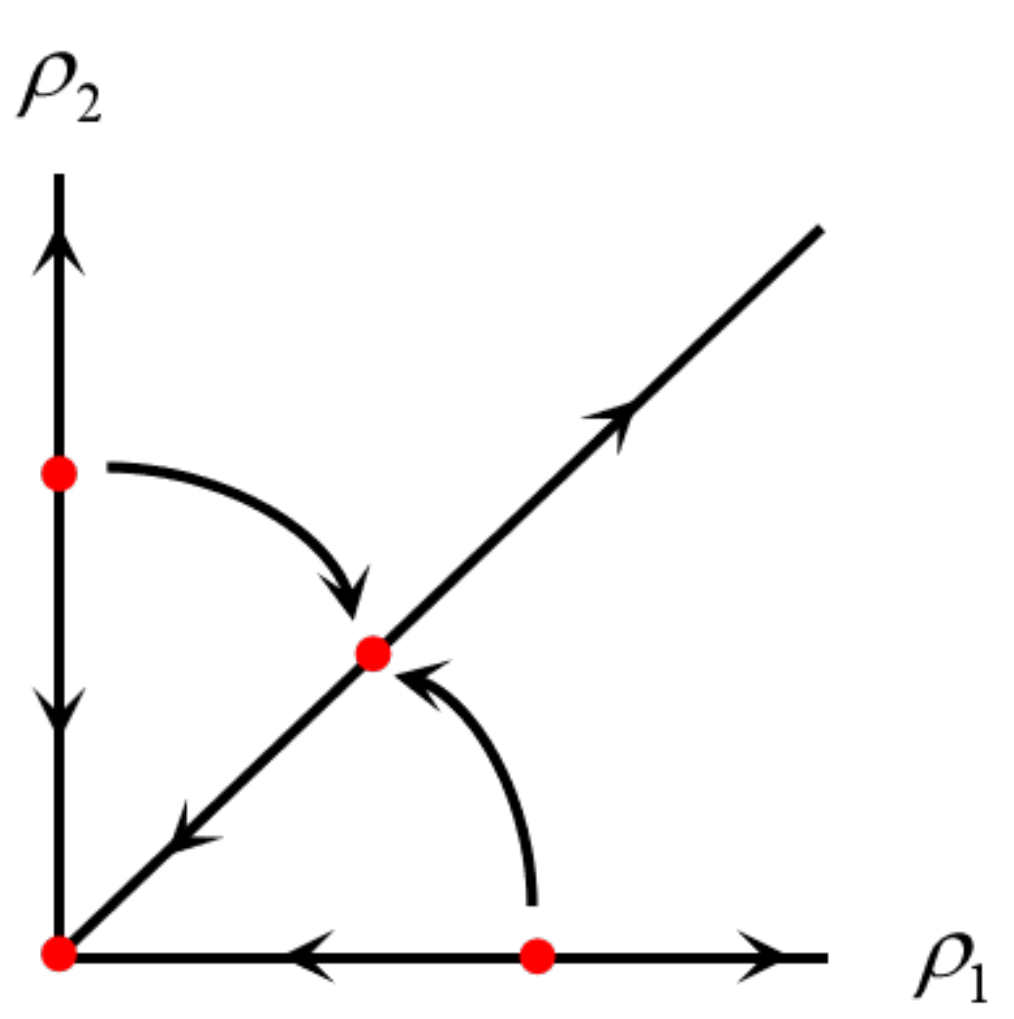}}}\end{minipage}             \\
{$a_1\mu>0$}                                          & \begin{minipage}[b]{0.13\columnwidth}\raisebox{-.5\height}{\centerline{\includegraphics[width=0.9\linewidth]{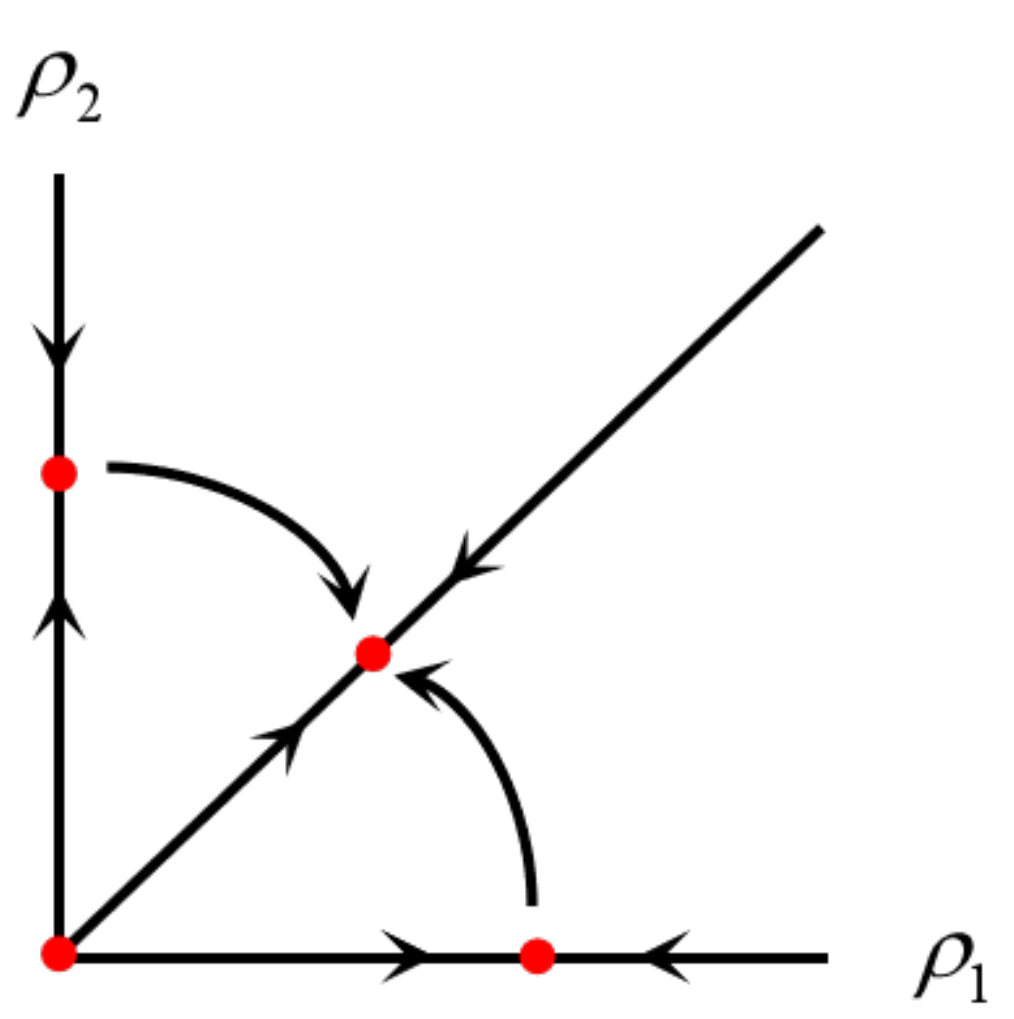}}}\end{minipage}                    &  \begin{minipage}[b]{0.13\columnwidth}\raisebox{-.5\height}{\centerline{\includegraphics[width=0.9\linewidth]{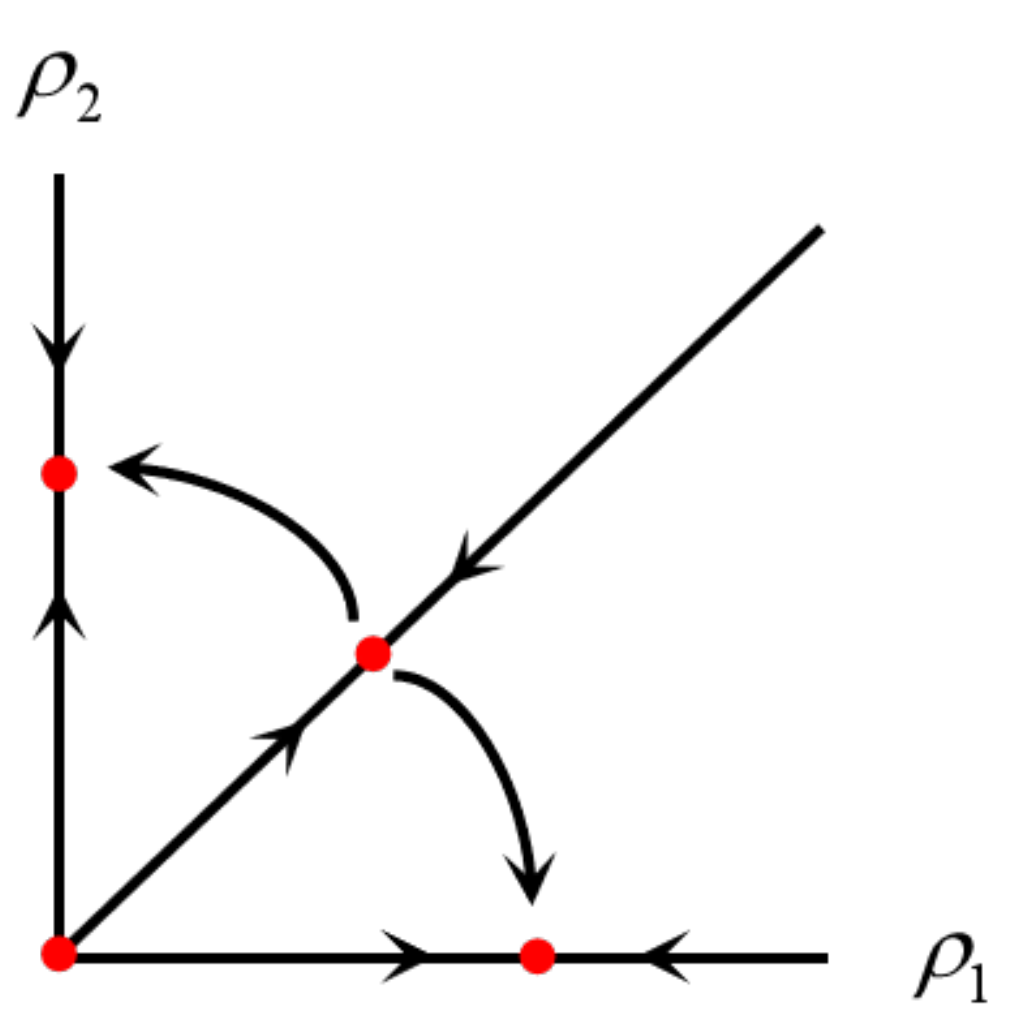}}}\end{minipage}
                                                      & \begin{minipage}[b]{0.13\columnwidth}\raisebox{-.5\height}{\centerline{\includegraphics[width=0.9\linewidth]{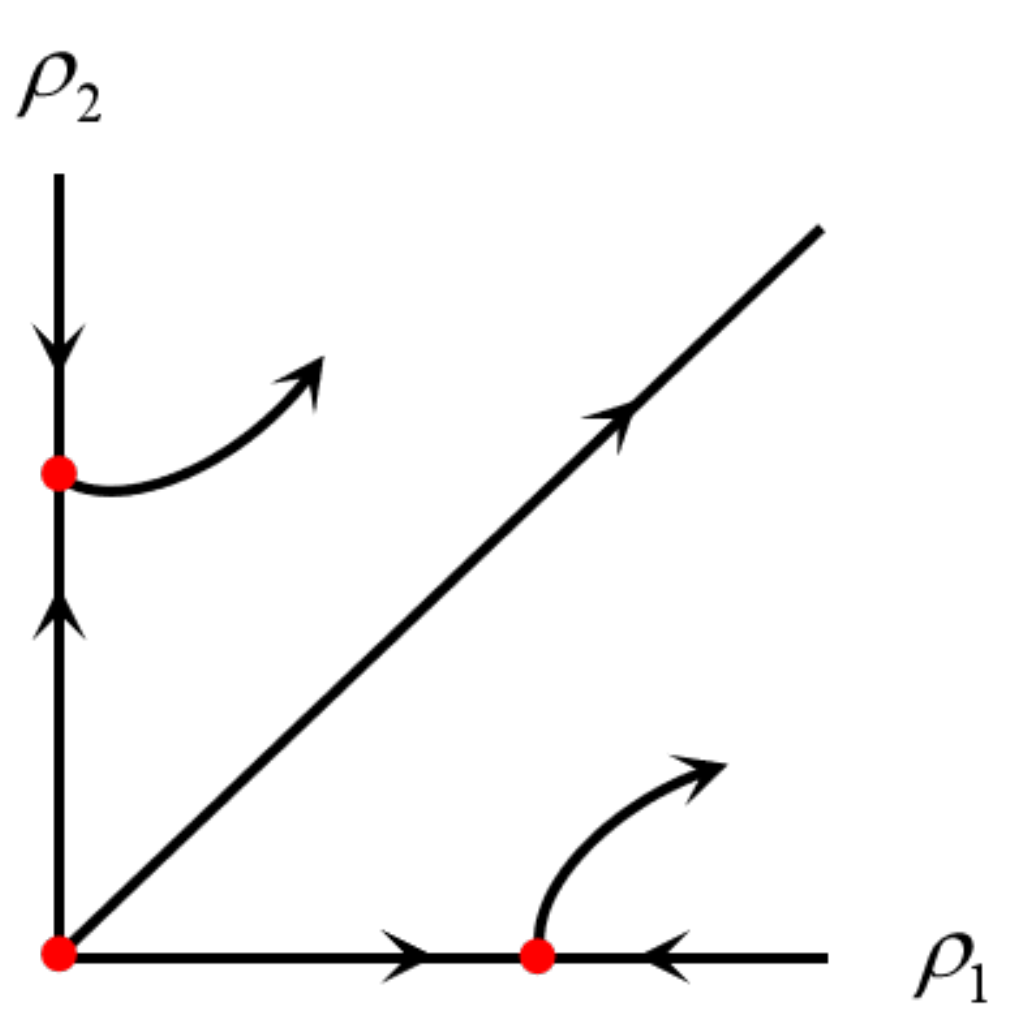}}}\end{minipage}                      & \begin{minipage}[b]{0.13\columnwidth}\raisebox{-.5\height}{\centerline{\includegraphics[width=0.9\linewidth]{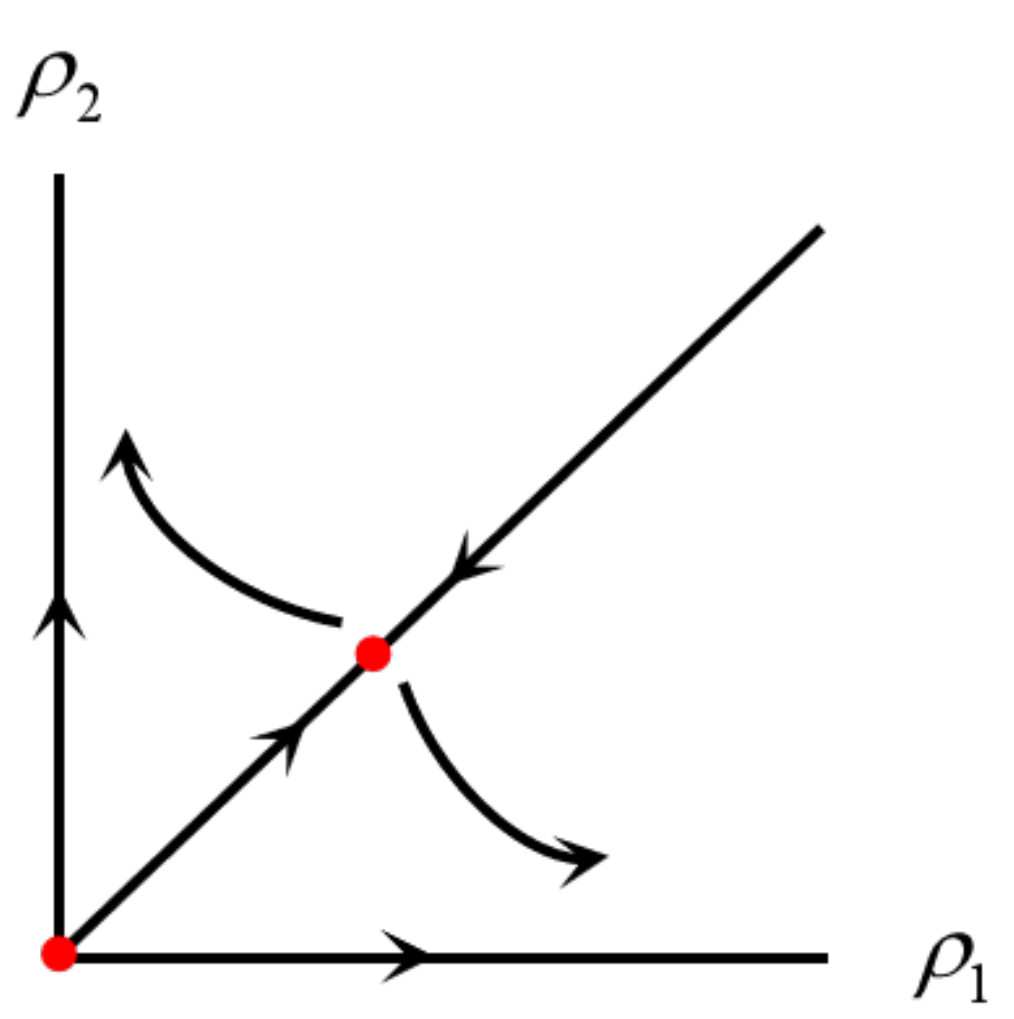}}}\end{minipage}
                                                      & \begin{minipage}[b]{0.13\columnwidth}\raisebox{-.5\height}{\centerline{\includegraphics[width=0.9\linewidth]{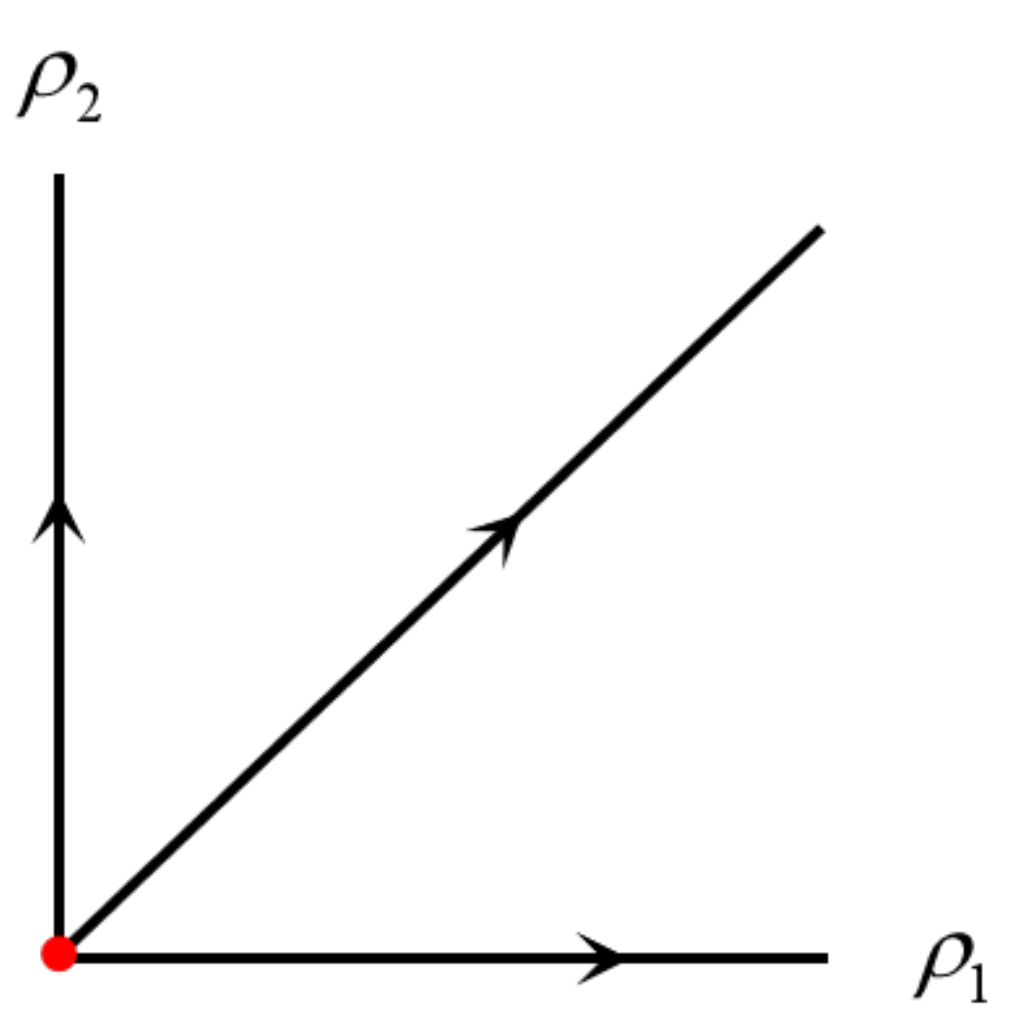}}}\end{minipage}                      & \begin{minipage}[b]{0.13\columnwidth}\raisebox{-.5\height}{\centerline{\includegraphics[width=0.9 \linewidth]{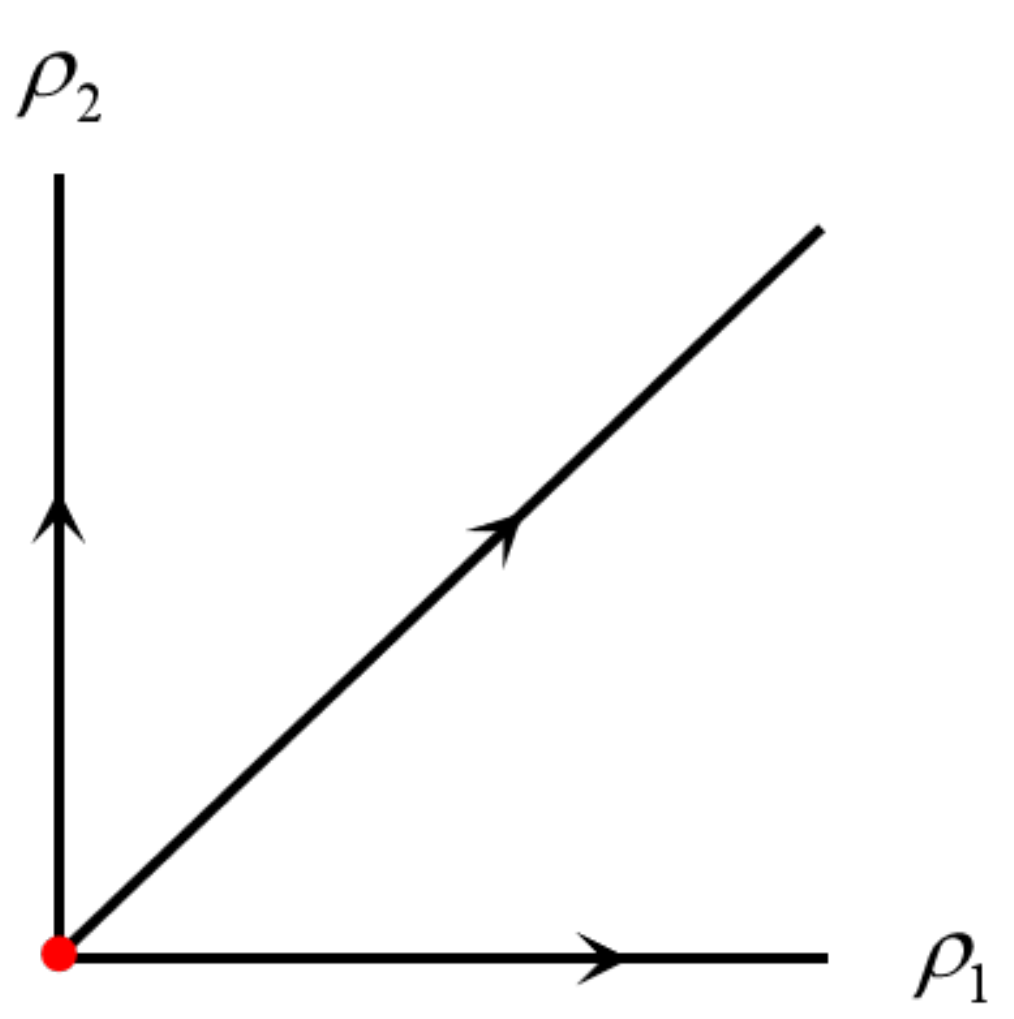}}}\end{minipage}
\end{tabular}
\end{ruledtabular}
\end{table}

Therefore, we can draw the following conclusions.
\begin{theorem}\label{rotating and standing}
We are mainly concerned with the properties corresponding to the following four equilibrium points of (\ref{rho}).\\
$\mathrm{(i)}$~ $(\rho_1,\rho_2)=(0,0)$ corresponds to the origin in the four-dimensional phase space and undergoes a
stationary solution, which is spatially homogeneous.\\
$\mathrm{(ii)}$~ $(\rho_1,\rho_2)=(0,\sqrt{\frac{-a_1\mu}{a_2}})$  corresponds to a periodic solution in the plane of $(z_2,z_3)$, which is spatially inhomogeneous.
At this point, the periodic solution restricted to the center subspace has the following approximate form
$$
U_t(\vartheta)(r,\theta) \approx
\sum_{i=1}^n{2|p_{1i}|\sqrt{\frac{-a_1\mu}{a_2}} J_n(\sqrt{\lambda_{nm}}r)\cos(\mathrm{Arg}(p_{1i})+\omega_{\hat{\lambda}}\vartheta+\omega_{\hat{\lambda}}t+n\theta) {e}_i},
$$
where $e_i$ is the $i$th unit coordinate vector of $\mathbb{R}^n$.
When $a_1\mu>0(<0),~a_2<0(>0)$, system undergoes a rotating wave .
Only when $a_1\mu>0,~a_2+a_3<0,~a_2-a_3>0$, the periodic solution of the equivariant Hopf bifurcation is orbitally asymptotically stable.\\
$\mathrm{(iii)}$~ $(\rho_1,\rho_2)=(\sqrt{\frac{-a_1\mu}{a_2}},0)$ corresponds to a a periodic solution in the plane of $(z_1,z_4)$, which is spatially inhomogeneous.
At this point, the periodic solution restricted to the center subspace has the following approximate form
$$
U_t(\vartheta)(r,\theta) \approx
\sum_{i=1}^n{2|p_{1i}|\sqrt{\frac{-a_1\mu}{a_2}} J_n(\sqrt{\lambda_{nm}}r)\cos(\mathrm{Arg}(p_{1i})+\omega_{\hat{\lambda}}\vartheta+\omega_{\hat{\lambda}}t-n\theta) {e}_i}.
$$
When $a_1\mu>0(<0),~a_2<0(>0)$, system undergoes a rotating wave in the opposite direction as that in (ii).
Its stability conditions are as same as $\mathrm{(ii)}$.\\
$\mathrm{(iv)}$~ $(\rho_1,\rho_2)=(\sqrt{\frac{-a_1\mu}{a_2+a_3}},\sqrt{\frac{-a_1\mu}{a_2+a_3}})$ corresponds to a periodic solution, which is spatially inhomogeneous.
At this point, the periodic solution restricted to the center subspace has the following approximate form
$$
U_t(\vartheta)(r,\theta) \approx
\sum_{i=1}^n{4|p_{1i}|\sqrt{\frac{-a_1\mu}{a_2+a_3}} J_n(\sqrt{\lambda_{nm}}r)\cos(\mathrm{Arg}(p_{1i})+\omega_{\hat{\lambda}}\vartheta+\omega_{\hat{\lambda}}t)\cos(n\theta) {e}_i}.
$$
When $a_1\mu>0(<0),~a_2+a_3<0(>0)$, system undergoes a standing wave.
Only when $a_1\mu>0,~a_2+a_3<0,~a_2-a_3<0$, the periodic solution of the equivariant Hopf bifurcation is orbitally asymptotically stable.
\end{theorem}
The proof is given in Appendix \ref{Proof of Theorem}.

\begin{corollary}\label{qici}
For $\lambda_{0m},m=0,1,2,\cdots$, the corresponding characteristic function is $\phi_{0m}^{c}$. Define $\hat{\phi}_{0m}^c=\frac{\phi_{0m}^{c}}{\|\phi_{0m}^{c}\|_{2,2}}$.
According to  \citep{Faria2000J},
after a similar calculation process shown above, the following normal form on the center manifold is obtained,
$$
\begin{aligned}
&\dot{z}_1=\mathrm{i}\omega_{\hat{\lambda}}z_1+{B}_{11}^{*}z_1\mu+{B}_{2100}^{*}z_1^2z_2+\cdots,\\
&\dot{z}_2=-\mathrm{i}\omega_{\hat{\lambda}}z_2+\overline{{B}^{*}_{11}}z_2\mu+\overline{{B}^{*}_{2100}}z_1z_2^2+\cdots.
\end{aligned}
$$
Introducing a set of polar coordinates, we can get
$$
\dot{\rho}=({a}_{1}^{*}\mu+{a}_2^{*}\rho^2)\rho+o(\mu^2\rho+|(\rho,\mu)|^4),
$$
with
$$
a_1^{*}=\mathrm{Re}\{B_{11}^{*}\},~a_2^{*}=\mathrm{Re}\{B_{2100}^{*}\},\\
$$
where the specific representation of $B_{11}^{*}$ and $B_{2100}^{*}$ is shown in \citep{Faria2000J}.

Besides, we get that \\
$\mathrm{(i)}$~ When ${a}_{2}^{*}<0(>0)$, the periodic solution is orbitally asymptotically stable(unstable).\\
$\mathrm{(ii)}$~ When ${a}_{1}^{*}{a}_{2}^{*}<0(>0)$, the bifurcation is supercritical(subcritical).
\end{corollary}

\subsection{Explicit formulas for a class of reaction-diffusion model with discrete time delay}\label{Explicit formulas}

In a specific model, it is necessary to calculate $A_{p_1p_2p_3p_4},~S_{y(0)z_k},~S_{y(-1)z_k},~k=1,2,3,4$ to determine the explicit expression of $B_{p_1p_2p_3p_4}$ in the normal form. Therefore, in order to provide a more general symbolic expression, we will consider a class of reaction-diffusion system with discrete time delay defined on a disk as follows
\begin{equation}\label{time-delayed reaction-diffusion system r theta}
\left\{\begin{array}{l}
\frac{\partial u(t, r, \theta)}{\partial t}=d_{1} \Delta_{r \theta} u(t, r, \theta)+F^{(1)}(u(t, r, \theta), v(t, r, \theta)),~(r, \theta) \in \mathbb{D},~t>0, \\
\frac{\partial v(t, r, \theta)}{\partial t}=d_{2} \Delta_{r \theta} v(t, r, \theta)+F^{(2)}(u(t, r, \theta), v(t, r, \theta), u(t-\tau, r, \theta), v(t-\tau, r, \theta)),~(r, \theta) \in \mathbb{D},~t>0, \\
\partial_{r} u(\cdot, R, \theta)=\partial_{r} v(\cdot, R, \theta)=0,~\theta \in [0,2\pi).
\end{array}\right.
\end{equation}
This type of model covers some predator-prey systems and chemical reaction models, etc. While in practice there are many ways to introduce the time delay $\tau$,  we demonstrate the critical method of analysis by including $\tau$ in the second equation for simplicity. Other types of systems can also refer to this process for calculation.

Assume that the model has a positive equilibrium point $E^*(u^*,v^*)$ and select the time delay $\tau$ as the bifurcation parameter. Letting $\bar{u}(t, r, \theta)=u(\tau t, r, \theta)-u^{*}, \bar{v}(t, r, \theta)=v(\tau t, r, \theta)-v^{*}$, we drop the bar for simplicity. Then system (\ref{time-delayed reaction-diffusion system r theta}) can be transformed into
\begin{equation}\label{time-delayed reaction-diffusion system ijkl}
 \left\{\begin{aligned}
\frac{\partial u(t, r, \theta)}{\partial t}=& \tau d_{1} \Delta_{r \theta} u(t, r, \theta)+\tau\left[a_{11}\left(u(t, r, \theta)+u^{*}\right)+a_{12}\left(v(t, r, \theta)+v^{*}\right)\right] \\
&+\tau \sum_{i+j \geq 2} \frac{1}{i ! j !} F_{i j}^{(1)}(0,0) u^{i}(t, r, \theta) v^{j}(t, r, \theta), \\
\frac{\partial v(t, r, \theta)}{\partial t}=& \tau d_{2} \Delta_{r \theta} v(t, r, \theta)+\tau\left[
a_{21}\left(u(t, r, \theta)+u^{*}\right)+a_{22}\left(v(t, r, \theta)+v^{*}\right)
\right] \\
+& \tau\left[
b_{21}\left(u(t-1, r, \theta)+u^{*}\right)+b_{22}\left(v(t-1, r, \theta)+v^{*}\right)
\right] \\
+& \tau \sum_{i+j+l+k \geq 2} \frac{1}{i ! j ! k ! l !} F_{i j k l}^{(2)}(0,0,0,0) u^{i}(t, r, \theta) v^{j}(t, r, \theta) u^{k}(t-1, r, \theta) v^{l}(t-1, r, \theta),
\end{aligned}\right.
\end{equation}
with
$$
\begin{gathered}
\left(\begin{array}{ll}
a_{11} & a_{12} \\
a_{21} & a_{22}
\end{array}\right)=\left(\begin{array}{cc}
\frac{\partial F^{(1)}\left(u^{*}, v^{*}\right)}{\partial u(t)} & \frac{\partial F^{(1)}\left(u^{*}, v^{*}\right)}{\partial v(t)} \\
\frac{\partial F^{(2)}\left(u^{*}, v^{*},u^{*}, v^{*}\right)}{\partial u(t)} & \frac{\partial F^{(2)}\left(u^{*}, v^{*},u^{*}, v^{*}\right)}{\partial v(t)}
\end{array}\right), \\
\left(\begin{array}{ll}
b_{11} & b_{12} \\
b_{21} & b_{22}
\end{array}\right)=\left(\begin{array}{cc}
0 & 0\\
\frac{\partial F^{(2)}\left(u^{*}, v^{*},u^{*}, v^{*}\right)}{\partial u(t-\tau)} & \frac{\partial F^{(2)}\left(u^{*}, v^{*},u^{*}, v^{*}\right)}{\partial v(t-\tau)}
\end{array}\right),\\
F_{i j}^{(1)}=\frac{\partial^{i+j} F^{(1)}}{\partial u^{i} \partial v^{j}}(0,0), \\
F_{i j k l}^{(2)}=\frac{\partial^{i+j+k+l} F^{(2)}}{\partial u^{i} \partial v^{j} \partial u^{k}(t-\tau) \partial v^{l}(t-\tau)}(0,0,0,0).
\end{gathered}
$$

Letting $\tau=\hat{\tau}+\mu$, where $\mu \in \mathbb{R}$ and $\hat{\tau}$ is the critical values at which Hopf bifurcations occur,  then system (\ref{time-delayed reaction-diffusion system ijkl}) can be written in an abstract form like (\ref{abstract functional differential equation}), where operators
${L}_0$ and $\tilde{F}$ are given, respectively, by
$$
L_0(\varphi)=(\hat{\tau}+\mu)\left(\begin{array}{c}
a_{11} \varphi_{1}(0)+a_{12} \varphi_{2}(0) \\
b_{21} \varphi_{1}(-1)+b_{22} \varphi_{2}(-1)+a_{21} \varphi_{1}(0)+a_{22} \varphi_{2}(0)
\end{array}\right),
$$
$$
\tilde{F}(\varphi, \mu)=(\hat{\tau}+\mu)\left(\begin{array}{c}\sum_{i+j \geq 2} \frac{1}{i ! j !} F_{i j}^{(1)}(0,0) \varphi_{1}^{i}(0) \varphi_{2}^{j}(0) \\ \sum_{i+j+k+l \geq 2} \frac{1}{i ! j ! k ! l !} F_{i j k l}^{(2)}(0,0,0,0) \varphi_{1}^{i}(0) \varphi_{2}^{j}(0) \varphi_{1}^{k}(-1) \varphi_{2}^{l}(-1)\end{array}\right),
$$
$$
\varphi=(\varphi_1,\varphi_2)^{\rm{T}} \in \mathscr{C}.
$$

Choosing $\xi=(1,p_0)^{T}$, with $p_0=\frac{\mathrm{i}\omega_{\hat{\lambda}}+d_1\lambda_{nm}-a_{11}}{a_{12}}$, we get that the bases of $P$ is
$$
\Phi_{r\theta}(\vartheta)=\left(\begin{array}{rrrr}
\mathrm{e}^{\mathrm{i}\omega_{\hat{\lambda}}\hat{\tau}\vartheta}\hat{\phi}_{nm}^c & \mathrm{e}^{-\mathrm{i}\omega_{\hat{\lambda}}\hat{\tau}\vartheta}\hat{\phi}_{nm}^c
& \mathrm{e}^{\mathrm{i}\omega_{\hat{\lambda}}\hat{\tau}\vartheta}\hat{\phi}_{nm}^s & \mathrm{e}^{-\mathrm{i}\omega_{\hat{\lambda}}\hat{\tau}\vartheta}\hat{\phi}_{nm}^s\\
p_0\mathrm{e}^{\mathrm{i}\omega_{\hat{\lambda}}\hat{\tau}\vartheta}\hat{\phi}_{nm}^c & \bar{p}_0\mathrm{e}^{-\mathrm{i}\omega_{\hat{\lambda}}\hat{\tau}\vartheta}\hat{\phi}_{nm}^c
& p_0\mathrm{e}^{\mathrm{i}\omega_{\hat{\lambda}}\hat{\tau}\vartheta}\hat{\phi}_{nm}^s & \bar{p}_0\mathrm{e}^{-\mathrm{i}\omega_{\hat{\lambda}}\hat{\tau}\vartheta}\hat{\phi}_{nm}^s
\end{array}\right),
$$
and the basis of $P^*$ can also be obtained.
The explicit formulas of $A_{p_1p_2p_3p_4},~S_{y(0)z_k},~S_{y(-1)z_k},~k=1,2,3,4$, and $h_{jp_1p_2p_3p_4}$ in the calculation of normal form are shown in Appendices \ref{A_{p_1p_2p_3p_4}}, \ref{S} and \ref{h}.

\section{Numerical simulations}\label{sec4}

\subsection{Numerical example 1: A diffusive Brusselator model with delayed feedback }\label{sec5.1}
In~\citep{Zuo2012J}, Zuo and Wei studied a diffusive Brusselator model with delayed feedback. We put this model on a disk with Neumann boundary conditions and perform some numerical simulations. The model in polar form now turns to be
\begin{equation}\label{Brusselator}
\left\{\begin{array}{l}
\frac{\partial u(t, r, \theta)}{\partial t}=d_{1} \Delta_{r \theta} u(t, r, \theta)+a-(b+1) u(t, r, \theta)+u^2(t, r, \theta) v(t, r, \theta), \\
\frac{\partial v(t, r, \theta)}{\partial t}=d_{2} \Delta_{r \theta} v(t, r, \theta)+b u(t, r, \theta)-u^2(t, r, \theta) v(t, r, \theta)+g (v(t-\tau, r, \theta)- v(t, r, \theta)),\\
~~~~~~~~~~~~~~~~~~~~~~~~~~~~~~~~~~~~~~~~~~~~~~~~~~~~~~~~~~~~~~~~~~~~~~~~~~~~~~~~~~~~~~~~~~~~~~~~~~~~~~~~~~~~~~(r, \theta) \in \mathbb{D},~t>0,\\
\partial_{r} u(\cdot, R, \theta)=\partial_{r} v(\cdot, R, \theta)=0,~\theta \in [0,2\pi).
\end{array}\right.
\end{equation}

\begin{figure}[htbp]
\centering
(a)\includegraphics[width=0.46\textwidth]{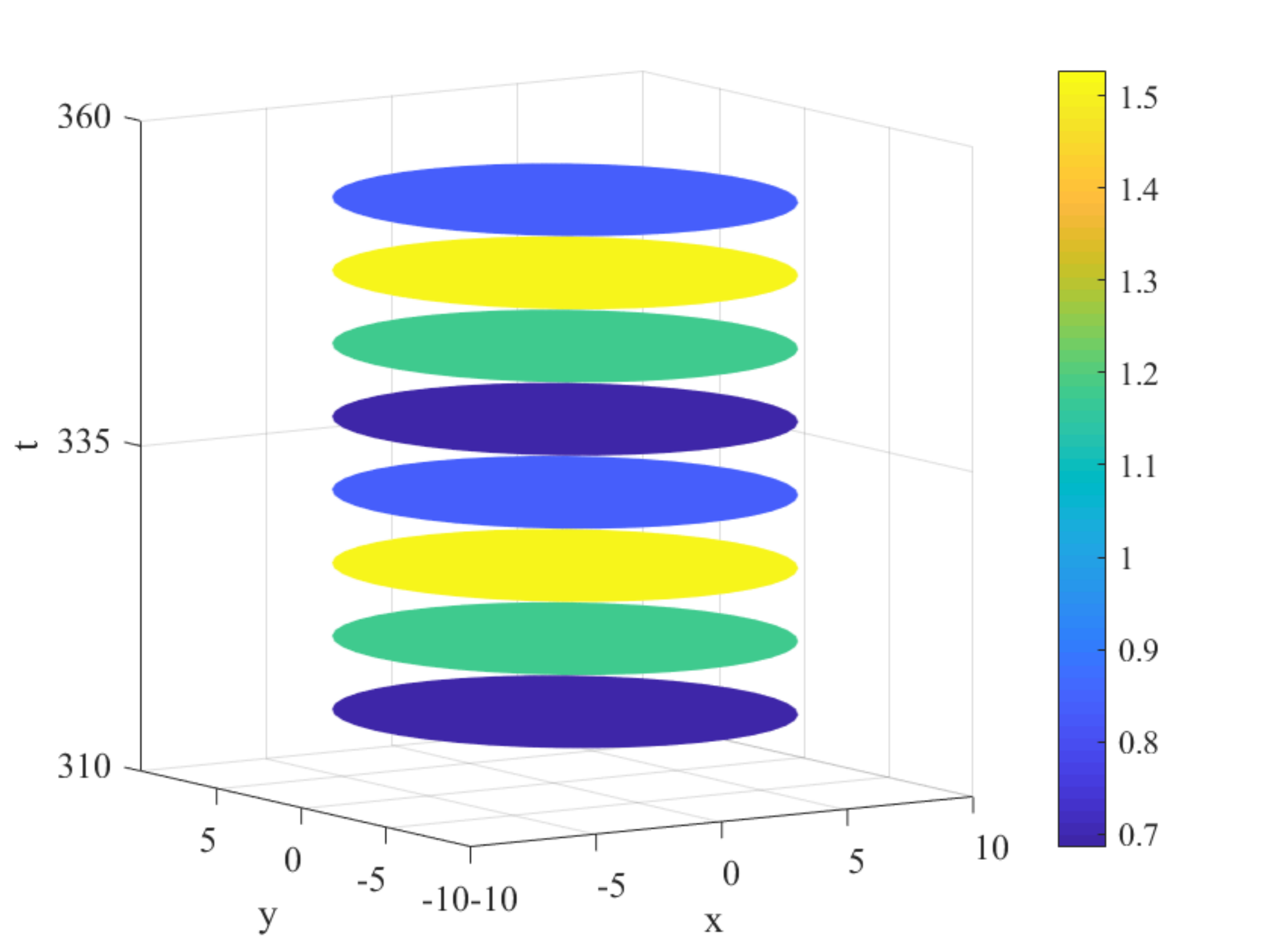}
(b)\includegraphics[width=0.46\textwidth]{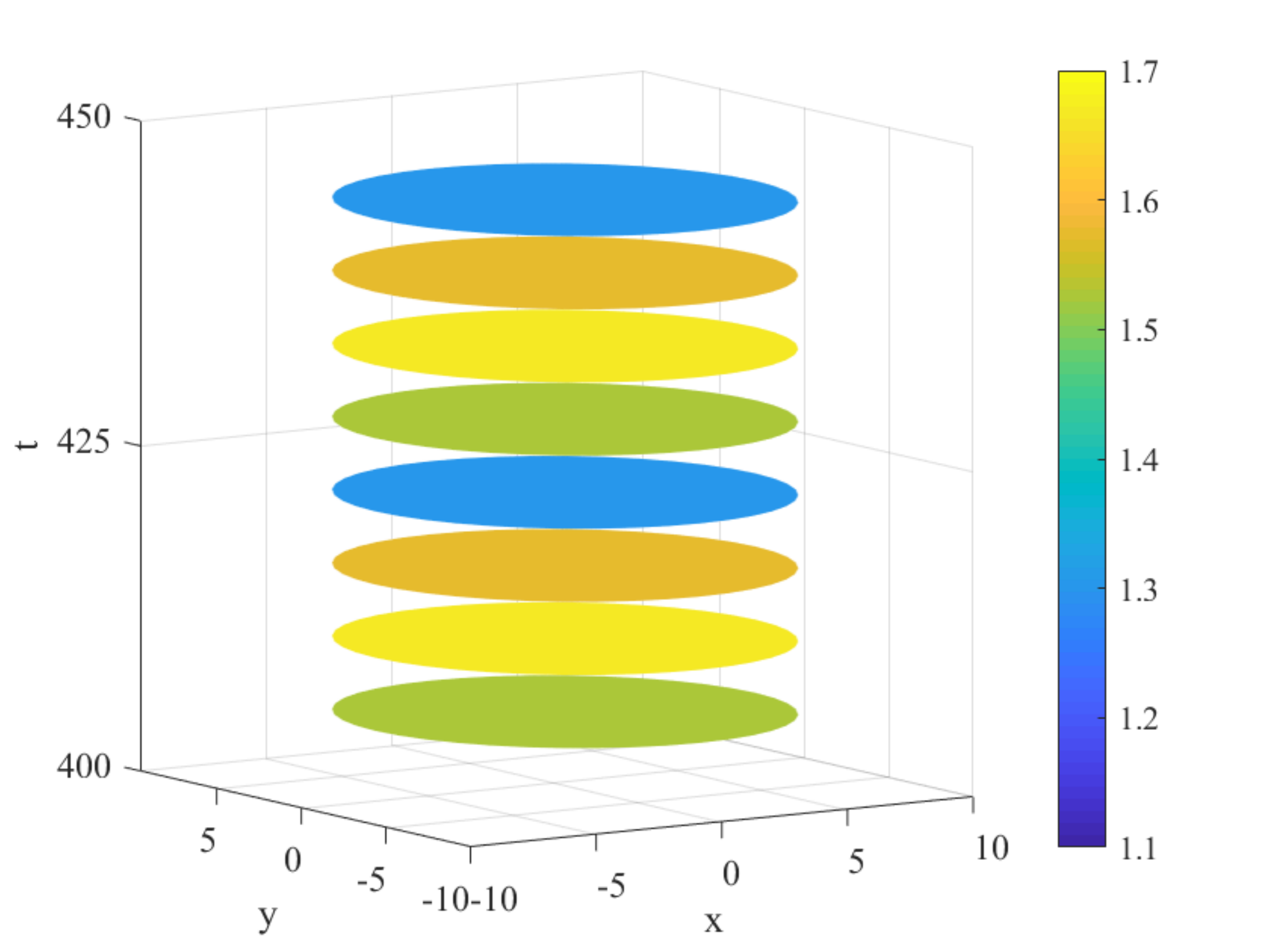}
\caption{Spatially homogeneous periodic solutions with parameters: $a=1,~b=1.5,~g=2,~d_1=2,~d_2=5,~R=10,\tau=2.~(a):u,~(b):v.$}
\label{unstable2}
\end{figure}

Fixing $a=1,~b=1.5,~g=2,~d_1=2,~d_2=5,~R=10$, we get that the unique positive equilibrium solution of the model is $(1,1.5)$. When $\hat{\lambda}=\lambda_{00}$, $\omega \approx 0.6166$ and $\hat{\tau} \approx 0.7128$. According to the common analysis of the standard Hopf bifurcation, we get that when $\tau \in \left[0,0.7128\right)$, $E^*$ is locally asymptotically stable. When $\tau$ increases from zero and crosses the critical value $\hat{\tau} \approx 0.7128$, a family of periodic solutions are bifurcated from $E^*$, which is spatially homogeneous
(see Figure \ref{unstable2} at $\tau=2$). By Corollary \ref{qici}, the Hopf bifurcation is supercritical and the periodic solutions are stable since ${a}_1^{*}{a}_2^{*} \approx -0.8264,~{a}_2^{*} \approx -0.6920$.

\subsection{Numerical example 2: A delayed predator-prey model with group defense and nonlocal competition}\label{sec5.2}
In \citep{Liu2020J}, the authors investigated a predator-prey model with group defence and nonlocal competition. Here, we use the method established above to investigate the dynamics of such a model on a disk.

\begin{equation}\label{nonlocal}
\left\{\begin{array}{l}
\frac{\partial u(t, r, \theta)}{\partial t}=d_{1} \Delta_{r \theta} u(t, r, \theta)+b u(t,r,\theta)\left(1-\frac{\hat{u}(t,r,\theta)}{K}\right)-a u^{\alpha}(t,r,\theta)v(t,r,\theta),~(r, \theta) \in \mathbb{D},~t>0, \\
\frac{\partial v(t, r, \theta)}{\partial t}=d_{2} \Delta_{r \theta} v(t, r, \theta)-d v(t,r,\theta)+ a eu^{\alpha}(t-\tau,r,\theta)v(t,r,\theta),~(r, \theta) \in \mathbb{D},~t>0, \\
\partial_{r} u(\cdot, R, \theta)=\partial_{r} v(\cdot, R, \theta)=0,~\theta \in [0,2\pi),
\end{array}\right.
\end{equation}
where the state $u$, $v$  and parameters are defined in \citep{Liu2020J}.  In particular, $\hat{u}(r,\theta,t)$ depicts the  nonlocal competition with the  form of
$$
\hat{u}(r,\theta,t)=\frac{1}{ \pi R^2} \int_{0}^{R} \int_{0}^{2 \pi} \bar{r} u\left(\bar{r},\bar{\theta},t\right)\mathrm{d} \bar{\theta} \mathrm{d} \bar{r}.
$$
According to subsection \ref{sec2.2}, characteristic equations of the linearization equation at the positive equilibrium point $E^*$ for system (\ref{nonlocal}) are
\begin{equation}\label{characteristic equations}
\left\{\begin{array}{l}
\gamma^2+P_{0m}\gamma+Q_{0m}-a_{12}b_{21}\mathrm{e}^{-\gamma \tau}=0,~m=0,1,2,\cdots,\\
\left(\gamma^2+\bar{P}_{nm}\gamma+\bar{Q}_{nm}-a_{12}b_{21}\mathrm{e}^{-\gamma \tau}\right)^2=0,~n=1,2,\cdots,~m=1,2,\cdots,
\end{array}\right.
\end{equation}
with
$$
\begin{aligned}
&P_{0m}=(d_1+d_2)\lambda_{0m}^2-a_{11}-c_{11},~{Q}_{0m}=(d_1\lambda_{0m}^2-a_{11}-c_{11})d_2\lambda_{0m}^2~m=0,1,2,\cdots;\\
&\bar{P}_{nm}=(d_1+d_2)\lambda_{nm}^2-a_{11},~\bar{Q}_{nm}=(d_1\lambda_{nm}^2-a_{11})d_2\lambda_{nm}^2,~n=1,2,\cdots,~m=1,2,\cdots,
\end{aligned}
$$
where the expressions of $a_{11},~a_{12},~b_{21},~c_{11}$ are shown in \citep{Liu2020J}.

\begin{figure}[htbp]
\includegraphics[width=0.74\textwidth]{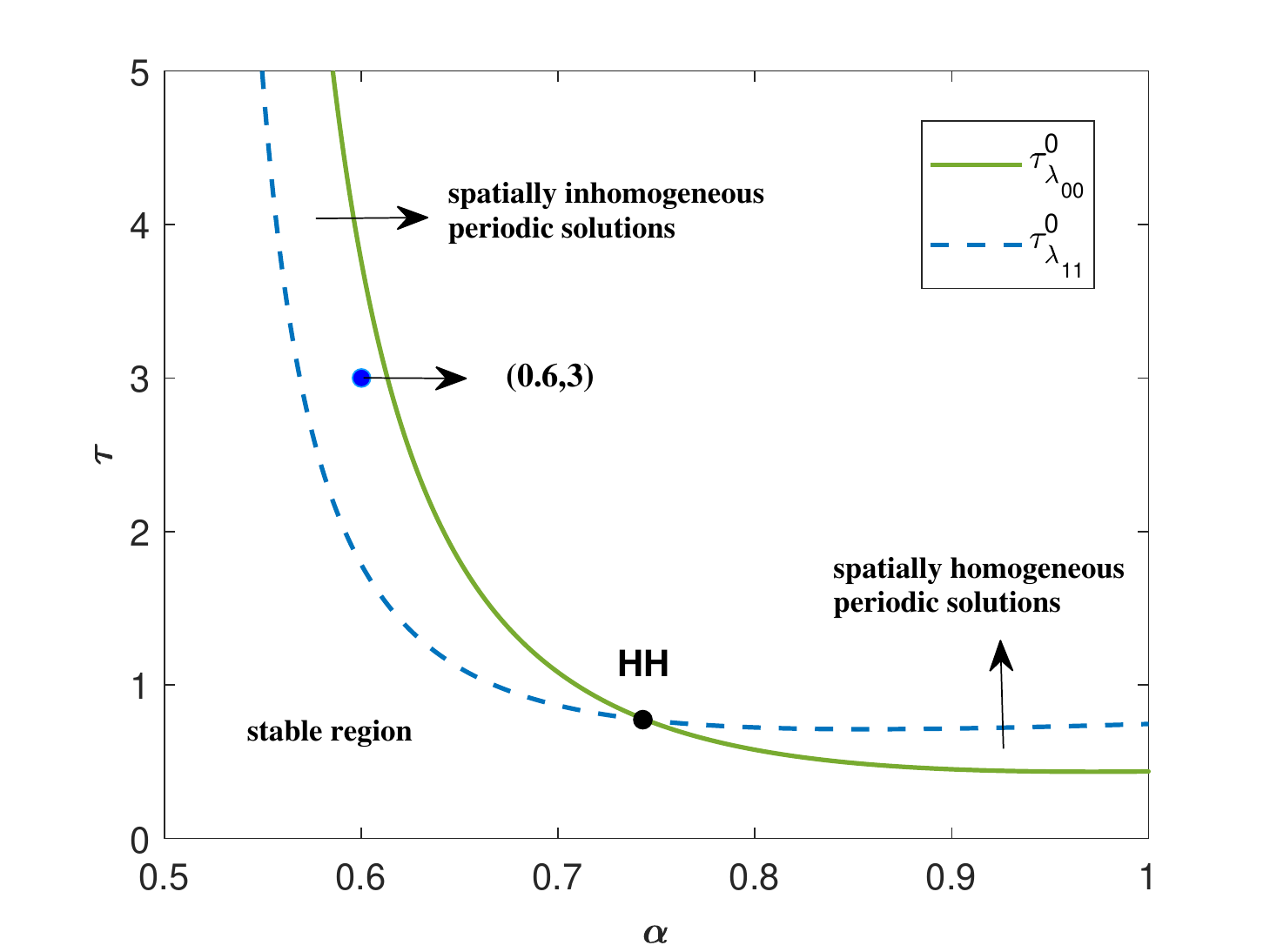}
\caption{Partial bifurcation curves on the $\alpha-\tau$ plane.}\label{nonlocalfenzhitu}
\end{figure}

\begin{figure}[htbp]
(a)\includegraphics[width=1\textwidth]{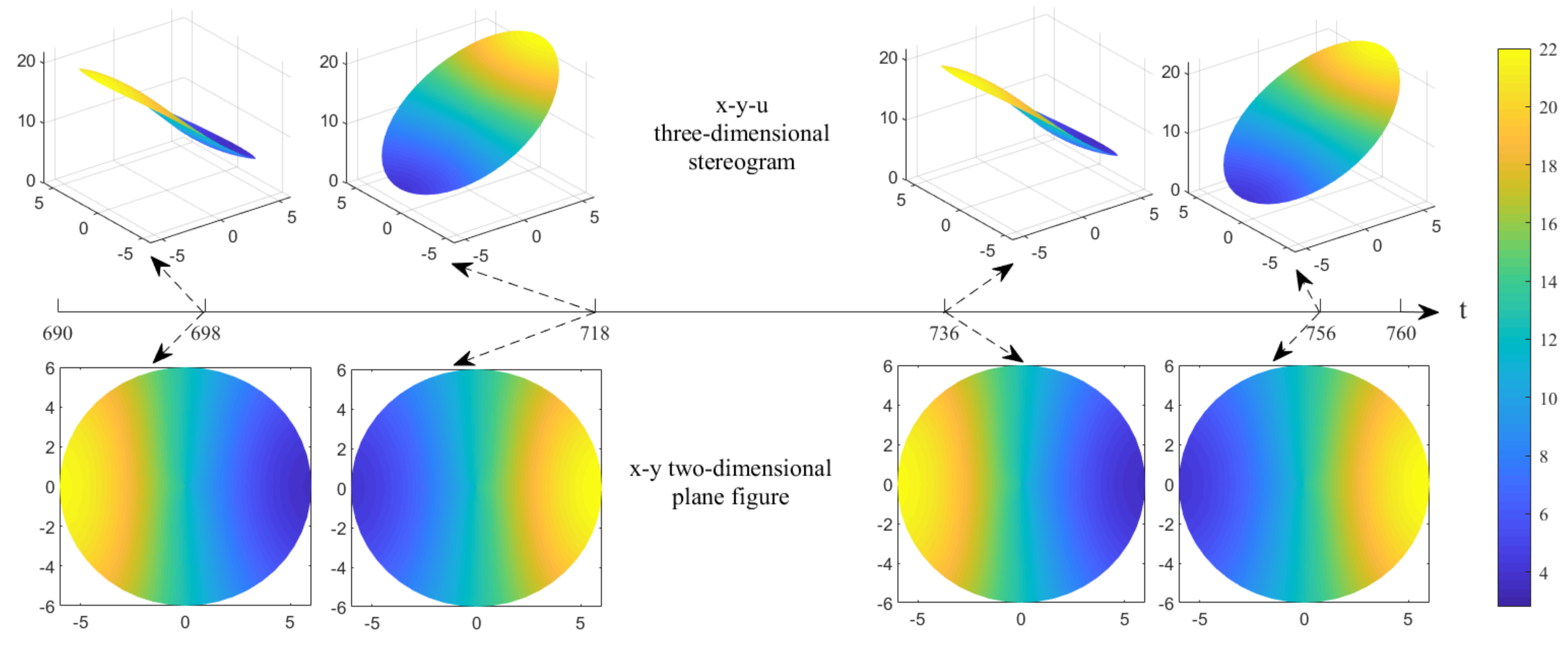}
(b)\includegraphics[width=1\textwidth]{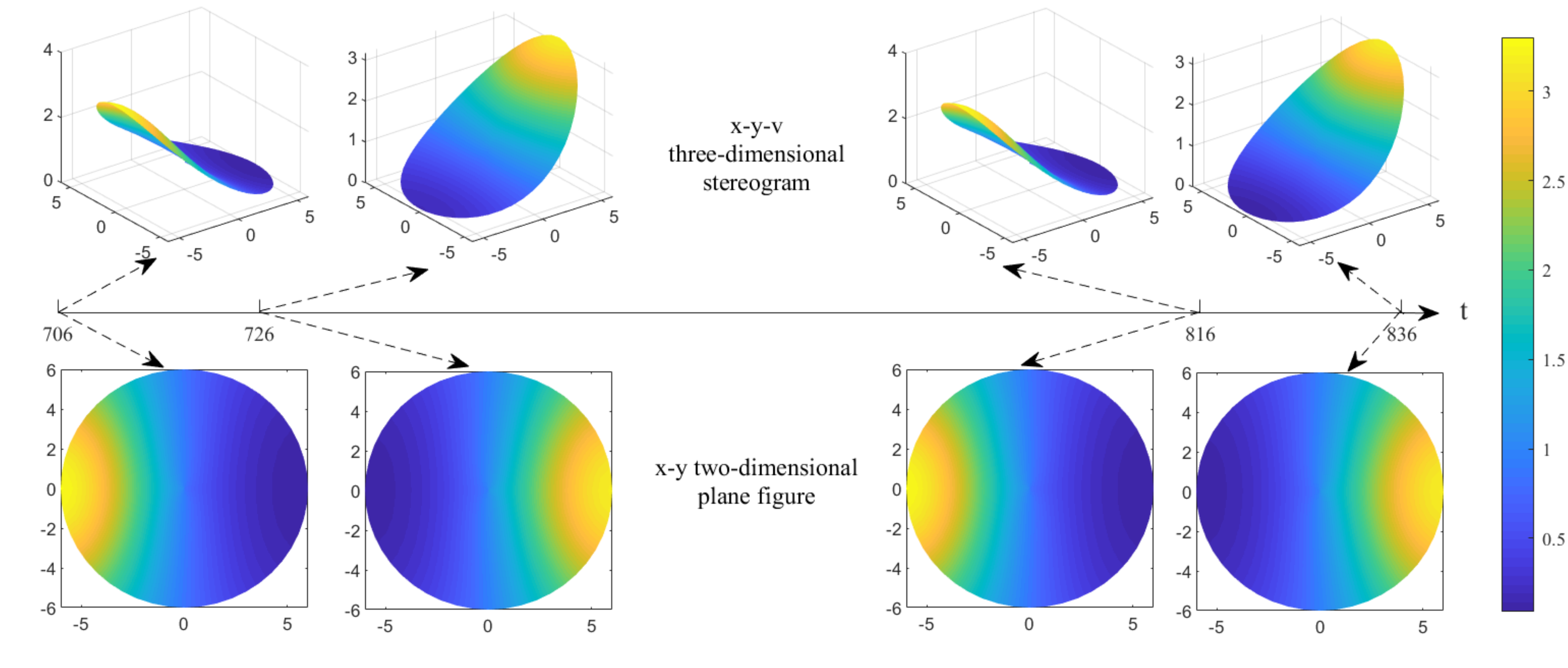}
\caption{The system produces standing waves with parameters: $b=0.25,~K=20,~a=0.3,~d=0.7,~e=0.5,\alpha=0.6,~d_1=0.3,~d_2=0.75,~R=6,~\tau=3$. Initial values are $u(t,r,\theta)=13.0320+0.01\cdot \cos t\cdot \cos r\cdot  \cos\theta,~v(t,r,\theta)=0.8108+0.01\cdot \cos t\cdot \cos r\cdot  \cos\theta,~t\in[-\tau,0)$. $(a):u,~(b):v$.}
\label{S-1}
\end{figure}

\begin{figure}[htbp]
\centering
(a)\includegraphics[width=1\textwidth]{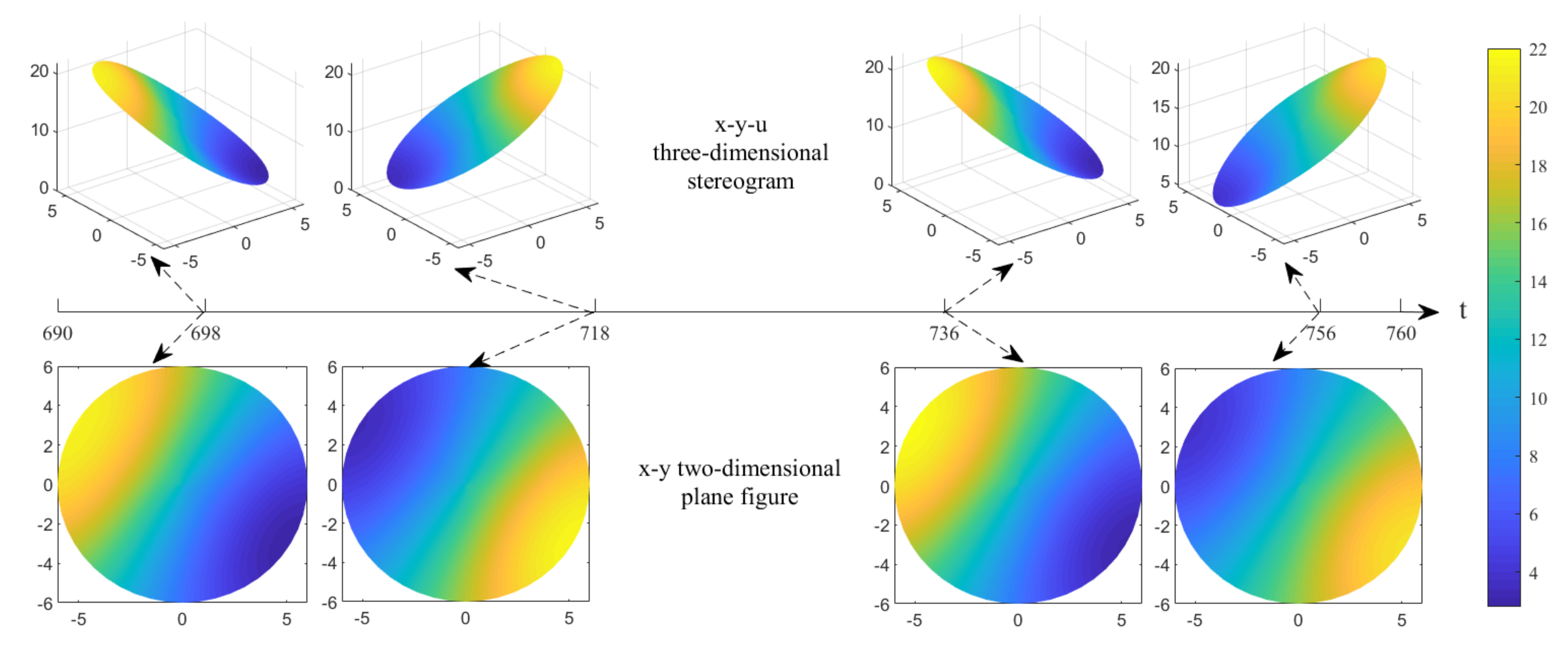}
(b)\includegraphics[width=1\textwidth]{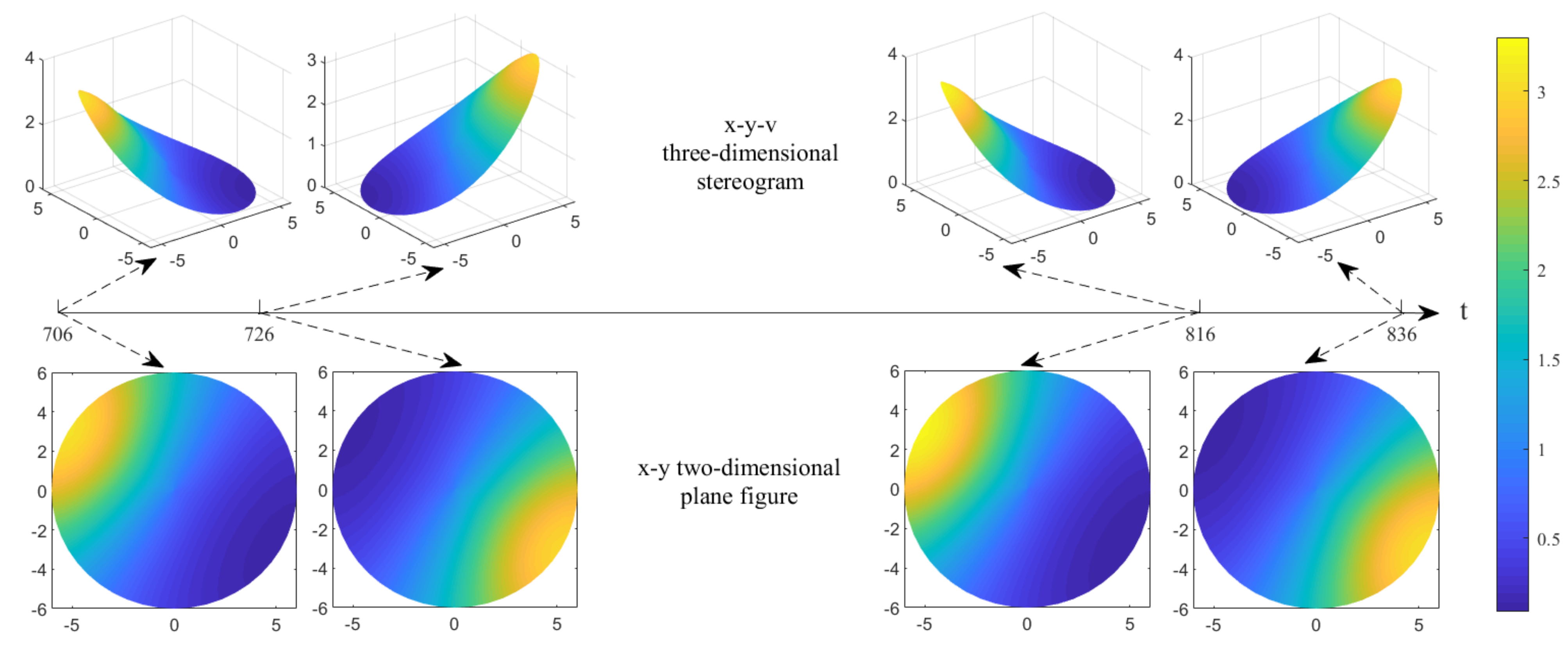}
\caption{The system produces standing waves with parameters: $b=0.25,~K=20,~a=0.3,~d=0.7,~e=0.5,\alpha=0.6,~d_1=0.3,~d_2=0.75,~R=6,~\tau=3$. Initial values are $u(t,r,\theta)=13.0320+0.01\cdot \cos t\cdot \cos r\cdot  \cos(\theta+\frac{\pi}{6}),~v(t,r,\theta)=0.8108+0.01\cdot \cos t\cdot \cos r\cdot  \cos(\theta+\frac{\pi}{6}),~t\in[-\tau,0)$. $(a):u,~(b):v$.}
\label{S-2}
\end{figure}

\begin{figure}[htbp]
\centering
(a)\includegraphics[width=1\textwidth]{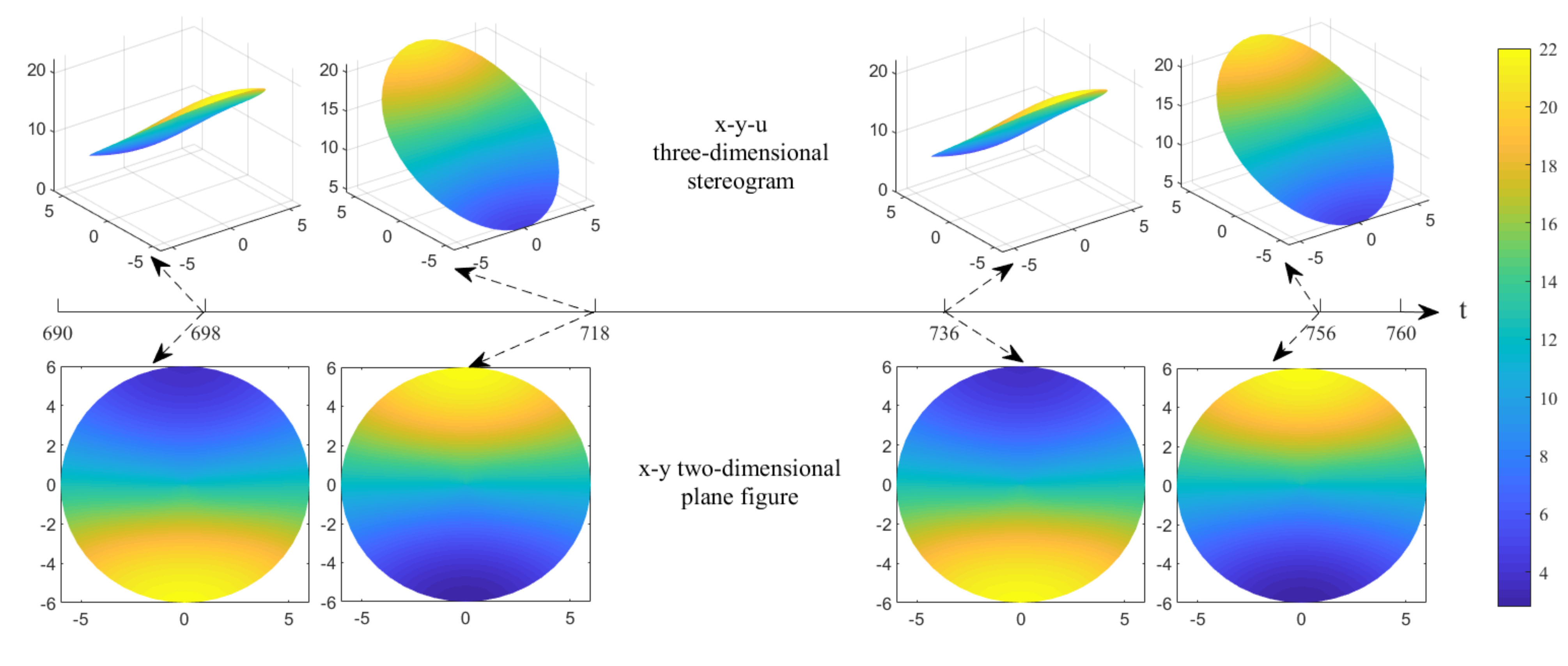}
(b)\includegraphics[width=1\textwidth]{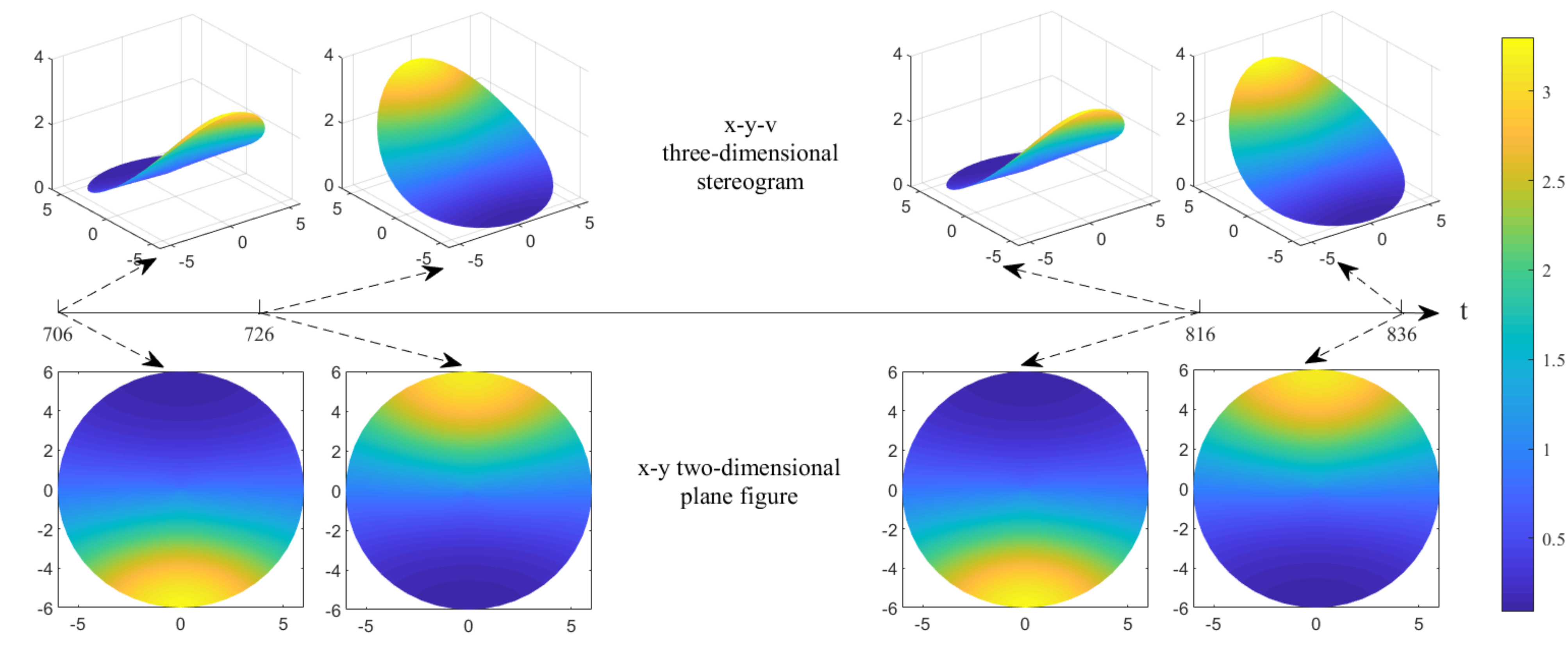}
\caption{The system produces standing waves with parameters: $b=0.25,~K=20,~a=0.3,~d=0.7,~e=0.5,\alpha=0.6,~d_1=0.3,~d_2=0.75,~R=6,~\tau=3$. Initial values are $u(t,r,\theta)=13.0320+0.01\cdot \cos t\cdot \cos r\cdot  \cos(\theta-\frac{\pi}{2}),~v(t,r,\theta)=0.8108+0.01\cdot \cos t\cdot \cos r\cdot  \cos(\theta-\frac{\pi}{2}),~t\in[-\tau,0)$. $(a):u,~(b):v$.}
\label{S-3}
\end{figure}

\begin{figure}[htbp]
\centering
(a)\includegraphics[width=1\textwidth]{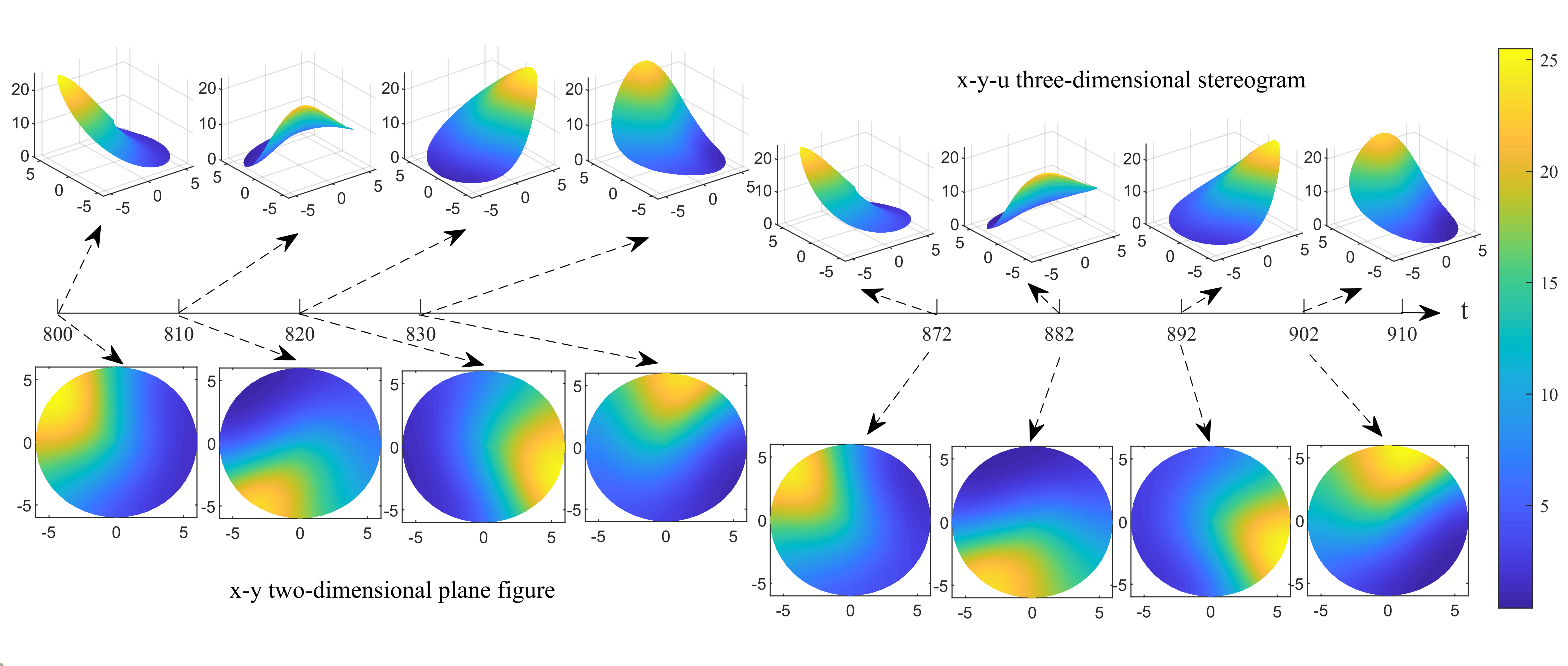}
(b)\includegraphics[width=1\textwidth]{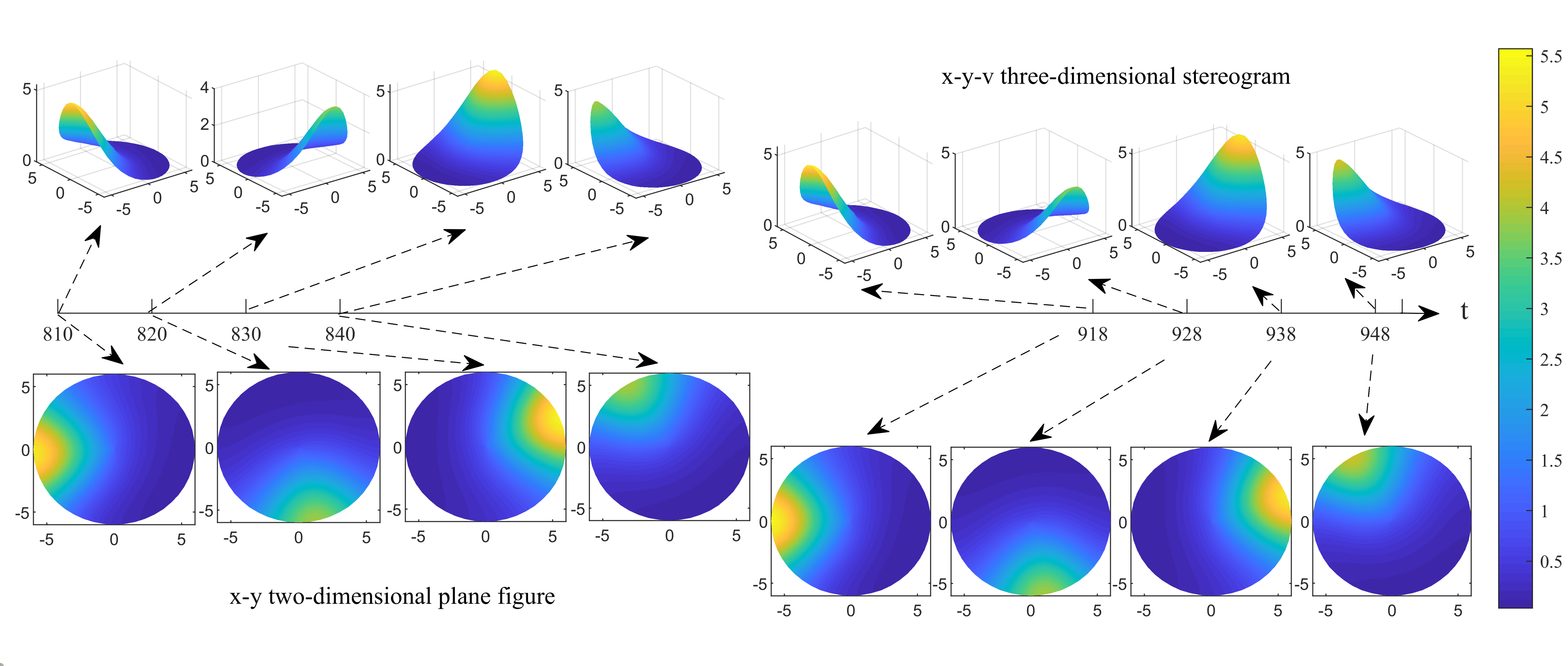}
\caption{The system produces rotating waves with parameters: $b=0.25,~K=20,~a=0.3,~d=0.7,~e=0.5,\alpha=0.6,~d_1=0.3,~d_2=0.75,~R=6,~\tau=3$. Initial values are $u(t,r,\theta)=13.0320+0.01\cdot \cos t\cdot \cos r\cdot  \sin\theta,~v(t,r,\theta)=0.8108+0.01\cdot \cos t\cdot \cos r\cdot  \cos\theta,~t\in[-\tau,0)$. $(a):u,~(b):v$.}
\label{T-1}
\end{figure}

\begin{figure}[htbp]
\centering
(a)\includegraphics[width=1\textwidth]{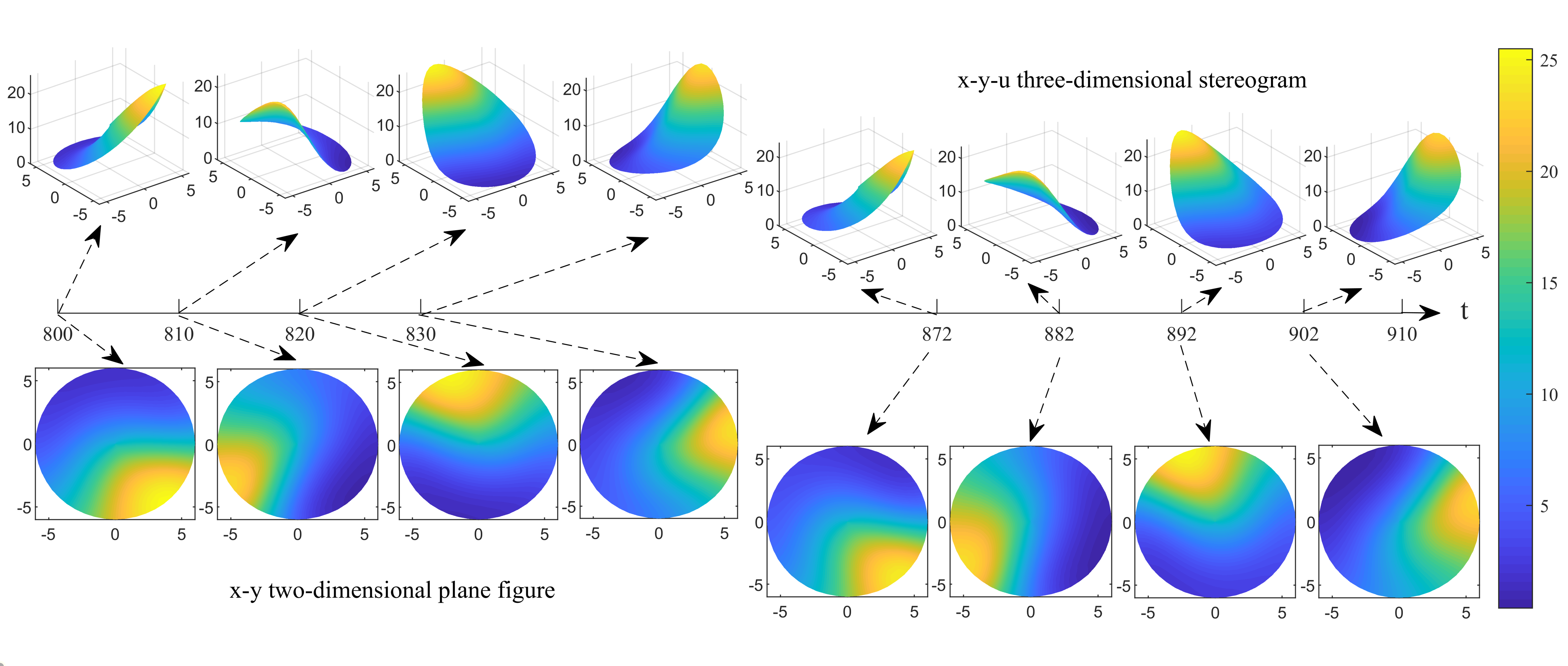}
(b)\includegraphics[width=1\textwidth]{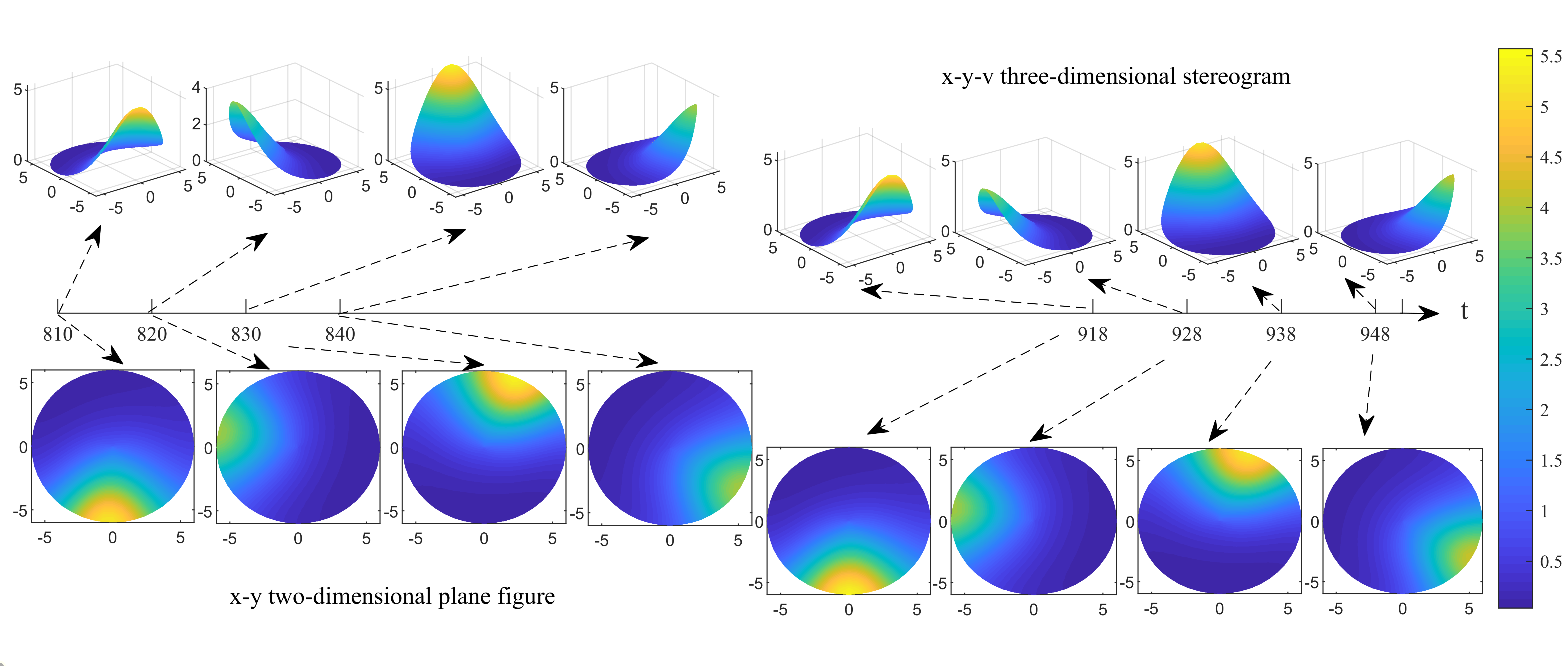}
\caption{The system produces counterpropagating waves with parameters: $b=0.25,~K=20,~a=0.3,~d=0.7,~e=0.5,\alpha=0.6,~d_1=0.3,~d_2=0.75,~R=6,~\tau=3$. Initial values are $u(t,r,\theta)=13.0320+0.01\cdot \cos t\cdot \cos r\cdot  \cos\theta,~v(t,r,\theta)=0.8108+0.01\cdot \cos t\cdot \cos r\cdot  \sin\theta,~t\in[-\tau,0)$. $(a):u,~(b):v$.}
\label{T-2}
\end{figure}

Fixing $b=0.25,~K=20,~a=0.3,~d=0.7,~e=0.5,~d_1=0.3,~d_2=0.75,~R=6$,~applying the same mathematical analysis method mentioned in \citep{Ruan2003J,Liu2020J}, at the unique positive constant steady solution, the first two bifurcation curves on the $\alpha-\tau$ plane  are shown in Figure \ref{nonlocalfenzhitu}. We select $\alpha=0.6$ on the plane of $\alpha-\tau$.
When $\hat{\lambda}=\lambda_{11}$, Hopf bifurcation occurs at $\hat{\tau}~=\tau_{\lambda_{11}}^0 \approx 1.7825$.
We know that when $\tau ~\textless~1.7825$, $E^*$ is locally asymptotically stable,
and when $\tau~\textgreater~1.7825$, $E^*$ is unstable. The bifurcation generated at this time is an equivariant hopf bifurcation.
It can be obtained through numerical calculation that $\mu=1.2175,~B_{11}\approx 0.0021-0.0911\mathrm{i},~B_{2001}\approx -0.1075+0.0745\mathrm{i},~B_{1110}\approx -0.1813+0.1620\mathrm{i}$.
Thus, $a_1\mu \approx 0.0026,~a_2 \approx -0.1075,~a_2 +a_3 \approx -0.2888,~a_2 -a_3 \approx 0.0738$, which corresponds to Case 2 when $a_1\mu>0$ in Table \ref{tab2}.
By Theorem \ref{rotating and standing}, we know that system possesses an unstable standing wave (see Figure \ref{S-1}-\ref{S-3}) and two stable rotating waves (see Figure \ref{T-1}-\ref{T-2}).
\begin{remark}
We can see from Figure \ref{nonlocalfenzhitu} that as $\alpha$ changes, a double Hopf point HH appears. Below the lower line is the stable region of the system where $\tau~\textless~ min\left\{\tau_{\lambda_{00}^0},\tau_{\lambda_{11}^{0}}\right\}$. And above the lower line where $\tau~\textgreater~ \left\{\tau_{\lambda_{00}^0},\tau_{\lambda_{11}^{0}}\right\}$ the system may produce spatially homogeneous or inhomogeneous period solutions. Investigating the detailed bifurcation sets might require studying an at least six dimensional center manifold.
\end{remark}

\begin{remark}
We can find that only when the initial value restricted to the center subspace satisfies $\rho_1=\rho_2$, the spatially inhomogeneous periodic solution is in the form of standing waves. For example, we select the initial value as $u(t,r,\theta)=u^*+\varsigma_1(t,r) \cdot \cos(\theta+\hat{\theta}),~t \in [-\tau,0);~v(t,r,\theta)=v^{*}+\varsigma_2(t,r)  \cdot \cos(\theta+\hat{\theta}),~t \in [-\tau,0)$, which has the following approximate form restricted to the center manifold
$$
U_0(\vartheta)(r,\theta) \approx 4(\mathrm{Re}\{\varsigma_1(\vartheta,r)\},\mathrm{Re}\{p_0\cdot\varsigma_2(\vartheta,r)\})^{\mathrm{T}}\cos(\theta+\hat{\theta}),
$$
with $z_1=z_2=z_3=z_4=1$. Thus, the spatially inhomogeneous periodic solutions in Figure \ref{S-1}-\ref{S-3} are in the form of standing waves.
No matter what value of $\hat{\theta}$ is taken, the simulation is a standing wave solution, which reflects the effect of  $O(2)$ equivariance.
However, when the initial values of $u$ and $v$ are chosen with other forms,  solutions of the system are attracted by one of two coexisting stable rotating waves (see Figure \ref{T-1}-\ref{T-2}), which may be clockwise (Figure \ref{T-2}) or counterclockwise (Figure \ref{T-1}).
\end{remark}


\newpage
\bibliography{Hopf}
\newpage
\appendix
\section{The specific calculation of $g_3^1(z,0,\mu)$ }\label{ Appendix A}
\subsection{The calculation of $\mathrm{Proj}_{\mathrm{Ker}(M_3^1)} f_3^1(z,0,0).$}

Writing $\tilde{F}_3(\Phi_{r\theta} z,\mu)$ as follows
$$
\tilde{F}_3(\Phi_{r\theta} z,0)=\sum_{p_1+p_2+p_3+p_4=3} A_{p_1p_2p_3p_4} \left(\hat{\phi}_{nm}^c\right)^{p_1+p_2}  \left({\hat{\phi}_{nm}^s}\right)^{p_3+p_4} z_1^{p_1}z_2^{p_2}z_3^{p_3}z_4^{p_4},
$$
then we have
$$
\begin{aligned}
f_3^1(z,0,0)&=\left(\begin{array}{cc}
\left\langle \tilde{F}_3(\Phi_{r\theta} z,0),\Psi_{r\theta}^1(0) \right\rangle \\
\left\langle \tilde{F}_3(\Phi_{r\theta} z,0),\Psi_{r\theta}^2(0) \right\rangle
\end{array}\right)\\
&=
\overline{\Psi(0)}
\left(\begin{array}{cc}
\sum_{p_1+p_2+p_3+p_4=3} A_{p_1p_2p_3p_4} \int_0^R \int_0^{2\pi} r \left(\hat{\phi}_{nm}^c\right)^{p_1+p_2}  \left({\hat{\phi}_{nm}^s}\right)^{p_3+p_4+1}\mathrm{d}\theta \mathrm{d}r z_1^{p_1}z_2^{p_2}z_3^{p_3}z_4^{p_4}\\
\sum_{p_1+p_2+p_3+p_4=3} A_{p_1p_2p_3p_4}  \int_0^R \int_0^{2\pi} r \left(\hat{\phi}_{nm}^c\right)^{p_1+p_2+1}  \left({\hat{\phi}_{nm}^s}\right)^{p_3+p_4} \mathrm{d}\theta \mathrm{d}r z_1^{p_1}z_2^{p_2}z_3^{p_3}z_4^{p_4}
\end{array}
\right).
\end{aligned}
$$

Noticing the fact
$$
\int_0^R\int_0^{2\pi} r \left(\hat{\phi}_{nm}^c\right)^{4} \mathrm{d} \theta \mathrm{d} r = \int_0^R\int_0^{2\pi} r \left(\hat{\phi}_{nm}^s\right)^{4} \mathrm{d} \theta \mathrm{d} r=0 ,
$$
$$
\int_0^R\int_0^{2\pi} r \left(\hat{\phi}_{nm}^c\right)^{3} \hat{\phi}_{nm}^s \mathrm{d} \theta \mathrm{d} r= \int_0^R\int_0^{2\pi} r \hat{\phi}_{nm}^c \left(\hat{\phi}_{nm}^s\right)^{3} \mathrm{d} \theta \mathrm{d} r=0,
$$
$$
\int_0^R\int_0^{2\pi} r \left(\hat{\phi}_{nm}^c\right)^{2} \left(\hat{\phi}_{nm}^s\right)^{2} \mathrm{d} \theta \mathrm{d} r \triangleq \mathrm{M}_{22} ,
$$
and the relationship of $\Phi_{r\theta}$ and $\Psi_{r\theta}$, we get (\ref{Part1}).

\subsection{The calculation of $\mathrm{Proj}_{\mathrm{Ker}(M_3^1)} \left(D_z f_2^1(z,0,0)U_2^1(z,0)\right).$}

We have
\begin{equation}\label{f21}
\begin{aligned}
f_2^1(z,0,0)&=\left(\begin{array}{cc}
\left\langle \tilde{F}_2(\Phi_{r\theta} z,0),\Psi_{r\theta}^1(0) \right\rangle \\
\left\langle \tilde{F}_2(\Phi_{r\theta} z,0),\Psi_{r\theta}^2(0)\right\rangle
\end{array}\right)\\
&=
\overline{\Psi(0)}
\left(\begin{array}{cc}
\sum_{p_1+p_2+p_3+p_4=2} A_{p_1p_2p_3p_4} \int_0^R \int_0^{2\pi} r \left(\hat{\phi}_{nm}^c\right)^{p_1+p_2}  \left({\hat{\phi}_{nm}^s}\right)^{p_3+p_4+1}\mathrm{d}\theta \mathrm{d}r z_1^{p_1}z_2^{p_2}z_3^{p_3}z_4^{p_4}\\
\sum_{p_1+p_2+p_3+p_4=2} A_{p_1p_2p_3p_4} \int_0^R \int_0^{2\pi} r \left(\hat{\phi}_{nm}^c\right)^{p_1+p_2+1}  \left({\hat{\phi}_{nm}^s}\right)^{p_3+p_4} \mathrm{d}\theta \mathrm{d}r z_1^{p_1}z_2^{p_2}z_3^{p_3}z_4^{p_4}
\end{array}
\right).
\end{aligned}
\end{equation}
Noticing the fact that
$$
\int_0^R\int_0^{2\pi} r \left(\hat{\phi}_{nm}^c\right)^{3} \mathrm{d} \theta \mathrm{d} r= \int_0^R\int_0^{2\pi} r \left(\hat{\phi}_{nm}^s\right)^{3} \mathrm{d} \theta \mathrm{d} r=0,
$$
$$
\int_0^R\int_0^{2\pi} r \left(\hat{\phi}_{nm}^c\right)^{2} \hat{\phi}_{nm}^s \mathrm{d} \theta \mathrm{d} r
=\int_0^R\int_0^{2\pi} r \hat{\phi}_{nm}^c \left(\hat{\phi}_{nm}^s\right)^2 \mathrm{d} \theta \mathrm{d} r=0 ,
$$
then we get $\mathrm{Proj}_{\mathrm{Ker}(M_3^1)} \left(D_z f_2^1(z,0,0)U_2^1(z,0)\right)=0$.

\subsection{The calculation of $\mathrm{Proj}_{\mathrm{Ker}(M_3^1)} \left(D_y f_2^1(z,0,0)U_2^2(z,0)\right).$}

Firstly, we calculate the $\rm{Fr\acute{e}chet}$ derivative $D_yf_2^1(z,0,0):Q_s\rightarrow \mathscr{X}_\mathbb{C}$.
By (\ref{F2z}) and (\ref{F2phi}),
$\tilde{F}_2(z,y,0)$ can be written as
\begin{equation}\label{tildeF2}
\begin{aligned}
\tilde{F}_2(z,y,0)&=S_2(\Phi_{r\theta} z,y)+o(z^2,y^2)\\
&=S_{yz_1}(y)z_1\hat{\phi}_{nm}^c+S_{yz_2}(y)z_2\hat{\phi}_{nm}^c+S_{yz_3}(y)z_3\hat{\phi}_{nm}^s+S_{yz_4}(y)z_4\hat{\phi}_{nm}^s+o(z^2,y^2),
\end{aligned}
\end{equation}
where $S_{yz_k}(k=1,2,3,4):Q_s\rightarrow \mathscr{X}_\mathbb{C}$ are linear operators, and
$$
S_{yz_k}(\varphi)=S_{y(0)z_k}(\varphi(0))+S_{y(-1)z_k}(\varphi(-1)).
$$

Let
$$
U_2^2(z,0)\triangleq h(z)=\sum_{j\ge 0}h_j(z)\hat{\phi}_{jk}(r,\theta),
$$
with
$$
h_{jk}(z)=\left(\begin{array}{cc}
h_{jk}^{(1)}(z)\\
h_{jk}^{(2)}(z)
\end{array}\right)=\sum_{p_1+p_2+p_3+p_4=2}\left(\begin{array}{cc}
h_{jkp_1p_2p_3p_4}^{(1)}(z)\\
h_{jkp_1p_2p_3p_4}^{(2)}(z)
\end{array}\right)z_1^{p_1}z_2^{p_2}z_3^{p_3}z_4^{p_4}.
$$
Therefore,
$$
\begin{aligned}
D_y\tilde{F}_2(z,0,0)\left(U_2^2(z,0)\right)&=\left(\begin{array}{cc}
\langle D_y\tilde{F}_2(z,0,0)\left(U_2^2(z,0)\right),\Psi_{r\theta}^1(0) \rangle\\
\langle D_y\tilde{F}_2(z,0,0)\left(U_2^2(z,0)\right),\Psi_{r\theta}^2(0) \rangle
\end{array}\right)\\
&=
\overline{\Psi(0)}
\left(\begin{array}{cc}
\sum_{j\ge 0}\left[\mathrm{M}_{jkcs}S_{yz_1}(h_{jk})z_1+\mathrm{M}_{jkcs}S_{yz_2}(h_{jk})z_2\right.\\
\left.+\mathrm{M}_{jkss}S_{yz_3}(h_{jk})z_3+\mathrm{M}_{jkss}S_{yz_4}(h_{jk})z_4\right]\\
\sum_{j\ge 0}[\mathrm{M}_{jkcc}S_{yz_1}(h_{jk})z_1+\mathrm{M}_{jkcc}S_{yz_2}(h_{jk})z_2\\
+\mathrm{M}_{jksc}S_{yz_3}(h_{jk})z_3+\mathrm{M}_{jksc}S_{yz_4}(h_{jk})z_4]
\end{array}
\right),
\end{aligned}
$$
where
$$
\mathrm{M}_{jksc}=\mathrm{M}_{jkcs}=\int_0^R \int_0^{2\pi} r \hat{\phi}_{jk}\hat{\phi}_{nm}^c\hat{\phi}_{nm}^s \mathrm{d}\theta\mathrm{d}r=\left\{\begin{array}{cccc}
\mathrm{M}_{0kcs}^c & \hat{\phi}_{0}~or~\hat{\phi}_{0k}^c,\\
0, & otherwise,
\end{array}\right.
$$
$$
\mathrm{M}_{jkss}=\int_0^R \int_0^{2\pi} r \hat{\phi}_{jk}\hat{\phi}_{nm}^s\hat{\phi}_{nm}^s \mathrm{d}\theta\mathrm{d}r=\left\{\begin{array}{cccc}
\mathrm{M}_{2nkss}^c,  & \hat{\phi}_{jk}=\hat{\phi}_{2nk}^c,\\
0, & otherwise,
\end{array}\right.
$$
$$
\mathrm{M}_{jkcc}=\int_0^R \int_0^{2\pi} r \hat{\phi}_{jk}\hat{\phi}_{nm}^c\hat{\phi}_{nm}^c \mathrm{d}\theta\mathrm{d}r=\left\{\begin{array}{cccc}
\mathrm{M}_{2nkcc}^s,  & \hat{\phi}_{jk}=\hat{\phi}_{2nk}^s,\\
0, & otherwise.
\end{array}\right.
$$

Moreover, we have
$$
\begin{aligned}
D_y\tilde{F}_2(z,0,0)\left(U_2^2(z,0)\right)
&=\bar{\Psi}(0)\left(\begin{array}{cccc}
N_1\\
N_2
\end{array}
\right),
\end{aligned}
$$
with
$$
\begin{aligned}
N_1=&\mathrm{M}_{0kcs}^c\left(S_{yz_1}(h_{0k}^{ccs})z_1+S_{yz_2}(h_{0k}^{ccs})z_2\right)
+\mathrm{M}_{2nkss}^c\left(S_{yz_3}(h_{2nk}^{css})z_3+S_{yz_4}(h_{2nk}^{css})z_4\right),\\
N_2=&\mathrm{M}_{2nkcc}^s\left(S_{yz_1}(h_{2nk}^{ccs})z_1+S_{yz_2}(h_{2nk}^{ccs})z_2\right)
+\mathrm{M}_{0kcs}^c\left(S_{yz_3}(h_{0k}^{ccs})z_3+S_{yz_4}(h_{0k}^{ccs})z_4\right).\\
\end{aligned}
$$
Thus,
\begin{equation}\label{part3}
\begin{aligned}
&\frac{1}{3!}\mathrm{Proj}_{\mathrm{Ker(M_3^1)}}\left(D_yf_2^1(z,0,0)U_2^2(z,0)\right)\\
=&\left(
\begin{array}{cccc}
E_{2100}^1z_1^2z_2+E_{2001}^1z_1^2z_4+E_{0120}^1z_3^2z_2+E_{0021}^1z_3^2z_4+E_{1110}^1z_1z_2z_3+E_{1011}^1z_1z_3z_4\\
\overline{E_{2100}^1}z_1z_2^2+\overline{E_{2001}^1}z_2^2z_3+\overline{E_{0120}^1}z_4^2z_1+\overline{E_{0021}^1}z_3z_4^2+\overline{E_{1110}^1}z_1z_2z_4+\overline{E_{1011}^1}z_2z_3z_4\\
E_{2100}^2z_1^2z_2+E_{2001}^2z_1^2z_4+E_{0120}^2z_3^2z_2+E_{0021}^2z_3^2z_4+E_{1110}^2z_1z_2z_3+E_{1011}^2z_1z_3z_4\\
\overline{E_{2100}^2}z_1z_2^2+\overline{E_{2001}^2}z_2^2z_3+\overline{E_{0120}^2}z_4^2z_1+\overline{E_{0021}^2}z_3z_4^2+\overline{E_{1110}^2}z_1z_2z_4+\overline{E_{1011}^2}z_2z_3z_4
\end{array}
\right),
\end{aligned}
\end{equation}
where
$$
\begin{aligned}
E_{2100}^1=\frac{1}{6}\overline{\Psi_1(0)}&\left[\mathrm{M}_{0kcs}^c\left(S_{yz_1}(h_{0k1100}^{ccs})+S_{yz_2}(h_{0k2000}^{ccs})\right)\right],\\
E_{2001}^1=\frac{1}{6}\overline{\Psi_1(0)}&\left[\mathrm{M}_{0kcs}^cS_{yz_1}(h_{0k1001}^{ccs})+\mathrm{M}_{2nkss}^cS_{yz_4}(h_{2nk2000}^{css})\right],\\
E_{0120}^1=\frac{1}{6}\overline{\Psi_1(0)}&\left[\mathrm{M}_{0kcs}^cS_{yz_2}(h_{0k0020}^{ccs}))+\mathrm{M}_{2nkss}^cS_{yz_3}(h_{2nk0110}^{css})\right],\\
E_{0021}^1=\frac{1}{6}\overline{\Psi_1(0)}&\left[\mathrm{M}_{2nkss}^c\left(S_{yz_3}(h_{2nk0011}^{css})+S_{yz_4}(h_{2nk0020}^{css})\right)\right],\\
E_{1110}^1=\frac{1}{6}\overline{\Psi_1(0)}&\left[\mathrm{M}_{0kcs}^c\left(S_{yz_1}(h_{0k0110}^{ccs})+S_{yz_2}(h_{0k1010}^{ccs})\right)+\mathrm{M}_{2nkss}^cS_{yz_3}(h_{2nk1100}^{css})\right],\\
E_{1011}^1=\frac{1}{6}\overline{\Psi_1(0)}&\left[\mathrm{M}_{0kcs}^cS_{yz_1}(h_{0k0011}^{ccs})
+\mathrm{M}_{2nkss}^c\left(S_{yz_3}(h_{2nk1001}^{css})+S_{yz_4}(h_{2nk1010}^{css})\right)\right],\\
E_{2100}^2=\frac{1}{6}\overline{\Psi_3(0)}&\left[\mathrm{M}_{2nkcc}^s\left(S_{yz_1}(h_{2nk1100}^{ccs})+S_{yz_2}(h_{2nk2000}^{ccs})\right)\right],\\
E_{2001}^2=\frac{1}{6}\overline{\Psi_3(0)}&\left[\mathrm{M}_{2nkcc}^sS_{yz_1}(h_{2nk1001}^{ccs})+\mathrm{M}_{0kcs}^cS_{yz_2}(h_{0k2000}^{ccs})\right],\\
E_{0120}^2=\frac{1}{6}\overline{\Psi_3(0)}&\left[\mathrm{M}_{2nkcc}^sS_{yz_2}(h_{2nk0020}^{ccs})+\mathrm{M}_{0kcs}^cS_{yz_3}(h_{0k0110}^{ccs})\right],\\
E_{0021}^2=\frac{1}{6}\overline{\Psi_3(0)}&\left[\mathrm{M}_{0kcs}^c\left(S_{yz_3}(h_{0k0011}^{ccs})+S_{yz_4}(h_{0k0020}^{ccs})\right)\right],\\
E_{1110}^2=\frac{1}{6}\overline{\Psi_3(0)}&\left[\mathrm{M}_{2nkcc}^s\left(S_{yz_1}(h_{2nk0110}^{ccs})+S_{yz_2}(h_{2nk1010}^{ccs})\right)+\mathrm{M}_{0kcs}^cS_{yz_3}(h_{0k1100}^{ccs})\right]\\
E_{1011}^2=\frac{1}{6}\overline{\Psi_3(0)}&\left[\mathrm{M}_{2nkcc}^sS_{yz_1}(h_{2nk0011}^{ccs})+\mathrm{M}_{0kcs}^c\left(S_{yz_3}(h_{0k1001}^{ccs})+S_{yz_4}(h_{0k1010}^{ccs})\right)\right].
\end{aligned}
$$

Now, we need to calculate
$$
\begin{aligned}
&h_{0k2000}^{ccs},~h_{0k1100}^{ccs},~h_{0k1010}^{ccs},~h_{0k1001}^{ccs},~h_{0k0110}^{ccs},~h_{0k0020}^{ccs},~h_{0k0011}^{ccs},\\
&h_{2nk2000}^{ccs},~h_{2nk1100}^{ccs},~h_{2nk1010}^{ccs},~h_{2nk1001}^{ccs},~h_{2nk0110}^{ccs},~h_{2nk0020}^{ccs},~h_{2nk0011}^{ccs},\\
&h_{2nk2000}^{css},~h_{2nk1100}^{css},~h_{2nk1010}^{css},~h_{2nk1001}^{css},~h_{2nk0110}^{css},~h_{2nk0020}^{css},~h_{2nk0011}^{css}.
\end{aligned}
$$

From (\ref{zt2}) and(\ref{Mj1z}), we get
$$
\begin{aligned}
M_2^2U_2^2(z,0)(\vartheta)&=M_2^2h(z)(\vartheta)\\
&=\left\{\begin{array}{cc}
D_zh(z)Bz-\tilde{D}_0\Delta h(0)-\tilde{L}_0(h(z)), & \vartheta=0,\\
D_zh(z)Bz-D_{\vartheta}h(z), & \vartheta \neq 0,\\
\end{array}\right.\\
&=\left\{\begin{array}{cc}
\sum_{j \ge 0}\left[D_zh_j(z)\hat{\phi}_{jk}(r,\theta)Bz-\tilde{D}_0\Delta h_j(z)\hat{\phi}_{jk}(r,\theta)-\tilde{L}_0(h_j(z)\hat{\phi}_{jk}(r,\theta))\right], & \vartheta=0,\\
\sum_{j \ge 0}\left[D_zh_j(z)\hat{\phi}_{jk}(r,\theta)Bz-D_{\vartheta}h_j(z)\hat{\phi}_{jk}(r,\theta)\right], & \vartheta \neq 0,\\
\end{array}\right.
\end{aligned}
$$
and
$$
f_2^2(z,0,0)=\left\{\begin{array}{cc}\begin{aligned}
\tilde{F}_2&(z,0,0)-\Phi_1(0)f_2^{1(1)}(z,0,0)\hat{\phi}_{nm}^s-\Phi_2(0)f_2^{1(2)}(z,0,0)\hat{\phi}_{nm}^s\\
&-\Phi_3(0)f_2^{1(1)}(z,0,0)\hat{\phi}_{nm}^c-\Phi_4(0)f_2^{1(1)}(z,0,0)\hat{\phi}_{nm}^c, & \vartheta=0,\\
-\Phi_1&(\vartheta)f_2^{1(1)}(z,0,0)\hat{\phi}_{nm}^s-\Phi_2(\vartheta)f_2^{1(2)}(z,0,0)\hat{\phi}_{nm}^s\\
&-\Phi_3(\vartheta)f_2^{1(1)}(z,0,0)\hat{\phi}_{nm}^c-\Phi_4(\vartheta)f_2^{1(1)}(z,0,0)\hat{\phi}_{nm}^c, & \vartheta \neq 0.
\end{aligned}\end{array}\right.
$$
Besides, we have
\begin{equation}\label{M22U22}
 \langle M_2^2\left(U_2^2(z,0)\right),\beta_{jk} \rangle= \langle f_2^2(z,0,0),\beta_{jk} \rangle,
\end{equation}
with $\beta_{jk}=\frac{\phi_{jk}}{\|\phi_{jk}\|}$.

Thus, the expressions of $h_{jp_1p_2p_3p_4}$ can be obtained. Due to the large number of expressions, we show the specific results in the Appendix F.

Noting the fact that
$$
\mathrm{M}_{2nkcc}^s=\mathrm{M}_{2nkss}^c,
$$
therefore, we have
$$
\begin{aligned}
&E_{2100}^1=E_{0021}^2,~E_{2001}^1=E_{0120}^2,~E_{0120}^1=E_{2001}^2,\\
&E_{0021}^1=E_{2100}^2,~E_{1110}^1=E_{1011}^2,~E_{1011}^1=E_{1110}^2,
\end{aligned}
$$
For simplification of notations, we rewrite (\ref{part3}) as (\ref{Part3}).

\section{Proof of Lemma \ref{reduce}}\label{Proof of Lemma}
By a smooth transformation
\begin{equation}\label{smooth transformation}
\begin{aligned}
&z_1=\zeta_1+b_1\zeta_1^2\zeta_2+b_2\zeta_3^2\zeta_2+b_3\zeta_3^2\zeta_4+b_4\zeta_1\zeta_3\zeta_4,\\
&z_2=\bar{z}_1,\\
&z_3=\zeta_1+b_1\zeta_3^2\zeta_4+b_2\zeta_1^2\zeta_4+b_3\zeta_1^2\zeta_2+b_4\zeta_1\zeta_2\zeta_3,\\
&z_4=\bar{z}_3,
\end{aligned}
\end{equation}
we have
\begin{equation}\label{zeta}
\begin{aligned}
&\zeta_1=z_1-b_1z_1^2z_2-b_2z_3^2z_2-b_3z_3^2z_4-b_4z_1z_3z_4+o(4),\\
&\zeta_2=\bar{\zeta}_1,\\
&\zeta_3=z_1-b_1z_3^2z_4-b_2z_1^2z_4-b_3z_1^2z_2-b_4z_1z_2z_3+o(4),\\
&\zeta_4=\bar{\zeta}_3.
\end{aligned}
\end{equation}
Then
$$
\begin{aligned}
\dot{\zeta}_1=&\dot{z}_1-2b_1z_1z_2\dot{z}_1-b_1z_1^2\dot{z}_2-2b_2z_2z_3\dot{z}_3-b_2z_3^2\dot{z}_2-2b_3z_3z_4\dot{z}_3-b_3z_3^2\dot{z}_4
-b_4z_1z_3\dot{z}_4-b_4z_1z_4\dot{z}_3-b_4z_3z_4\dot{z}_1+o(4)\\
=&(\mathrm{i}\omega_{\hat{\lambda}}+B_{11}\mu)z_1+B_{2100}z_1^2z_2^2+B_{1011}z_1z_3z_4+B_{2001}z_1^2z_4+B_{0120}z_3^2z_2+B_{0021}z_3^2z_4+B_{1110}z_1z_2z_3\\
&-3b_1(\mathrm{i}\omega_{\hat{\lambda}}+B_{11}\mu)z_1^2z_2-3b_2(\mathrm{i}\omega_{\hat{\lambda}}+B_{11}\mu)z_3^2z_2-3b_3(\mathrm{i}\omega_{\hat{\lambda}}+B_{11}\mu)z_3^2z_4
-3b_4(\mathrm{i}\omega_{\hat{\lambda}}+B_{11}\mu)z_1z_3z_4+o(4).\\
\end{aligned}
$$
Let
$$
\begin{aligned}
&b_1=\frac{B_{2100}}{3(\mathrm{i}\omega_{\hat{\lambda}}+B_{11}\mu)},~b_2=\frac{B_{0120}}{3(\mathrm{i}\omega_{\hat{\lambda}}+B_{11}\mu)},\\
&b_3=\frac{B_{0021}}{3(\mathrm{i}\omega_{\hat{\lambda}}+B_{11}\mu)},~b_4=\frac{B_{1011}}{3(\mathrm{i}\omega_{\hat{\lambda}}+B_{11}\mu)},
\end{aligned}
$$
then
$$
\begin{aligned}
\dot{\zeta}_1=&(\mathrm{i}\omega_{\hat{\lambda}}+B_{11}\mu)z_1+B_{2001}z_1^2z_4^2+B_{1110}z_1z_2z_3+o(4)\\
=&(\mathrm{i}\omega_{\hat{\lambda}}+B_{11}\mu)\zeta_1+B_{2001}\zeta_1^2\zeta_4^2+B_{1110}\zeta_1\zeta_2\zeta_3+o(4).
\end{aligned}
$$
The same is true for $\dot{\zeta}_2, \dot{\zeta}_3$ and $\dot{\zeta}_4$ so that (\ref{z1z2z3z4}) is established.

\section{Proof of Theorem \ref{rotating and standing}}\label{Proof of Theorem}

We only need to prove the approximate expressions of rotating and standing wave solutions reduced to the center subspace, and the rest of the theorem can be easily obtained from previous analysis.

By (\ref{basis of P}), (\ref{fenjie}),(\ref{chi}), we get
$$
\begin{aligned}
U_t(\vartheta)(r,\theta) \approx &\xi \mathrm{e}^{\mathrm{i}\omega_{\hat{\lambda}}\vartheta} J_n(\sqrt{\lambda_{nm}}r)\mathrm{e}^{\mathrm{i}n\theta} \rho_1\mathrm{e}^{\mathrm{i}\chi_1(t)}
+\bar{\xi} \mathrm{e}^{-\mathrm{i}\omega_{\hat{\lambda}}\vartheta} J_n(\sqrt{\lambda_{nm}}r)\mathrm{e}^{\mathrm{i}n\theta} \rho_2\mathrm{e}^{-\mathrm{i}\chi_2(t)}\\
&+\xi \mathrm{e}^{\mathrm{i}\omega_{\hat{\lambda}}\vartheta} J_n(\sqrt{\lambda_{nm}}r)\mathrm{e}^{-\mathrm{i}n\theta} \rho_2\mathrm{e}^{\mathrm{i}\chi_2(t)}
+\bar{\xi} \mathrm{e}^{-\mathrm{i}\omega_{\hat{\lambda}}\vartheta} J_n(\sqrt{\lambda_{nm}}r)\mathrm{e}^{-\mathrm{i}n\theta} \rho_1\mathrm{e}^{-\mathrm{i}\chi_1(t)}
\end{aligned}
$$
with $\xi=(p_{11},p_{12},\cdots,p_{1n})^{\mathrm{T}}$.
For simplicity, we also rewrite $p_{1i}$ in the form of a complex angle as $p_{1i}=|p_{1i}|\mathrm{e}^{\mathrm{i} \mathrm{Arg}(p_{1i})}$ in the subsequent calculations.

For $(\rho_1,\rho_2)=(0,\sqrt{\frac{-a_1\mu}{a_2}})$,
$$
\begin{aligned}
U_t(\vartheta)(r,\theta) \approx &\xi \mathrm{e}^{\mathrm{i}\omega_{\hat{\lambda}}\vartheta} J_n(\sqrt{\lambda_{nm}}r)\mathrm{e}^{\mathrm{i}n\theta} \rho_1\mathrm{e}^{\mathrm{i}\chi_1(t)}
+\bar{\xi} \mathrm{e}^{-\mathrm{i}\omega_{\hat{\lambda}}\vartheta} J_n(\sqrt{\lambda_{nm}}r)\mathrm{e}^{-\mathrm{i}n\theta} \rho_1\mathrm{e}^{-\mathrm{i}\chi_1(t)}\\
 \approx &
\sum_{i=1}^n{2|p_{1i}|\sqrt{\frac{-a_1\mu}{a_2}} J_n(\sqrt{\lambda_{nm}}r)\cos(\mathrm{Arg}(p_{1i})+\omega_{\hat{\lambda}}\vartheta+\omega_{\hat{\lambda}}t+n\theta) {e}_i}.
\end{aligned}
$$
where $e_i$ is the $i$th unit coordinate vector of $\mathbb{R}^n$. This corresponds to the form of a rotating wave solution in the plane of $(z_2,z_3)$.

For $(\rho_1,\rho_2)=(\sqrt{\frac{-a_1\mu}{a_2}},0)$,
$$
\begin{aligned}
U_t(\vartheta)(r,\theta) \approx &\bar{\xi} \mathrm{e}^{-\mathrm{i}\omega_{\hat{\lambda}}\vartheta} J_n(\sqrt{\lambda_{nm}}r)\mathrm{e}^{\mathrm{i}n\theta} \rho_2\mathrm{e}^{-\mathrm{i}\chi_2(t)}
+\xi \mathrm{e}^{\mathrm{i}\omega_{\hat{\lambda}}\vartheta} J_n(\sqrt{\lambda_{nm}}r)\mathrm{e}^{-\mathrm{i}n\theta} \rho_2\mathrm{e}^{\mathrm{i}\chi_2(t)}\\
 \approx &
\sum_{i=1}^n{2|p_{1i}|\sqrt{\frac{-a_1\mu}{a_2}} J_n(\sqrt{\lambda_{nm}}r)\cos(\mathrm{Arg}(p_{1i})+\omega_{\hat{\lambda}}\vartheta+\omega_{\hat{\lambda}}t-n\theta) {e}_i},
\end{aligned}
$$
which corresponds to the form of a rotating wave solution in the opposite direction in the plane of $(z_1,z_4)$.

For $(\rho_1,\rho_2)=(\sqrt{\frac{-a_1\mu}{a_2+a_3}},\sqrt{\frac{-a_1\mu}{a_2+a_3}})$,
$$
\begin{aligned}
U_t(\vartheta)(r,\theta) \approx &
\sum_{i=1}^n{2|p_{1i}|\sqrt{\frac{-a_1\mu}{a_2+a_3}}  J_n(\sqrt{\lambda_{nm}}r)\cos(\mathrm{Arg}(p_{1i})+\omega_{\hat{\lambda}}\vartheta+\omega_{\hat{\lambda}}t+n\theta) {e}_i}\\
+&
\sum_{i=1}^n{2|p_{1i}|\sqrt{\frac{-a_1\mu}{a_2+a_3}}  J_n(\sqrt{\lambda_{nm}}r)\cos(\mathrm{Arg}(p_{1i})+\omega_{\hat{\lambda}}\vartheta+\omega_{\hat{\lambda}}t-n\theta) {e}_i}\\
 \approx &\sum_{i=1}^n{4|p_{1i}|\sqrt{\frac{-a_1\mu}{a_2+a_3}} J_n(\sqrt{\lambda_{nm}}r)\cos(\mathrm{Arg}(p_{1i})+\omega_{\hat{\lambda}}\vartheta+\omega_{\hat{\lambda}}t)\cos(n\theta) {e}_i},
\end{aligned}
$$
which means when $n\theta=\frac{\pi}{2}$ or $n\theta=\frac{3\pi}{2}$, the form of the solution does not change over time. In other words, in a two-dimensional plane, the image of the solution has a fixed axis,
thus, it corresponds to the form of a standing wave solution.

\section{The calculation formula for $A_{p_1p_2p_3p_4}$}\label{A_{p_1p_2p_3p_4}}
\subsection{The calculation formula for $A_{p_1p_2p_3p_4}(p_1+p_2+p_3+p_4=2)$}
$$
\begin{aligned}
A_{2000}=2\hat{\tau}\left(\begin{array}{cccc}
A_{2000}^1\\
A_{2000}^2
\end{array}\right)&,~
A_{1100}=2\hat{\tau}\left(\begin{array}{cccc}
A_{1100}^1\\
A_{1100}^2
\end{array}\right),\\
A_{1010}=2\hat{\tau}\left(\begin{array}{cccc}
A_{1010}^1\\
A_{1010}^2
\end{array}\right)&,~
A_{1001}=2\hat{\tau}\left(\begin{array}{cccc}
A_{1001}^1\\
A_{1001}^2
\end{array}\right),\\
A_{0020}=2\hat{\tau}\left(\begin{array}{cccc}
A_{0020}^1\\
A_{0020}^2
\end{array}\right)&,~
A_{0011}=2\hat{\tau}\left(\begin{array}{cccc}
A_{0011}^1\\
A_{0011}^2
\end{array}\right),\\
A_{0200}=\overline{A_{2000}}&,~A_{0101}=\overline{A_{1010}},\\
A_{0110}=\overline{A_{1001}}&,~A_{0002}=\overline{A_{0020}},
\end{aligned}
$$
with
$$
\begin{aligned}
A_{2000}^1=&F_{20}^{(1)}+F_{11}^{(1)}p_0+F_{02}^{(1)}p_0^2,\\
A_{1100}^1=&2F_{20}^{(1)}+F_{11}^{(1)}(p_0+\bar{p}_0)+2F_{02}^{(1)}p_0\bar{p}_0,\\
A_{1010}^1=&2F_{20}^{(1)}+2F_{11}^{(1)}p_0+2F_{02}^{(1)}p_0^2,\\
A_{1001}^1=&2F_{20}^{(1)}+F_{11}^{(1)}(p_0+\bar{p}_0)+2F_{02}^{(1)}p_0\bar{p}_0,\\
A_{0020}^1=&F_{20}^{(1)}+F_{11}^{(1)}{p}_0+F_{02}^{(1)}p_0^2,\\
A_{0011}^1=&2F_{20}^{(1)}+F_{11}^{(1)}(p_0+\bar{p}_0)+2F_{02}^{(1)}p_0\bar{p}_0,\\
A_{2000}^2=&F_{2000}^{(2)}+F_{1100}^{(2)}p_0+F_{0200}^{(2)}p_0^2+F_{0020}^{(2)}\mathrm{e}^{-2\mathrm{i}\omega_{\hat{\lambda}}\hat{\tau}}
+F_{0011}^{(2)}p_0\mathrm{e}^{-2\mathrm{i}\omega_{\hat{\lambda}}\hat{\tau}}+F_{0002}^{(2)}p_0^2\mathrm{e}^{-2\mathrm{i}\omega_{\hat{\lambda}}\hat{\tau}}\\
&+F_{1010}^{(2)}\mathrm{e}^{-\mathrm{i}\omega_{\hat{\lambda}}\hat{\tau}}+(F_{1001}^{(2)}+F_{0110}^{(2)})p_0\mathrm{e}^{-\mathrm{i}\omega_{\hat{\lambda}}\hat{\tau}}
+F_{0101}^{(2)}p_0^2\mathrm{e}^{-\mathrm{i}\omega_{\hat{\lambda}}\hat{\tau}},\\
A_{1100}^2=&2F_{2000}^{(2)}+F_{1100}^{(2)}(p_0+\bar{p}_0)+2F_{0200}^{(1)}p_0\bar{p}_0+2F_{0020}^{(2)}
+F_{0011}^{(2)}(p_0+\bar{p}_0)+2F_{0002}^{(2)}p_0\bar{p}_0
+F_{1010}^{(2)}(\mathrm{e}^{-\mathrm{i}\omega_{\hat{\lambda}}\hat{\tau}}+\mathrm{e}^{\mathrm{i}\omega_{\hat{\lambda}}\hat{\tau}})\\
&+F_{1001}^{(2)}(p_0\mathrm{e}^{-\mathrm{i}\omega_{\hat{\lambda}}\hat{\tau}}+\bar{p}_0\mathrm{e}^{\mathrm{i}\omega_{\hat{\lambda}}\hat{\tau}})
+F_{0110}^{(2)}(p_0\mathrm{e}^{\mathrm{i}\omega_{\hat{\lambda}}\hat{\tau}}+\bar{p}_0\mathrm{e}^{-\mathrm{i}\omega_{\hat{\lambda}}\hat{\tau}})
+F_{0101}^{(2)}(p_0\bar{p}_0\mathrm{e}^{\mathrm{i}\omega_{\hat{\lambda}}\hat{\tau}}+p_0\bar{p}_0\mathrm{e}^{-\mathrm{i}\omega_{\hat{\lambda}}\hat{\tau}}),\\
A_{1010}^2=&2F_{2000}^{(2)}+2F_{1100}^{(2)}p_0+2F_{0200}^{(1)}p_0^2+2F_{0020}^{(2)}\mathrm{e}^{-2\mathrm{i}\omega_{\hat{\lambda}}\hat{\tau}}
+2F_{0011}^{(2)}p_0\mathrm{e}^{-2\mathrm{i}\omega_{\hat{\lambda}}\hat{\tau}}+2F_{0002}^{(2)}p_0^2\mathrm{e}^{-2\mathrm{i}\omega_{\hat{\lambda}}\hat{\tau}}\\
&+2F_{1010}^{(2)}\mathrm{e}^{-\mathrm{i}\omega_{\hat{\lambda}}\hat{\tau}}+2(F_{1001}^{(2)}+F_{0110}^{(2)})p_0\mathrm{e}^{-\mathrm{i}\omega_{\hat{\lambda}}\hat{\tau}}
+2F_{0101}^{(2)}p_0^2\mathrm{e}^{-\mathrm{i}\omega_{\hat{\lambda}}\hat{\tau}},\\
A_{1001}^2=&2F_{2000}^{(2)}+F_{1100}^{(2)}(p_0+\bar{p}_0)+2F_{0200}^{(1)}p_0\bar{p}_0+2F_{0020}^{(2)}
+F_{0011}^{(2)}(p_0+\bar{p}_0)+2F_{0002}^{(2)}p_0\bar{p}_0
+F_{1010}^{(2)}(\mathrm{e}^{-\mathrm{i}\omega_{\hat{\lambda}}\hat{\tau}}+\mathrm{e}^{\mathrm{i}\omega_{\hat{\lambda}}\hat{\tau}})\\
&+F_{1001}^{(2)}(p_0\mathrm{e}^{-\mathrm{i}\omega_{\hat{\lambda}}\hat{\tau}}+\bar{p}_0\mathrm{e}^{\mathrm{i}\omega_{\hat{\lambda}}\hat{\tau}})
+F_{0110}^{(2)}(p_0\mathrm{e}^{\mathrm{i}\omega_{\hat{\lambda}}\hat{\tau}}+\bar{p}_0\mathrm{e}^{-\mathrm{i}\omega_{\hat{\lambda}}\hat{\tau}})
+F_{0101}^{(2)}(p_0\bar{p}_0\mathrm{e}^{\mathrm{i}\omega_{\hat{\lambda}}\hat{\tau}}+p_0\bar{p}_0\mathrm{e}^{-\mathrm{i}\omega_{\hat{\lambda}}\hat{\tau}}),\\
A_{0020}^2=&F_{2000}^{(2)}+F_{1100}^{(2)}p_0+F_{0200}^{(2)}p_0^2+F_{0020}^{(2)}\mathrm{e}^{-2\mathrm{i}\omega_{\hat{\lambda}}\hat{\tau}}
+F_{0011}^{(2)}p_0\mathrm{e}^{-2\mathrm{i}\omega_{\hat{\lambda}}\hat{\tau}}+F_{0002}^{(2)}p_0^2\mathrm{e}^{-2\mathrm{i}\omega_{\hat{\lambda}}\hat{\tau}}\\
&+F_{1010}^{(2)}\mathrm{e}^{-\mathrm{i}\omega_{\hat{\lambda}}\hat{\tau}}+(F_{1001}^{(2)}+F_{1001}^{(2)})p_0\mathrm{e}^{-\mathrm{i}\omega_{\hat{\lambda}}\hat{\tau}}
+F_{0101}^{(2)}p_0^2\mathrm{e}^{-\mathrm{i}\omega_{\hat{\lambda}}\hat{\tau}},\\
A_{0011}^2=&2F_{2000}^{(2)}+F_{1100}^{(2)}(p_0+\bar{p}_0)+2F_{0200}^{(1)}p_0\bar{p}_0+2F_{0020}^{(2)}
+F_{0011}^{(2)}(p_0+\bar{p}_0)+2F_{0002}^{(2)}p_0\bar{p}_0
+F_{1010}^{(2)}(\mathrm{e}^{-\mathrm{i}\omega_{\hat{\lambda}}\hat{\tau}}+\mathrm{e}^{\mathrm{i}\omega_{\hat{\lambda}}\hat{\tau}})\\
&+F_{1001}^{(2)}(p_0\mathrm{e}^{-\mathrm{i}\omega_{\hat{\lambda}}\hat{\tau}}+\bar{p}_0\mathrm{e}^{\mathrm{i}\omega_{\hat{\lambda}}\hat{\tau}})
+F_{0110}^{(2)}(p_0\mathrm{e}^{\mathrm{i}\omega_{\hat{\lambda}}\hat{\tau}}+\bar{p}_0\mathrm{e}^{-\mathrm{i}\omega_{\hat{\lambda}}\hat{\tau}})
+F_{0101}^{(2)}(p_0\bar{p}_0\mathrm{e}^{\mathrm{i}\omega_{\hat{\lambda}}\hat{\tau}}+p_0\bar{p}_0\mathrm{e}^{-\mathrm{i}\omega_{\hat{\lambda}}\hat{\tau}}).
\end{aligned}
$$
\subsection{The calculation formula for $A_{p_1p_2p_3p_4}(p_1+p_2+p_3+p_4=3)$}
$$
\begin{aligned}
A_{2100}=6\hat{\tau}\left(\begin{array}{cccc}
A_{2100}^1\\
A_{2100}^2
\end{array}\right)&,~
A_{2010}=6\hat{\tau}\left(\begin{array}{cccc}
A_{2010}^1\\
A_{2010}^2
\end{array}\right),\\
A_{2001}=6\hat{\tau}\left(\begin{array}{cccc}
A_{2001}^1\\
A_{2001}^2
\end{array}\right)&,~
A_{1020}=6\hat{\tau}\left(\begin{array}{cccc}
A_{1020}^1\\
A_{1020}^2
\end{array}\right),\\
A_{0120}=6\hat{\tau}\left(\begin{array}{cc}
A_{0120}^1\\
A_{0120}^2
\end{array}\right)&,~
A_{0021}=6\hat{\tau}\left(\begin{array}{cc}
A_{0021}^1\\
A_{0021}^2
\end{array}\right),\\
\end{aligned}
$$
$$
\begin{aligned}
A_{1110}=6\hat{\tau}\left(\begin{array}{cc}
A_{1110}^1\\
A_{1110}^2
\end{array}\right)&,~
A_{1011}=6\hat{\tau}\left(\begin{array}{cc}
A_{1011}^1\\
A_{1011}^2
\end{array}\right),\\
A_{1200}=\overline{A_{2100}}&,~A_{0210}=\overline{A_{2010}},\\
A_{0201}=\overline{A_{2001}}&,~A_{0002}=\overline{A_{0020}},
\end{aligned}
$$
with
$$
\begin{aligned}
A_{2100}^1=&3F_{30}^{(1)}+(\bar{p}_0+2p_0)F_{21}^{(1)}+(p_0^2+2p_0\bar{p}_0)F_{12}^{(1)}+3p_0^2\bar{p}_0F_{03}^{(1)},\\
A_{2010}^1=&3F_{30}^{(1)}+3p_0F_{21}^{(1)}+3p_0^2F_{12}^{(1)}+3p_0^3F_{03}^{(1)},\\
A_{2001}^1=&3F_{30}^{(1)}+(\bar{p}_0+2p_0)F_{21}^{(1)}+(p_0^2+2p_0\bar{p}_0)F_{12}^{(1)}+3p_0^2\bar{p}_0F_{03}^{(1)},\\
A_{1020}^1=&3F_{30}^{(1)}+3p_0F_{21}^{(1)}+3p_0^2F_{12}^{(1)}+3p_0^3F_{03}^{(1)},\\
A_{0120}^1=&3F_{30}^{(1)}+(\bar{p}_0+2p_0)F_{21}^{(1)}+(p_0^2+2p_0\bar{p}_0)F_{12}^{(1)}+3p_0^2\bar{p}_0F_{03}^{(1)},\\
A_{0021}^1=&3F_{30}^{(1)}+(\bar{p}_0+2p_0)F_{21}^{(1)}+(p_0^2+2p_0\bar{p}_0)F_{12}^{(1)}+3p_0^2\bar{p}_0F_{03}^{(1)},\\
A_{1110}^1=&6F_{30}^{(1)}+(2\bar{p}_0+4p_0)F_{21}^{(1)}+(2p_0^2+4p_0\bar{p}_0)F_{12}^{(1)}+6p_0^2\bar{p}_0F_{03}^{(1)},\\
A_{1011}^1=&6F_{30}^{(1)}+(2\bar{p}_0+4p_0)F_{21}^{(1)}+(2p_0^2+4p_0\bar{p}_0)F_{12}^{(1)}+6p_0^2\bar{p}_0F_{03}^{(1)},\\
A_{2100}^2=&3F_{3000}^{(2)}+(\bar{p}_0+2p_0)F_{2100}^{(2)}+(p_0^2+2p_0\bar{p}_0)F_{1200}^{(2)}+3p_0^2\bar{p}_0F_{0300}^{(2)}\\
&+\left(3F_{0030}^{(2)}+(\bar{p}_0+2p_0)F_{0021}^{(2)}+(p_0^2+2p_0\bar{p}_0)F_{0012}^{(2)}+3p_0^2\bar{p}_0F_{0003}^{(2)}\right)\mathrm{e}^{-\mathrm{i}\omega_{\hat{\lambda}}\hat{\tau}}\\
&+F_{2010}^{(2)}(\mathrm{e}^{\mathrm{i}\omega_{\hat{\lambda}}\hat{\tau}}+2\mathrm{e}^{-\mathrm{i}\omega_{\hat{\lambda}}\hat{\tau}})
+F_{2001}^{(2)}(\bar{p}_0\mathrm{e}^{\mathrm{i}\omega_{\hat{\lambda}}\hat{\tau}}+2p_0\mathrm{e}^{-\mathrm{i}\omega_{\hat{\lambda}}\hat{\tau}})
+F_{0210}^{(2)}(2p_0\bar{p}_0\mathrm{e}^{-\mathrm{i}\omega_{\hat{\lambda}}\hat{\tau}}+p_0^2\mathrm{e}^{\mathrm{i}\omega_{\hat{\lambda}}\hat{\tau}})\\
&+F_{0201}^{(2)}(2p_0^2\bar{p}_0\mathrm{e}^{-\mathrm{i}\omega_{\hat{\lambda}}\hat{\tau}}+p_0^2\bar{p}_0\mathrm{e}^{\mathrm{i}\omega_{\hat{\lambda}}\hat{\tau}})
+F_{1020}^{(2)}(\mathrm{e}^{-2\mathrm{i}\omega_{\hat{\lambda}}\hat{\tau}}+2)
+F_{1002}^{(2)}(p_0^2\mathrm{e}^{-2\mathrm{i}\omega_{\hat{\lambda}}\hat{\tau}}+2p_0\bar{p}_0)\\
&+F_{0120}^{(2)}(\bar{p}_0\mathrm{e}^{-2\mathrm{i}\omega_{\hat{\lambda}}\hat{\tau}}+2p_0)
+F_{0102}^{(2)}(p_0^2\bar{p}_0\mathrm{e}^{-2\mathrm{i}\omega_{\hat{\lambda}}\hat{\tau}}+2p_0^2\bar{p}_0),\\
A_{2010}^2=&3F_{3000}^{(2)}+3p_0F_{2100}^{(2)}+3p_0^2F_{1200}^{(2)}+3p_0^3F_{0300}^{(2)}\\
&+\left(3F_{0030}^{(2)}+3p_0F_{0021}^{(2)}+3p_0^2F_{0012}^{(2)}+3p_0^3F_{0003}^{(2)}\right)\mathrm{e}^{-3\mathrm{i}\omega_{\hat{\lambda}}\hat{\tau}}\\
&+\left(3F_{2010}^{(2)}+3p_0F_{2001}^{(2)}+3p_0^2F_{0210}^{(2)}
+3p_0^3F_{0201}^{(2)}\right)\mathrm{e}^{-\mathrm{i}\omega_{\hat{\lambda}}\hat{\tau}}\\
&+\left(3F_{1020}^{(2)}+3p_0F_{0120}^{(2)}+3p_0F_{1002}^{(2)}+
+3p_0^3F_{0102}^{(2)}\right)\mathrm{e}^{-2\mathrm{i}\omega_{\hat{\lambda}}\hat{\tau}},\\
A_{2001}^2=&3F_{3000}^{(2)}+(\bar{p}_0+2p_0)F_{2100}^{(2)}+(p_0^2+2p_0\bar{p}_0)F_{1200}^{(2)}+3p_0^2\bar{p}_0F_{0300}^{(2)}\\
&+\left(3F_{0030}^{(2)}+(\bar{p}_0+2p_0)F_{0021}^{(2)}+(p_0^2+2p_0\bar{p}_0)F_{0012}^{(2)}+3p_0^2\bar{p}_0F_{0003}^{(2)}\right)\mathrm{e}^{-\mathrm{i}\omega_{\hat{\lambda}}\hat{\tau}}\\
&+F_{2010}^{(2)}(\mathrm{e}^{\mathrm{i}\omega_{\hat{\lambda}}\hat{\tau}}+2\mathrm{e}^{-\mathrm{i}\omega_{\hat{\lambda}}\hat{\tau}})
+F_{2001}^{(2)}(\bar{p}_0\mathrm{e}^{\mathrm{i}\omega_{\hat{\lambda}}\hat{\tau}}+2p_0\mathrm{e}^{-\mathrm{i}\omega_{\hat{\lambda}}\hat{\tau}})
+F_{0210}^{(2)}(2p_0\bar{p}_0\mathrm{e}^{-\mathrm{i}\omega_{\hat{\lambda}}\hat{\tau}}+p_0^2\mathrm{e}^{\mathrm{i}\omega_{\hat{\lambda}}\hat{\tau}})\\
&+F_{0201}^{(2)}(2p_0^2\bar{p}_0\mathrm{e}^{-\mathrm{i}\omega_{\hat{\lambda}}\hat{\tau}}+p_0^2\bar{p}_0\mathrm{e}^{\mathrm{i}\omega_{\hat{\lambda}}\hat{\tau}})
+F_{1020}^{(2)}(\mathrm{e}^{-2\mathrm{i}\omega_{\hat{\lambda}}\hat{\tau}}+2)
+F_{1002}^{(2)}(p_0^2\mathrm{e}^{-2\mathrm{i}\omega_{\hat{\lambda}}\hat{\tau}}+2p_0\bar{p}_0)\\
&+F_{0120}^{(2)}(\bar{p}_0\mathrm{e}^{-2\mathrm{i}\omega_{\hat{\lambda}}\hat{\tau}}+2p_0)
+F_{0102}^{(2)}(p_0^2\bar{p}_0\mathrm{e}^{-2\mathrm{i}\omega_{\hat{\lambda}}\hat{\tau}}+2p_0^2\bar{p}_0),\\
\end{aligned}
$$
$$
\begin{aligned}
A_{1020}^2=&3F_{3000}^{(2)}+3p_0F_{2100}^{(2)}+3p_0^2F_{1200}^{(2)}+3p_0^3F_{0300}^{(2)}\\
&+\left(3F_{0030}^{(2)}+3p_0F_{0021}^{(2)}+3p_0^2F_{0012}^{(2)}+3p_0^3F_{0003}^{(2)}\right)\mathrm{e}^{-3\mathrm{i}\omega_{\hat{\lambda}}\hat{\tau}}\\
&+\left(3F_{2010}^{(2)}+3p_0F_{2001}^{(2)}+3p_0^2F_{0210}^{(2)}
+3p_0^3F_{0201}^{(2)}\right)\mathrm{e}^{-\mathrm{i}\omega_{\hat{\lambda}}\hat{\tau}}\\
&+\left(3F_{1020}^{(2)}+3p_0F_{0120}^{(2)}+3p_0F_{1002}^{(2)}+
+3p_0^3F_{0102}^{(2)}\right)\mathrm{e}^{-2\mathrm{i}\omega_{\hat{\lambda}}\hat{\tau}},\\
A_{0120}^2=&3F_{3000}^{(2)}+(\bar{p}_0+2p_0)F_{2100}^{(2)}+(p_0^2+2p_0\bar{p}_0)F_{1200}^{(2)}+3p_0^2\bar{p}_0F_{0300}^{(2)}\\
&+\left(3F_{0030}^{(2)}+(\bar{p}_0+2p_0)F_{0021}^{(2)}+(p_0^2+2p_0\bar{p}_0)F_{0012}^{(2)}+3p_0^2\bar{p}_0F_{0003}^{(2)}\right)\mathrm{e}^{-\mathrm{i}\omega_{\hat{\lambda}}\hat{\tau}}\\
&+F_{2010}^{(2)}(\mathrm{e}^{\mathrm{i}\omega_{\hat{\lambda}}\hat{\tau}}+2\mathrm{e}^{-\mathrm{i}\omega_{\hat{\lambda}}\hat{\tau}})
+F_{2001}^{(2)}(\bar{p}_0\mathrm{e}^{\mathrm{i}\omega_{\hat{\lambda}}\hat{\tau}}+2p_0\mathrm{e}^{-\mathrm{i}\omega_{\hat{\lambda}}\hat{\tau}})
+F_{0210}^{(2)}(2p_0\bar{p}_0\mathrm{e}^{-\mathrm{i}\omega_{\hat{\lambda}}\hat{\tau}}+p_0^2\mathrm{e}^{\mathrm{i}\omega_{\hat{\lambda}}\hat{\tau}})\\
&+F_{0201}^{(2)}(2p_0^2\bar{p}_0\mathrm{e}^{-\mathrm{i}\omega_{\hat{\lambda}}\hat{\tau}}+p_0^2\bar{p}_0\mathrm{e}^{\mathrm{i}\omega_{\hat{\lambda}}\hat{\tau}})
+F_{1020}^{(2)}(\mathrm{e}^{-2\mathrm{i}\omega_{\hat{\lambda}}\hat{\tau}}+2)
+F_{1002}^{(2)}(p_0^2\mathrm{e}^{-2\mathrm{i}\omega_{\hat{\lambda}}\hat{\tau}}+2p_0\bar{p}_0)\\
&+F_{0120}^{(2)}(\bar{p}_0\mathrm{e}^{-2\mathrm{i}\omega_{\hat{\lambda}}\hat{\tau}}+2p_0)
+F_{0102}^{(2)}(p_0^2\bar{p}_0\mathrm{e}^{-2\mathrm{i}\omega_{\hat{\lambda}}\hat{\tau}}+2p_0^2\bar{p}_0),\\
A_{0021}^2=&3F_{3000}^{(2)}+(\bar{p}_0+2p_0)F_{2100}^{(2)}+(p_0^2+2p_0\bar{p}_0)F_{1200}^{(2)}+3p_0^2\bar{p}_0F_{0300}^{(2)}\\
&+\left(3F_{0030}^{(2)}+(\bar{p}_0+2p_0)F_{0021}^{(2)}+(p_0^2+2p_0\bar{p}_0)F_{0012}^{(2)}+3p_0^2\bar{p}_0F_{0003}^{(2)}\right)\mathrm{e}^{-\mathrm{i}\omega_{\hat{\lambda}}\hat{\tau}}\\
&+F_{2010}^{(2)}(\mathrm{e}^{\mathrm{i}\omega_{\hat{\lambda}}\hat{\tau}}+2\mathrm{e}^{-\mathrm{i}\omega_{\hat{\lambda}}\hat{\tau}})
+F_{2001}^{(2)}(\bar{p}_0\mathrm{e}^{\mathrm{i}\omega_{\hat{\lambda}}\hat{\tau}}+2p_0\mathrm{e}^{-\mathrm{i}\omega_{\hat{\lambda}}\hat{\tau}})
+F_{0210}^{(2)}(2p_0\bar{p}_0\mathrm{e}^{-\mathrm{i}\omega_{\hat{\lambda}}\hat{\tau}}+p_0^2\mathrm{e}^{\mathrm{i}\omega_{\hat{\lambda}}\hat{\tau}})\\
&+F_{0201}^{(2)}(2p_0^2\bar{p}_0\mathrm{e}^{-\mathrm{i}\omega_{\hat{\lambda}}\hat{\tau}}+p_0^2\bar{p}_0\mathrm{e}^{\mathrm{i}\omega_{\hat{\lambda}}\hat{\tau}})
+F_{1020}^{(2)}(\mathrm{e}^{-2\mathrm{i}\omega_{\hat{\lambda}}\hat{\tau}}+2)
+F_{1002}^{(2)}(p_0^2\mathrm{e}^{-2\mathrm{i}\omega_{\hat{\lambda}}\hat{\tau}}+2p_0\bar{p}_0)\\
&+F_{0120}^{(2)}(\bar{p}_0\mathrm{e}^{-2\mathrm{i}\omega_{\hat{\lambda}}\hat{\tau}}+2p_0)
+F_{0102}^{(2)}(p_0^2\bar{p}_0\mathrm{e}^{-2\mathrm{i}\omega_{\hat{\lambda}}\hat{\tau}}+2p_0^2\bar{p}_0),\\
A_{1110}^2=&6F_{3000}^{(2)}+(2\bar{p}_0+4p_0)F_{2100}^{(2)}+(2p_0^2+4p_0\bar{p}_0)F_{1200}^{(2)}+6p_0^2\bar{p}_0F_{0300}^{(2)}\\
&+\left(6F_{0030}^{(2)}+(2\bar{p}_0+4p_0)F_{0021}^{(2)}+(2p_0^2+4p_0\bar{p}_0)F_{0012}^{(2)}+6p_0^2\bar{p}_0F_{0003}^{(2)}\right)\mathrm{e}^{-\mathrm{i}\omega_{\hat{\lambda}}\hat{\tau}}\\
&+F_{2010}^{(2)}(2\mathrm{e}^{\mathrm{i}\omega_{\hat{\lambda}}\hat{\tau}}+4\mathrm{e}^{-\mathrm{i}\omega_{\hat{\lambda}}\hat{\tau}})
+F_{2001}^{(2)}(2\bar{p}_0\mathrm{e}^{\mathrm{i}\omega_{\hat{\lambda}}\hat{\tau}}+4p_0\mathrm{e}^{-\mathrm{i}\omega_{\hat{\lambda}}\hat{\tau}})
+F_{0210}^{(2)}(4p_0\bar{p}_0\mathrm{e}^{-\mathrm{i}\omega_{\hat{\lambda}}\hat{\tau}}+2p_0^2\mathrm{e}^{\mathrm{i}\omega_{\hat{\lambda}}\hat{\tau}})\\
&+F_{0201}^{(2)}(4p_0^2\bar{p}_0\mathrm{e}^{-\mathrm{i}\omega_{\hat{\lambda}}\hat{\tau}}+2p_0^2\bar{p}_0\mathrm{e}^{\mathrm{i}\omega_{\hat{\lambda}}\hat{\tau}})
+F_{1020}^{(2)}(2\mathrm{e}^{-2\mathrm{i}\omega_{\hat{\lambda}}\hat{\tau}}+4)
+F_{1002}^{(2)}(2p_0^2\mathrm{e}^{-2\mathrm{i}\omega_{\hat{\lambda}}\hat{\tau}}+4p_0\bar{p}_0)\\
&+F_{0120}^{(2)}(2\bar{p}_0\mathrm{e}^{-2\mathrm{i}\omega_{\hat{\lambda}}\hat{\tau}}+4p_0)
+F_{0102}^{(2)}(2p_0^2\bar{p}_0\mathrm{e}^{-2\mathrm{i}\omega_{\hat{\lambda}}\hat{\tau}}+4p_0^2\bar{p}_0),\\
A_{1011}^2=&6F_{3000}^{(2)}+(2\bar{p}_0+4p_0)F_{2100}^{(2)}+(2p_0^2+4p_0\bar{p}_0)F_{1200}^{(2)}+6p_0^2\bar{p}_0F_{0300}^{(2)}\\
&+\left(6F_{0030}^{(2)}+(2\bar{p}_0+4p_0)F_{0021}^{(2)}+(2p_0^2+4p_0\bar{p}_0)F_{0012}^{(2)}+6p_0^2\bar{p}_0F_{0003}^{(2)}\right)\mathrm{e}^{-\mathrm{i}\omega_{\hat{\lambda}}\hat{\tau}}\\
&+F_{2010}^{(2)}(2\mathrm{e}^{\mathrm{i}\omega_{\hat{\lambda}}\hat{\tau}}+4\mathrm{e}^{-\mathrm{i}\omega_{\hat{\lambda}}\hat{\tau}})
+F_{2001}^{(2)}(2\bar{p}_0\mathrm{e}^{\mathrm{i}\omega_{\hat{\lambda}}\hat{\tau}}+4p_0\mathrm{e}^{-\mathrm{i}\omega_{\hat{\lambda}}\hat{\tau}})
+F_{0210}^{(2)}(4p_0\bar{p}_0\mathrm{e}^{-\mathrm{i}\omega_{\hat{\lambda}}\hat{\tau}}+2p_0^2\mathrm{e}^{\mathrm{i}\omega_{\hat{\lambda}}\hat{\tau}})\\
&+F_{0201}^{(2)}(4p_0^2\bar{p}_0\mathrm{e}^{-\mathrm{i}\omega_{\hat{\lambda}}\hat{\tau}}+2p_0^2\bar{p}_0\mathrm{e}^{\mathrm{i}\omega_{\hat{\lambda}}\hat{\tau}})
+F_{1020}^{(2)}(2\mathrm{e}^{-2\mathrm{i}\omega_{\hat{\lambda}}\hat{\tau}}+4)
+F_{1002}^{(2)}(2p_0^2\mathrm{e}^{-2\mathrm{i}\omega_{\hat{\lambda}}\hat{\tau}}+4p_0\bar{p}_0)\\
&+F_{0120}^{(2)}(2\bar{p}_0\mathrm{e}^{-2\mathrm{i}\omega_{\hat{\lambda}}\hat{\tau}}+4p_0)
+F_{0102}^{(2)}(2p_0^2\bar{p}_0\mathrm{e}^{-2\mathrm{i}\omega_{\hat{\lambda}}\hat{\tau}}+4p_0^2\bar{p}_0).
\end{aligned}
$$

\section{The calculation formula for $S_{y(0)z_k},~S_{y(-1)z_k},~k=1,2,3,4$}\label{S}
$$
\begin{aligned}
&S_{y(0)z_1}=\left(\begin{array}{cc}
2 F_{20}^{(1)}+F_{11}^{(1)}p_0 & F_{11}^{(1)}+2 F_{02}^{(1)}p_0\\
2 F_{2000}^{(2)}+F_{1100}^{(2)}p_0+F_{1010}^{(2)}\mathrm{e}^{-\mathrm{i}\omega_{\hat{\lambda}}\hat{\tau}}+F_{1001}^{(2)}p_0\mathrm{e}^{-\mathrm{i}\omega_{\hat{\lambda}}\hat{\tau}} & F_{1100}^{(2)}+2 F_{0200}^{(2)}p_0+F_{0110}^{(2)}\mathrm{e}^{-\mathrm{i}\omega_{\hat{\lambda}}\hat{\tau}}+F_{0101}^{(2)}p_0\mathrm{e}^{-\mathrm{i}\omega_{\hat{\lambda}}\hat{\tau}}\\
\end{array}\right),\\
&S_{y(0)z_2}=\left(\begin{array}{cccc}
2 F_{20}^{(1)}+F_{11}^{(1)}\bar{p}_0 & F_{11}^{(1)}+2 F_{02}^{(1)}\bar{p}_0\\
2 F_{2000}^{(2)}+F_{1100}^{(2)}\bar{p}_0+F_{1010}^{(2)}\mathrm{e}^{\mathrm{i}\omega_{\hat{\lambda}}\hat{\tau}}+F_{1001}^{(2)}\bar{p}_0\mathrm{e}^{\mathrm{i}\omega_{\hat{\lambda}}\hat{\tau}} & F_{1100}^{(2)}+2 F_{0200}^{(2)}\bar{p}_0+F_{0110}^{(2)}\mathrm{e}^{\mathrm{i}\omega_{\hat{\lambda}}\hat{\tau}}+F_{0101}^{(2)}\bar{p}_0\mathrm{e}^{\mathrm{i}\omega_{\hat{\lambda}}\hat{\tau}}\\
\end{array}\right),\\
&S_{y(0)z_3}=\left(\begin{array}{cccc}
2 F_{20}^{(1)}+F_{11}^{(1)}p_0 & F_{11}^{(1)}+2 F_{02}^{(1)}p_0\\
2 F_{2000}^{(2)}+F_{1100}^{(2)}p_0+F_{1010}^{(2)}\mathrm{e}^{-\mathrm{i}\omega_{\hat{\lambda}}\hat{\tau}}+F_{1001}^{(2)}p_0\mathrm{e}^{-\mathrm{i}\omega_{\hat{\lambda}}\hat{\tau}} & F_{1100}^{(2)}+2 F_{0200}^{(2)}p_0+F_{0110}^{(2)}\mathrm{e}^{-\mathrm{i}\omega_{\hat{\lambda}}\hat{\tau}}+F_{0101}^{(2)}p_0\mathrm{e}^{-\mathrm{i}\omega_{\hat{\lambda}}\hat{\tau}}\\
\end{array}\right),\\
&S_{y(0)z_4}=\left(\begin{array}{cccc}
2 F_{20}^{(1)}+F_{11}^{(1)}\bar{p}_0 & F_{11}^{(1)}+2 F_{02}^{(1)}\bar{p}_0\\
2 F_{2000}^{(2)}+F_{1100}^{(2)}\bar{p}_0+F_{1010}^{(2)}\mathrm{e}^{\mathrm{i}\omega_{\hat{\lambda}}\hat{\tau}}+F_{1001}^{(2)}\bar{p}_0\mathrm{e}^{\mathrm{i}\omega_{\hat{\lambda}}\hat{\tau}} & F_{1100}^{(2)}+2 F_{0200}^{(2)}\bar{p}_0+F_{0110}^{(2)}\mathrm{e}^{\mathrm{i}\omega_{\hat{\lambda}}\hat{\tau}}+F_{0101}^{(2)}\bar{p}_0\mathrm{e}^{\mathrm{i}\omega_{\hat{\lambda}}\hat{\tau}}\\
\end{array}\right),\\
\end{aligned}
$$
$$
\begin{aligned}
&S_{y(-1)z_1}=\left(\begin{array}{cccc}
0 & 0\\
2 F_{0020}^{(2)}\mathrm{e}^{-\mathrm{i}\omega_{\hat{\lambda}}\hat{\tau}}+F_{0011}^{(2)}p_0\mathrm{e}^{-\mathrm{i}\omega_{\hat{\lambda}}\hat{\tau}}+F_{1010}^{(2)}+F_{0110}^{(2)}p_0 & F_{0011}^{(2)}\mathrm{e}^{-\mathrm{i}\omega_{\hat{\lambda}}\hat{\tau}}+2 F_{0002}^{(2)}p_0\mathrm{e}^{-\mathrm{i}\omega_{\hat{\lambda}}\hat{\tau}}+F_{1001}^{(2)}+F_{0101}^{(2)}p_0\\
\end{array}\right),\\
&S_{y(-1)z_2}=\left(\begin{array}{cccc}
0 & 0\\
2 F_{0020}^{(2)}\mathrm{e}^{\mathrm{i}\omega_{\hat{\lambda}}\hat{\tau}}+F_{0011}^{(2)}\bar{p}_0\mathrm{e}^{\mathrm{i}\omega_{\hat{\lambda}}\hat{\tau}}+F_{1010}^{(2)}+F_{0110}^{(2)}\bar{p}_0 & F_{0011}^{(2)}\mathrm{e}^{\mathrm{i}\omega_{\hat{\lambda}}\hat{\tau}}+2 F_{0002}^{(2)}\bar{p}_0\mathrm{e}^{\mathrm{i}\omega_{\hat{\lambda}}\hat{\tau}}+F_{1001}^{(2)}+F_{0101}^{(2)}\bar{p}_0\\
\end{array}\right),\\
&S_{y(-1)z_3}=\left(\begin{array}{cccc}
0 & 0\\
2 F_{0020}^{(2)}\mathrm{e}^{-\mathrm{i}\omega_{\hat{\lambda}}\hat{\tau}}+F_{0011}^{(2)}p_0\mathrm{e}^{-\mathrm{i}\omega_{\hat{\lambda}}\hat{\tau}}+F_{1010}^{(2)}+F_{0110}^{(2)}p_0 & F_{0011}^{(2)}\mathrm{e}^{-\mathrm{i}\omega_{\hat{\lambda}}\hat{\tau}}+2 F_{0002}^{(2)}p_0\mathrm{e}^{-\mathrm{i}\omega_{\hat{\lambda}}\hat{\tau}}+F_{1001}^{(2)}+F_{0101}^{(2)}p_0\\
\end{array}\right),\\
&S_{y(-1)z_4}=\left(\begin{array}{cccc}
0 & 0\\
2 F_{0020}^{(2)}\mathrm{e}^{\mathrm{i}\omega_{\hat{\lambda}}\hat{\tau}}+F_{0011}^{(2)}\bar{p}_0\mathrm{e}^{\mathrm{i}\omega_{\hat{\lambda}}\hat{\tau}}+F_{1010}^{(2)}+F_{0110}^{(2)}\bar{p}_0 & F_{0011}^{(2)}\mathrm{e}^{\mathrm{i}\omega_{\hat{\lambda}}\hat{\tau}}+2 F_{0002}^{(2)}\bar{p}_0\mathrm{e}^{\mathrm{i}\omega_{\hat{\lambda}}\hat{\tau}}+F_{1001}^{(2)}+F_{0101}^{(2)}\bar{p}_0\\
\end{array}\right).
\end{aligned}
$$

\section{The calculation formula for $h_{jp_1p_2p_3p_4}$}\label{h}
$$
\begin{aligned}
&h_{0k2000}^{ccs}(\vartheta)=-\mathrm{M}_{0kcs}^c\mathrm{e}^{2\mathrm{i}\omega_{\hat{\lambda}}\vartheta}\left[-2\mathrm{i}\omega_{\hat{\lambda}}-\lambda_{0k}\tilde{D}_0+\tilde{L}_0(\mathrm{e}^{2\mathrm{i}\omega_{\hat{\lambda}}}\cdot I_d)\right]^{-1}A_{2000},\\
&h_{0k1100}^{ccs}(\vartheta)=-\mathrm{M}_{0kcs}^c\left[-\lambda_{0k}\tilde{D}_0+\tilde{L}_0(I_d)\right]^{-1}A_{1100},\\
&h_{0k1010}^{ccs}(\vartheta)=-\mathrm{M}_{0kcs}^c\mathrm{e}^{2\mathrm{i}\omega_{\hat{\lambda}}\vartheta}\left[-2\mathrm{i}\omega_{\hat{\lambda}}-\lambda_{0k}\tilde{D}_0+\tilde{L}_0(\mathrm{e}^{2\mathrm{i}\omega_{\hat{\lambda}}}\cdot I_d)\right]^{-1}A_{1010},\\
&h_{0k1001}^{ccs}(\vartheta)=-\mathrm{M}_{0kcs}^c\left[-\lambda_{0k}\tilde{D}_0+\tilde{L}_0(I_d)\right]^{-1}A_{1001},\\
&h_{0k0110}^{ccs}(\vartheta)=-\mathrm{M}_{0kcs}^c\left[-\lambda_{0k}\tilde{D}_0+\tilde{L}_0(I_d)\right]^{-1}A_{0110},\\
&h_{0k0020}^{ccs}(\vartheta)=-\mathrm{M}_{0kcs}^c\mathrm{e}^{2\mathrm{i}\omega_{\hat{\lambda}}\vartheta}\left[-2\mathrm{i}\omega_{\hat{\lambda}}-\lambda_{0k}\tilde{D}_0+\tilde{L}_0(\mathrm{e}^{2\mathrm{i}\omega_{\hat{\lambda}}}\cdot I_d)\right]^{-1}A_{0020},\\
&h_{0k0011}^{ccs}(\vartheta)=-\mathrm{M}_{0kcs}^c\left[-\lambda_{0k}\tilde{D}_0+\tilde{L}_0(I_d)\right]^{-1}A_{0011},\\
\end{aligned}
$$
where $k=0,1,2 \cdots,$

$$
\begin{aligned}
&h_{2nk2000}^{ccs}(\vartheta)=-\mathrm{M}_{2nkcc}^s\mathrm{e}^{2\mathrm{i}\omega_{\hat{\lambda}}\vartheta}\left[-2\mathrm{i}\omega_{\hat{\lambda}}-\lambda_{2nk}\tilde{D}_0+\tilde{L}_0(\mathrm{e}^{2\mathrm{i}\omega_{\hat{\lambda}}}\cdot I_d)\right]^{-1}A_{2000},\\
&h_{2nk1100}^{ccs}(\vartheta)=-\mathrm{M}_{2nkcc}^s\left[-\lambda_{2nk}\tilde{D}_0+\tilde{L}_0(I_d)\right]^{-1}A_{1100},\\
&h_{2nk1010}^{ccs}(\vartheta)=-\mathrm{M}_{2nkcc}^s\mathrm{e}^{2\mathrm{i}\omega_{\hat{\lambda}}\vartheta}\left[-2\mathrm{i}\omega_{\hat{\lambda}}-\lambda_{2nk}\tilde{D}_0+\tilde{L}_0(\mathrm{e}^{2\mathrm{i}\omega_{\hat{\lambda}}}\cdot I_d)\right]^{-1}A_{1010},\\
&h_{2nk1001}^{ccs}(\vartheta)=-\mathrm{M}_{2nkcc}^s\left[-\lambda_{2nk}\tilde{D}_0+\tilde{L}_0(I_d)\right]^{-1}A_{1001},\\
&h_{2nk0110}^{ccs}(\vartheta)=-\mathrm{M}_{2nkcc}^s\left[-\lambda_{2nk}\tilde{D}_0+\tilde{L}_0(I_d)\right]^{-1}A_{0110},\\
&h_{2nk0020}^{ccs}(\vartheta)=-\mathrm{M}_{2nkcc}^s\mathrm{e}^{2\mathrm{i}\omega_{\hat{\lambda}}\vartheta}\left[-2\mathrm{i}\omega_{\hat{\lambda}}-\lambda_{2nk}\tilde{D}_0+\tilde{L}_0(\mathrm{e}^{2\mathrm{i}\omega_{\hat{\lambda}}}\cdot I_d)\right]^{-1}A_{0020},\\
&h_{2nk0011}^{ccs}(\vartheta)=-\mathrm{M}_{2nkcc}^s\left[-\lambda_{2nk}\tilde{D}_0+\tilde{L}_0(I_d)\right]^{-1}A_{0011},\\
&h_{2nk2000}^{css}(\vartheta)=-\mathrm{M}_{2nkss}^c\mathrm{e}^{2\mathrm{i}\omega_{\hat{\lambda}}\vartheta}\left[-2\mathrm{i}\omega_{\hat{\lambda}}-\lambda_{2nk}\tilde{D}_0+\tilde{L}_0(\mathrm{e}^{2\mathrm{i}\omega_{\hat{\lambda}}}\cdot I_d)\right]^{-1}A_{2000},\\
&h_{2nk1100}^{css}(\vartheta)=-\mathrm{M}_{2nkss}^c\left[-\lambda_{2nk}\tilde{D}_0+\tilde{L}_0(I_d)\right]^{-1}A_{1100},\\
&h_{2nk1010}^{css}(\vartheta)=-\mathrm{M}_{2nkss}^c\mathrm{e}^{2\mathrm{i}\omega_{\hat{\lambda}}\vartheta}\left[-2\mathrm{i}\omega_{\hat{\lambda}}-\lambda_{2nk}\tilde{D}_0+\tilde{L}_0(\mathrm{e}^{2\mathrm{i}\omega_{\hat{\lambda}}}\cdot I_d)\right]^{-1}A_{1010},\\
&h_{2nk1001}^{css}(\vartheta)=-\mathrm{M}_{2nkss}^c\left[-\lambda_{2nk}\tilde{D}_0+\tilde{L}_0(I_d)\right]^{-1}A_{1001},\\
&h_{2nk0110}^{css}(\vartheta)=-\mathrm{M}_{2nkss}^c\left[-\lambda_{2nk}\tilde{D}_0+\tilde{L}_0(I_d)\right]^{-1}A_{0110},\\
&h_{2nk0020}^{css}(\vartheta)=-\mathrm{M}_{2nkss}^c\mathrm{e}^{2\mathrm{i}\omega_{\hat{\lambda}}\vartheta}\left[-2\mathrm{i}\omega_{\hat{\lambda}}-\lambda_{2nk}\tilde{D}_0+\tilde{L}_0(\mathrm{e}^{2\mathrm{i}\omega_{\hat{\lambda}}}\cdot I_d)\right]^{-1}A_{0020},\\
&h_{2nk0011}^{css}(\vartheta)=-\mathrm{M}_{2nkss}^c\left[-\lambda_{2nk}\tilde{D}_0+\tilde{L}_0(I_d)\right]^{-1}A_{0011}.\\
\end{aligned}
$$
where $k=1,2,\cdots$.

\end{document}